\documentclass[aop,MSNbibl,secthm,citesort,dvips]{arximspdf}
\usepackage{graphicx}

\doi{10.1214/10-AOP605}
\volume{39}
\issue{6}
\pubyear{2011}
\firstpage{2119}
\lastpage{2177}

\makeatletter

\newtheorem{lemma}[thm]{Lemma}
\newtheorem{corollary}[thm]{Corollary}
\newtheorem{conjecture}[thm]{Conjecture}
\newtheorem {prop}[thm]{Proposition}

\newproclaim{remarks}{Remarks}
\newproclaim{remark}{Remark}

\newcommand{\midd}{\mid}

\def\l{\lambda}
\def\eps{\varepsilon}
\def\E{{\mathbb{E}}}
\def\P{{\mathbb{P}}}
\def\R{{\mathbb{R}}}
\def\Z{{\mathbb{Z}}}
\def\N{{\mathbb{N}}}
\def\re{{e}}
\def\F{{\mathcal{F}}}
\def\ud{{d}}
\def\1{{\mathbf{1}}}
\def\tX{{\tilde X}}
\def\tY{{\tilde Y}}
\def\tZ{{\tilde Z}}
\def\tD{{\tilde\Delta}}
\def\Var{\operatorname{\mathbb{V}\mathrm{ar}}}

\newcommand{\bfrac}[2]{#1/(#2)}
\newcommand{\afrac}[2]{(#1)/(#2)}
\newcommand{\gfrac}[2]{#1/#2}

\newcommand{\sgn}{\operatorname{sgn}}

\newcommand{\Polya}{P\'olya\ }

\makeatother

\begin{document}
\begin{frontmatter}

\title{The simple harmonic urn}
\runtitle{The simple harmonic urn}

\begin{aug}
\author[A]{\fnms{Edward} \snm{Crane}\ead[label=e1]{edward.crane@gmail.com}},
\author[A]{\fnms{Nicholas} \snm{Georgiou}\ead[label=e2]{n.georgiou@gmail.com}},
\author[B]{\fnms{Stanislav} \snm{Volkov}\corref{}\ead[label=e3]{s.volkov@bristol.ac.uk}},
\author[C]{\fnms{Andrew~R.} \snm{Wade}\ead[label=e4]{andrew.wade@strath.ac.uk}}
and
\author[A]{\fnms{Robert J.} \snm{Waters}\ead[label=e5]{rob@aquae.org.uk}}
\runauthor{E. Crane et al.}
\affiliation{University of Bristol, University of Bristol, University
of Bristol, University~of~Strathclyde and University of Bristol}
\address[A]{E. Crane\\
N. Georgiou\\
R. J. Waters\\
Heilbronn Institute for Mathematical Research\\
Department of Mathematics\\
University of Bristol\\
Bristol BS8~1TW\\
United Kingdom\\
\printead{e1}\\
\hphantom{E-mail: }\printead*{e2}\\
\hphantom{E-mail: }\printead*{e5}} %adresu isvedimo komanda gale!
\address[B]{S. Volkov\\
Department of Mathematics\\
University of Bristol\\
Bristol BS8~1TW\\
United Kingdom\\
\printead{e3}\\}
\address[C]{A. R. Wade\\
Department of Mathematics and Statistics\\
University of Strathclyde\\
Glasgow G1 1XH\\
United Kingdom\\
\printead{e4}}
\end{aug}

% HISTORY:
\received{\smonth{11} \syear{2009}}
\revised{\smonth{8} \syear{2010}}

% ABSTRACT
%
\begin{abstract}
We study a generalized P\'olya urn model with two types of ball.
If the drawn ball is red, it is replaced together with a black
ball, but if the drawn ball is black it is replaced and a red ball
is thrown out of the urn. When only black balls remain, the
roles of the colors are swapped and the process restarts. We
prove that the resulting Markov chain is transient but that if we
throw out a ball every time the colors swap, the process is
recurrent. We show that the embedded process obtained by observing
the number of balls in the urn at the swapping times has a scaling
limit that is essentially the square of a Bessel diffusion. We
consider an oriented percolation model naturally associated with
the urn process, and obtain detailed information about its
structure, showing that the open subgraph is an infinite tree with
a single end. We also study a natural continuous-time embedding of
the urn process that demonstrates the relation to the simple
harmonic oscillator; in this setting, our transience result
addresses an open problem in the recurrence theory of
two-dimensional linear birth and death processes due to Kesten and
Hutton. We obtain results on the area swept out by the process. We
make use of connections between the urn process and birth--death
processes, a uniform renewal process, the Eulerian numbers, and
Lamperti's problem on processes with asymptotically small drifts;
we prove some new results on some of these classical objects that
may be of independent interest. For instance, we give sharp new
asymptotics for the first two moments of the counting function of
the uniform renewal process. Finally, we discuss some related
models of independent interest, including a ``Poisson earthquakes''
Markov chain on the homeomorphisms of the plane.
\end{abstract}

% KEYWORDS
%
\begin{keyword}[class=AMS]
\kwd[Primary ]{60J10}
\kwd[; secondary ]{60J25}
\kwd{60K05}
\kwd{60K35}.
\end{keyword}
\begin{keyword}
\kwd{Urn model}
\kwd{recurrence classification}
\kwd{oriented percolation}
\kwd{uniform renewal process}
\kwd{two-dimensional linear birth and death process}
\kwd{Bessel process}
\kwd{coupling}
\kwd{Eulerian numbers}.
\end{keyword}

\end{frontmatter}

%s1 ###
\section{Introduction}
\label{secIntro}

Urn models have a venerable history in probability theory, with
classical contributions
having been made by the Bernoullis and Laplace, among others. The
modern view of many urn
models is as prototypical reinforced stochastic processes. Classical
urn schemes
were often employed as ``thought experiments'' in which to frame
statistical questions; as stochastic
processes, urn models have wide-ranging applications in economics, the
physical sciences,
and statistics. There is a large literature on urn
models and their applications---see, for example, the monographs
\cite{JK,Mahmoud} and the surveys \cite{BK,Review}---and some important contributions have been made in the last few
years: see, for example, \cite{Ja,Fl}.

A generalized \Polya urn with 2 types of ball, or 2 colors, is a
discrete-time Markov chain $(X_n,Y_n)_{n
\in\Z_+}$ on $\Z_+^2$, where $\Z_+ :=\{0,1,2,\ldots\}$. The
possible transitions of the chain are specified
by a $2 \times2$ \textit{reinforcement matrix} $A = (a_{ij})_{i,j=1}^2$
and the transition probabilities depend on the current state:
%
%e1 ###
\begin{eqnarray}\label{eqtrprob}
\P\bigl((X_{n+1},Y_{n+1})=(X_n+a_{11},Y_n+a_{12}) \bigr) &=& \frac
{X_n}{X_n+Y_n}, \nonumber\\[-8pt]\\[-8pt]
\P\bigl((X_{n+1},Y_{n+1})=(X_n+a_{21},Y_n+a_{22}) \bigr) &=& \frac
{Y_n}{X_n+Y_n}.\nonumber
\end{eqnarray}
This process can be viewed as an urn which at time $n$ contains $X_n$
red balls and $Y_n$ black balls. At each
stage, a ball is drawn from the urn at random, and then returned
together with $a_{i1}$ red balls and
$a_{i2}$ black balls, where $i=1$ if the chosen ball is red and $i=2$
if it is black.

A fundamental problem is to study the long-term behavior of
$(X_n,Y_n)$, defined by (\ref{eqtrprob}),
or some function
thereof, such as the fraction of red balls~$X_n/\allowbreak(X_n+Y_n)$. In many cases,
coarse asymptotics for such
quantities are gover\-ned
by the eigenvalues of the reinforcement matrix $A$ (see, e.g., \cite{AK} or \cite{AN}, Section V.9).
However, there are some interesting special cases (see, e.g.,~\cite{PV1999}), and analysis of
finer behavior is in several cases still an open problem.

A large body of asymptotic theory is known under
various conditions on~$A$ and its eigenvalues.
Often it is assumed
that all $a_{ij}\ge0$, for example, $A=$
$ \bigl[
{{1 \atop 0} \enskip {0 \atop 1}}
\bigr]$ specifies the \textit{standard \Polya urn},
while $A=  \bigl[
{{a\atop b }\enskip {b \atop a}}
\bigr]$ with $a, b
>0$ specifies a \textit{Friedman urn}.

In general, the entries $a_{ij}$ may be negative, meaning that balls
can be
thrown away as well as added, but nevertheless in the literature \emph
{tenability} is usually imposed. This
is the condition that regardless of the stochastic path taken by the process,
it is never required to remove a ball of a color not currently present
in the urn.
For example, the \textit{Ehrenfest
urn},
which models the diffusion of a gas between two chambers of a box,
is tenable despite its reinforcement matrix $ \bigl[
{{ -1 \atop 1 }\enskip{ 1 \atop -1}}
\bigr]$ having some negative entries.

Departing from tenability, the \textit{OK Corral model} is the 2-color
urn with reinforcement matrix
$ \bigl[
{{0 \atop -1}\enskip{ -1 \atop 0}}
\bigr]$. This\vspace*{1pt} model for destructive competition was studied by Williams
and McIlroy \cite{WM} and Kingman \cite{K99} (and
earlier as a~stochastic version of Lanchester's combat model; see,
e.g., \cite{watson}
and references therein). Kingman and Volkov \cite{KV} showed that the OK
Corral model can be viewed as a time-reversed Friedman urn with $a=0$
and $b=1$.

In this paper, we will study the 2-color urn model with reinforcement matrix
%
%e2 ###
\begin{equation}
\label{eqharmm}
A = \left[
\matrix {0 & 1 \cr -1 & 0
}
\right].
\end{equation}
To reiterate the urn model, at each time period
we draw a ball at random from the urn; if it is red, we replace it
and add an additional black ball, if it is black we replace it and
throw out a red ball.
The eigenvalues of $A$ are
$\pm i$, corresponding to the ordinary differential equation
$\dot{{\bf v}}=A{\bf v}$,
which governs the phase diagram of the simple harmonic oscillator. This
explains the name \textit{simple harmonic urn}.
Na\"{\i}vely, one might hope
that the behavior of the Markov chain is closely related to the paths
in the phase diagram.
We will see that it is, but that the exact behavior is somewhat more subtle.

%s2 ###
\section{Exact formulation of the model and main results}
\label{secResults}

%s2.1 ###
\subsection{The simple harmonic urn process}

The definition of the process~given by the transition probabilities
(\ref{eqtrprob}) and
the matrix (\ref{eqharmm}) only makes sense for $X_n,Y_n\ge0$;
however, it is easy to see that
almost surely (a.s.) $X_n < 0$ eventually.
Therefore, we reformulate the process
$(X_n,Y_n)$ rigorously as follows.

\begin{figure}

\includegraphics{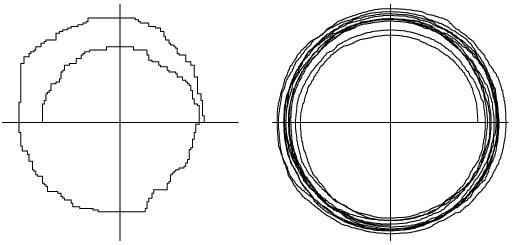}

\caption{Two sample trajectories of the simple harmonic urn process,
starting at $(50,0)$ and running for about $600$ steps
(\textup{left}) and starting at $(1000,0)$ and running for $100\mbox{,}000$ steps
(\textup{right}).} \label{Fig2}
\end{figure}

For $z_0 \in\N:= \{1,2,\ldots\}$ take $(X_0,Y_0) = (z_0,0)$;
we start on the positive $x$-axis for convenience but
the choice of initial state does not affect any of our asymptotic
results.
For $n \in\Z_+$, given $(X_n,Y_n) = (x,y) \in\Z^2 \setminus\{
(0,0)\}$,
we define the transition law of the process by
%
%e3 ###
\begin{equation}
\label{tranprobs}
(X_{n+1}, Y_{n+1} )=
\cases{
\bigl(x,y+\sgn(x)\bigr), &\quad  with probability $\displaystyle
\frac{|x|}{|x|+|y|}$,\cr
\bigl(x- \sgn(y),y\bigr), &\quad with probability $\displaystyle \frac{|y|}{|x|+|y|}$,
}
\end{equation}
where $\sgn(x)=-1,0,1$ if $x<0$, $x=0$, $x>0$, respectively.
The process~$(X_n,\allowbreak Y_n)_{n \in\Z_+}$
is an irreducible Markov chain with state-space
$\Z^2 \setminus\{ (0,0) \}$.
See Figure~\ref{Fig2} for some simulated trajectories of the simple harmonic urn
process.

%f1 ###

Let $\nu_0:=0$, and recursively define stopping times
\[
\nu_k := \min\{n>\nu_{k-1}\dvtx X_n Y_n =0\}\qquad  (k \in\N),
\]
where throughout the paper we adopt the usual convention $\min
\varnothing\,{:=}\,\infty$.~Thus, $(\nu_k)_{k \in\N}$ is the sequence
of times at which the process visits one of the axes.

It is easy to see that every $\nu_k$ is almost surely finite.
Moreover, by
construction, the process $(X_{\nu_k},Y_{\nu_k})_{k \in\N}$ visits
in cyclic (anticlockwise)
order
the half-lines \mbox{$\{y>0\}$,} $\{x<0\}$, $\{y<0\}$, $\{x>0\}$.
It is natural
(and fruitful) to consider the embedded process $(Z_k)_{k \in\Z_+}$ obtained
by taking
$Z_0 := z_0$ and $Z_k := | X_{\nu_k} | + | Y_{\nu_k}|$ ($k \in\N$).

If $(X_n,Y_n)$ is viewed
as a random walk on $\Z^2$, the process $Z_k$
is the embedded process of the
distances from $0$ at the instances of hitting the axes.
To interpret the process $(X_n,Y_n)$ as the urn model
described in Section \ref{secIntro},
we need a slight modification to the description there.
Starting with $z_0$ red balls, we run the process as described in
Section \ref{secIntro},
so the process traverses the first quadrant via an up/left path
until the red balls run out (i.e., we first hit the half-line $\{ y >0\}
$). Now we interchange the
roles of the red and black balls, and we still use $y$ to count the
black balls, but we switch to using $-x$ to count the number of red balls.
Now the process traverses the second quadrant via a left/down path
until the
black balls run out, and so on. In the urn model, $Z_k$ is the number
of balls remaining in the urn
when the urn becomes monochromatic for the $k$th time ($k \in\N$).

The strong Markov property and the transition law of $(X_n,Y_n)$
imply that~$Z_k$ is an irreducible Markov chain on $\N$.
Since our two Markov chains just described are irreducible, there
is the usual recurrence/transience dichotomy, in that either
the process is recurrent, meaning that with probability $1$ it returns
infinitely
often to any finite subset of the state space, or it is transient,
meaning that with probability $1$ it eventually escapes to infinity.
Our main question is whether the process $Z_k$
is recurrent or transient. It is easy to see that, by the nature of the
embedding, this also
determines whether the urn model $(X_n, Y_n)$ is recurrent or transient.

\begin{thm}
\label{th:main1}
$\!\!\!\!\!$The process $Z_k$ is transient; hence so is the process~$(X_n, Y_n)$.
\end{thm}

Exploiting a connection between the increments of the
process $Z_k$ and a~renewal process whose inter-arrival times
are uniform on $(0,1)$ will enable us to prove the following basic\vadjust{\goodbreak}
result.\vfill\eject

\begin{thm}
\label{th:zinc}
Let $n \in\N$. Then
%
%e4 ###
\begin{equation}
\label{zinc}
\E[ Z_{k+1} \midd Z_k = n ]
= n + \tfrac{2}{3} + O ( \re^{\alpha_1 n} )
\end{equation}
as $n \to\infty$, where $\alpha_1 + \beta_1 i = -(2.088843\dots) +
(7.461489\dots)i$
is a root of $\lambda- 1 + \re^{-\lambda} = 0$.
\end{thm}

The error term in (\ref{zinc}) is sharp, and we obtain it from
new (sharp) asymptotics for the uniform renewal process:
see Lemma \ref{lem: asymptotic} and Corollary \ref{uniform},
which improve on known results. To prove Theorem \ref{th:main1},
we need more than Theorem \ref{th:zinc}: we need to know about
the second moments of the increments of $Z_k$, amongst other things; see
Section \ref{S:zinc}. In fact, we prove Theorem \ref{th:main1}
using martingale arguments applied to $h(Z_k)$ for a well-chosen
function $h$; the analysis of the function $h(Z_k)$ rests
on a recurrence relation satisfied by the transition probabilities
of $Z_k$, which are related to the Eulerian numbers (see Section \ref
{secRubin}).

%s2.2 ###
\subsection{The leaky simple harmonic urn}
\label{S:leaky}

In fact the transience demonstrated in Theorem \ref{th:main1} is
rather delicate, as one can see by simulating the process. To
illustrate this, we consider a slight modification of the process,
which we call the \textit{leaky} simple harmonic urn. Suppose that
each time the roles of the colors are reversed, the addition of
the next ball of the new color causes one ball of the other
color to leak out of the urn; subsequently the usual simple
harmonic urn transition law applies. If the total number of balls
in the urn ever falls to one, then this modified rule causes the
urn to become monochromatic at the next step, and again it
contains only one ball. Thus, there will only be one ball in total
at all subsequent times, although it will alternate in color. We
will see that the system almost surely does reach this steady
state, and we obtain almost sharp tail bounds on the time that it
takes to do so. The leaky simple harmonic urn arises naturally in
the context of a percolation model associated to the simple
harmonic urn process, defined in Section \ref{S: percolation}
below.

As we did for the simple harmonic urn, we will represent the leaky urn
by a Markov chain $(X'_n, Y'_n)$.
For this version of the model,
it turns out to be more convenient to start just above the axis;
we take $(X'_0,Y'_0)=(z_0,1)$, where $z_0 \in\N$.
The distribution of $(X'_{n+1},Y'_{n+1})$ depends only on $(X'_n, Y'_n)
= (x,y)$.

If $x y \neq0$, the transition law is the same as that of the simple harmonic
urn process. The difference is when $x = 0$ or $y=0$; then the transition
law is
\begin{eqnarray*}
( X'_{n+1} , Y'_{n+1} ) & = &\bigl( - \sgn(y) , y - \sgn(y)
\bigr) \qquad (x = 0), \\
( X'_{n+1} , Y'_{n+1} ) & = &\bigl( x - \sgn(x) , \sgn(x) \bigr)\qquad\hspace*{10pt}  (y = 0) .
\end{eqnarray*}
Now $(X'_n,Y'_n)$ is a reducible Markov chain whose state-space
has two communicating classes, the closed class $\mathcal{C} = \{ (x,y)
\in\Z^2 \dvtx | x | + | y| = 1\}$
and the class $\{ (x,y) \in\Z^2 \dvtx |x| + |y| \geq2 \}$; if the
process enters the closed
class $\mathcal{C}$ it remains there for ever, cycling round the origin.
Let $\tau$ be the hitting time of the set $\mathcal{C}$, that is
\[
\tau:= \inf\{ n \in\mathbb{Z} \dvtx |X'_n| + |Y'_n| = 1 \} .
\]

\begin{thm}
\label{th:leaky} For the leaky urn, $\P( \tau< \infty) =1$.
Moreover, for any $\eps>0$, $\E[ \tau^{1-\eps}] < \infty$ but $\E
[ \tau^{1+\eps}] = \infty$.
\end{thm}

In contrast, Theorem \ref{th:main1} implies that the analogue of
$\tau$ for the ordinary urn process has $\P(\tau= \infty) >0$ if
$z_0 \geq2$.

%s2.3 ###
\subsection{The noisy simple harmonic urn}
\label{S:noisy}

In view of Theorems \ref{th:main1} and \ref{th:leaky}, it is
natural to ask about the properties of the hitting time $\tau$ if
at the time when the balls of one color run out we only discard a
ball of the other color with some probability $p \in(0,1)$. For
which $p$ is $\tau$ a.s. finite? (Answer: for $p \geq1/3$; see
Corollary \ref{noisyclassification} below.)

We consider the following natural generalization of the model
specified by~(\ref{tranprobs}) in order to probe more
precisely the recurrence/transience transition. We call this
generalization the \textit{noisy} simple harmonic urn process.
In a~sense that we will describe,
this model includes the leaky urn and also the
intermittent leaky urn mentioned at the start of this section.
The basic idea is to throw out (or add)
a random number of balls at each
time we are at an axis, generalizing the idea of the leaky urn.
It is more convenient here to work with irreducible Markov chains,
so we introduce a ``barrier'' for our process. We now describe the model
precisely.

Let $\kappa,\kappa_1,\kappa_2,\ldots$ be a sequence of
i.i.d. $\Z$-valued
random variables such that
%
%e5 ###
\begin{equation}
\label{kappabound}
\E\bigl[ \re^{\lambda|\kappa|} \bigr] < \infty
\end{equation}
for some $\lambda>0$, so in particular $\E[\kappa]$ is finite.
We now define the Markov chain $(\tX_n,\tY_n)_{n \in\Z_+}$ for the
noisy urn process. As for the leaky urn,
we start one step above the axis:
let $z_0 \in\N$, and take $(\tX_0,\tY_0) = (z_0,1)$.
For $n \in\Z_+$, given $(\tX_n,\tY_n) = (x,y) \in\Z^2 \setminus\{
(0,0)\}$,
we define the transition law as follows. If $xy \neq0$, then
\[
(\tX_{n+1}, \tY_{n+1} )=
\cases{
\bigl(x,y+\sgn(x)\bigr), &\quad with probability $\displaystyle
\frac{|x|}{|x|+|y|}$,\cr
\bigl(x- \sgn(y),y\bigr), &\quad with probability $\displaystyle \frac{|y|}{|x|+|y|}$,
}
\]
while if $x = 0$ or $y=0$ we have
\begin{eqnarray*}
(\tX_{n+1} , \tY_{n+1} ) & = & \bigl( - \sgn(y), \sgn(y) \max(1, |y| -
\kappa_n) \bigr)\qquad  (x=0), \\
(\tX_{n+1} , \tY_{n+1} ) & = & \bigl(\sgn(x) \max(1, |x| - \kappa_n),
\sgn(x) \bigr)\qquad \hspace*{10pt} (y=0).
\end{eqnarray*}
In other words, the transition law is the same
as (\ref{tranprobs}) except when
the process is on an axis at time $n$, in which case
instead of just moving one step away in the anticlockwise
perpendicular direction
it also moves an additional distance $\kappa_n$
parallel to the axis toward the origin
(stopping distance 1 away if it would otherwise
reach the next axis or overshoot).
Then $(\tX_n, \tY_n)_{n \in\Z_+}$ is an irreducible Markov
chain on $\Z^2 \setminus\{(0,0)\}$. The case where
$\P( \kappa= 0 ) = 1$ corresponds
to the original process $(X_n, Y_n)$ starting one unit later in time.

A fundamental random variable is the first passage time
to within distance~1 of the origin:
\[
\tau:= \min\{ n \in\Z_+ \dvtx | \tilde X_n | + | \tilde Y_n | = 1 \}=
\min\{ n \in\Z_+ \dvtx (\tX_n, \tY_n) \in\mathcal{C} \} .
\]
Define a sequence of stopping times $\tilde\nu_k$
by setting $\tilde\nu_0 := -1$ and for $k \in\N$,
\[
\tilde\nu_k := \min\{ n > \tilde\nu_{k-1} \dvtx \tX_n \tY_n = 0 \} .
\]
As an analogue of $Z_k$, set $\tZ_0 := z_0$ and
for $k \in\N$ define
\[
\tZ_k := \max\{ | \tX_{1+\tilde\nu_k }| , | \tY_{1+\tilde\nu_k
} | \} =
| \tX_{1+\tilde\nu_k }| + | \tY_{1+\tilde\nu_k } | -1 ;
\]
then $(\tZ_k)_{k \in\Z_+}$ is an irreducible Markov chain
on $\N$.
Define the return-time to
the state $1$ by
%
%e6 ###
\begin{equation}
\label{returntime}
\tau_q := \min\{ k \in\N\dvtx \tZ_k = 1 \},
\end{equation}
where the subscript $q$ signifies the fact that a time unit is one
traversal of a~quadrant here. By our embedding, $\tau= \tilde\nu
_{\tau_q}$.

Note that in the case $\P( \kappa=0)=1$, $(\tZ_k)_{k \in\Z_+}$ has
the same distribution as the
original $(Z_k)_{ k \in\Z_+}$. The noisy urn with $\P(\kappa=1)=1$
coincides with the leaky urn
described in Section \ref{S:leaky} up until the time $\tau$ (at which
point the leaky urn becomes
trapped in $\mathcal{C}$). Similarly, the embedded process $\tZ_k$ with
$\P(\kappa=1)=1$ coincides
with the process of distances from the origin of the leaky urn at the
times that it visits the axes,
up until time $\tau_q$ (at which point the leaky urn remains at
distance $1$ forever). Thus, in the
$\P(\kappa=1)=1$ cases of all the results that follow in this
section, $\tau$ and $\tau_q$ can be
taken to be defined in terms of the leaky urn $(X'_n,Y'_n)$.

The next result thus includes Theorem \ref{th:main1}
and the first part of Theorem~\ref{th:leaky}
as special cases.

\begin{thm}
\label{thm:mixed}
Suppose that $\kappa$ satisfies (\ref{kappabound}). Then
the process $\tZ_k$ is:
\begin{itemize}[(iii)]
\item[(i)] transient if $\E[\kappa] < 1/3$;
\item[(ii)] null-recurrent if $1/3 \leq\E[\kappa] \leq2/3$;
\item[(iii)] positive-recurrent if $\E[\kappa] > 2/3$.
\end{itemize}
\end{thm}

Of course, part (i) means that $\P(\tau_q < \infty)<1$,
part (ii) that $\P(\tau_q < \infty)=1$ but $\E[\tau_q] = \infty$,
and part (iii) that $\E[ \tau_q ]< \infty$. We can in fact obtain
more information
about the tails of $\tau_q$.

\begin{thm}
\label{thm:moments}
Suppose that $\kappa$ satisfies (\ref{kappabound})
and $\E[\kappa] \geq1/3$. Then
$\E[ \tau_q^p ] < \infty$ for $p < 3 \E[\kappa] - 1$ and
$\E[ \tau_q^p ] = \infty$ for $p > 3 \E[\kappa] - 1$.
\end{thm}

It should be possible, with some extra work, to show that $\E[
\tau_q^p ] = \infty$ when $p = 3 \E[\kappa] - 1$, using the
sharper results of \cite{ai} in place of those from \cite{aim}
that we use below in the proof of Theorem \ref{thm:moments}.

In the recurrent case, it is of interest to obtain more detailed
results on the tail of~$\tau$
(note that there is a change of time between $\tau$ and $\tau_q$).
We obtain the following upper and lower bounds, which are close to sharp.

\begin{thm}
\label{thm:moments2}
Suppose that $\kappa$ satisfies (\ref{kappabound}) and $\E[\kappa]
\geq1/3$.\vspace*{1pt}
Then
$\E[ \tau^p] < \infty$ for $p < \frac{3 \E[ \kappa] -1}{2}$ and
$\E[ \tau^p] = \infty$ for $p > \frac{3 \E[ \kappa] -1}{2}$.
\end{thm}

Theorems \ref{thm:mixed} and \ref{thm:moments2} have
an immediate corollary for the noisy urn process $(\tilde X_n, \tilde Y_n)$.

\begin{corollary}
\label{noisyclassification} Suppose that $\kappa$ satisfies
(\ref{kappabound}). The noisy simple harmonic urn process $(\tilde
X_n, \tilde Y_n)$ is recurrent if $\E[\kappa] \geq1/3$ and
transient if $\E[ \kappa] < 1/3$. Moreover, the process is
null-recurrent if $1/3 \leq\E[\kappa] < 1$ and positive-recurrent
if $\E[ \kappa] > 1$.
\end{corollary}

This result is close to sharp but leaves open the question of whether
the process
is null- or positive-recurrent when $\E[ \kappa]=1$ (we suspect the former).

We also study the distributional limiting behavior of
$\tZ_k$ in the appropriate scaling regime when
$\E[\kappa] < 2/3$. Again the case $\P(\kappa= 0) =1$ reduces
to the original~$Z_k$.

\begin{thm}
\label{thm:difflimit}
Suppose that $\kappa$ satisfies (\ref{kappabound})
and that $\E[\kappa] < 2/3$.
Let $(D_t)_{t \in[0,1]}$ be a diffusion process taking values in $\R
_+ := [0,\infty)$
with $D_0 = 0$ and
infinitesimal mean $\mu(x)$ and variance $\sigma^2(x)$ given for $x
\in\R_+$
by
\[
\mu(x) = \tfrac{2}{3} -\E[\kappa],\qquad  \sigma^2 (x) = \tfrac{2}{3} x.
\]
Then as $k \to\infty$,
\[
( k^{-1} \tZ_{ k t} )_{t \in[0,1]} \to(D_t)_{t \in[0,1]} ,
\]
where the convergence is in the sense of finite-dimensional distributions.
Up to multiplication by a scalar,
$D_t$ is the square of a Bessel process with parameter $4 - 6 \E
[\kappa] > 0$.
\end{thm}

Since a Bessel process with parameter $\gamma\in\N$
has the same law as the norm of a $\gamma$-dimensional Brownian motion,
Theorem \ref{thm:difflimit} says, for example,
that if $\E[\kappa] =0$ (e.g., for the original urn process)
the scaling limit of $\tZ_t$ is a scalar multiple of the norm-square
of $4$-dimensional Brownian motion,
while if $\E[ \kappa] = 1/2$
the scaling limit is a scalar multiple of the square
of one-dimensional Brownian motion.

To finish this section, consider the
\textit{area} swept out by the path of the noisy simple harmonic
urn on its first excursion (i.e., up to time $\tau$).
Additional motivation for studying this
random quantity is provided by the percolation model of Section~\ref{S: percolation}. Formally, for $n \in\N$ let $T_n$ be the area
of the triangle with vertices $(0,0)$, $(\tilde X_{n-1}, \tilde Y_{n-1})$,
and $(\tilde X_{n}, \tilde Y_{n})$, and define $A := \sum_{n=1}^{\tau
} T_n$.

\begin{thm}
\label{thm:area}
Suppose that $\kappa$ satisfies (\ref{kappabound}).
\begin{itemize}[(ii)]
\item[(i)] Suppose that $\E[\kappa]< 1/3$. Then $\P( A = \infty) > 0$.
\item[(ii)] Suppose that $\E[\kappa]\ge1/3$. Then $\E[ A^p ] <
\infty$ for $p < \frac{3 \E[\kappa] -1}{3}$.
\end{itemize}
\end{thm}

In particular, part (ii) gives us information about the leaky urn
model, which corresponds to the case where $\P(\kappa=1)=1$, at
least up until the hitting time of the closed cycle; we can still
make sense of the area swept out by the leaky urn up to this
hitting time. We then have $\E[A^p] < \infty$ for $p<2/3$, a
result of significance for the percolation model of the next
section. We suspect that the bounds in Theorem~\ref{thm:area}(ii) are
tight. We do not prove this but have the following result in the
case $\P(\kappa=1)=1$.

\begin{thm}
\label{thm:areamean}
Suppose $\P( \kappa=1) =1$ (or equivalently take the leaky urn). Then
$\E[ A] = \infty$.
\end{thm}

%s2.4 ###
\subsection{A percolation model}\label{S: percolation}

Associated to the simple harmonic urn is a~per\-colation model
which we describe in this section. The percolation model, as well as
being of
interest in its own right, couples many different instances
of the simple harmonic urn, and exhibits naturally an instance
of the leaky version of the urn in terms of the planar dual percolation
model. Our results
on the simple harmonic urn will enable us to establish
some interesting properties of the percolation model.

The simple harmonic urn can be viewed as a spatially inhomogeneous
random walk
on a directed graph whose vertices are $\Z^2 \setminus\{ (0,0)\}$; we
make this statement
more precise shortly.
In this section, we will view
the simple harmonic urn process not as a random path through a
predetermined directed graph
but as a deterministic path through a random directed graph.
To do this, it is helpful to consider a slightly larger state-space
which keeps track
of the number of times that the urn's path has wound around the origin.
We construct
this state-space as the vertex set of a graph $G$ that is embedded in
the Riemann surface
$\mathcal{R}$ of the complex logarithm, which is the universal cover of
$\mathbb{R}^2 \setminus\{(0,0)\}$.
To construct $G$, we take the usual square-grid lattice and delete the
vertex at the origin to obtain a graph on the vertex set $\mathbb
{Z}^2\setminus\{(0,0)\}$.
Make this into a directed graph by orienting each edge in the direction
of increasing argument;
the paths of the simple harmonic urn only
ever traverse edges in this direction. Leave undirected those edges
along any of the coordinate axes;
the paths of the simple harmonic urn never traverse these edges.
Finally, we let $G$ be the lift of this graph to the covering surfa\-ce~$\mathcal{R}$.

%f2 ###
\begin{figure}

\includegraphics{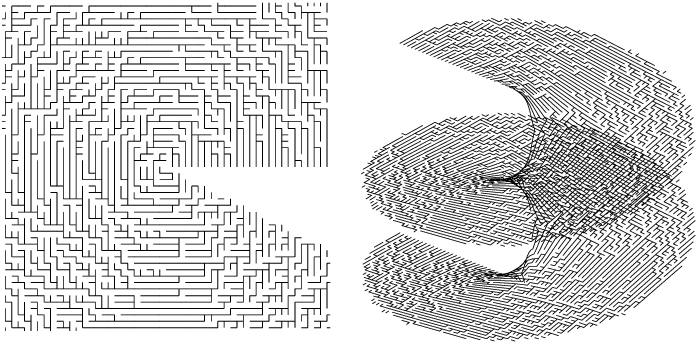}

\caption{Simulated realizations of the simple harmonic urn
percolation model: on a single sheet of $\mathcal{R}$ (\textup{left})
and on a larger section (\textup{right}).} \label{percpic}
\vspace*{-3pt}
\end{figure}

We will interpret a path of the simple harmonic urn as the unique
oriented path from some starting vertex through a random subgraph $H$
of $G$.
For each vertex $v$ of $G$, the graph $H$ has precisely one of the
out-edges from $v$ that are in $G$. If the projection of $v$ to
$\mathbb{Z}^2$ is $(x,y)$,
then the graph $H$ contains the edge from $v$ that projects onto the
edge from $(x,y)$ to $(x-\sgn(y), y)$ with probability
$|y|/(|x|+|y|)$, and otherwise it contains
the edge from $v$ that projects onto the edge from $(x,y)$ to $(x, y +
\sgn(x))$. These choices are to be made independently for all vertices
$v$ of $G$.
In particular, $H$ does not have any edges that project onto either of
the coordinate axes.
The random directed graph $H$ is an oriented percolation model that
encodes a coupling of many different paths of the simple harmonic urn.
To make this precise, let $v_0$ be any vertex of $G$.
Then there is a unique oriented path $v_0, v_1, v_2, \dots$ through~$H$.
That is, $(v_i,v_{i+1})$ is an edge of $H$ for each $i \ge0$.
Let the projection of~$v_i$ from $\mathcal{R}$ to $\mathbb{R}^2$ be
the point $(X_i, Y_i)$.
Then the sequence $(X_i, Y_i)_{i=0}^{\infty}$ is a~sample of the
simple harmonic urn process.
If $w_0$ is another vertex of $G$, with unique oriented path $w_0, w_1,
w_2, \ldots,$ then its projection to~$\mathbb{Z}^2$ is also a sample
path of the simple harmonic urn process, but we will show (see Theorem
\ref{thm:perc} below)
that with probability one the two paths eventually couple, which is to
say that there exist random finite $m \ge0$ and $n \ge0$ such that
for all $i \ge0$ we have $v_{i+m} = w_{i+n}$. Thus, the percolation
model encodes\vadjust{\goodbreak} many coalescing copies of the simple
harmonic urn process. Next, we show that it also encodes many copies of
the leaky urn of Section~\ref{S:leaky}.

We construct another random graph $H'$ that is the dual percolation
model to $H$. We begin with the planar dual of the square-grid lattice,
which is another square-grid lattice with vertices at the points
$(m+1/2,n+1/2)$, $m,n \in\mathbb{Z}$. We orient all the edges in the
direction of decreasing argument, and lift to the covering surface
$\mathcal{R}$ to obtain the dual graph $G'$. Now let $H'$ be the directed
subgraph of $G'$ that consists of all those edges of $G'$ that do not
cross an edge of $H$. It turns out that $H'$ can be viewed as an
oriented percolation model that encodes a coupling of many different
paths of the leaky simple harmonic urn.

To explain this, we define a mapping $\Phi$ from the vertices of $G'$
to $\mathbb{Z}^2$. Let $(x,y)$ be the coordinates of the projection of
$v \in G'$ to the shifted square lattice $\mathbb{Z}^2 + (1/2,1/2)$. Then
\[
\Phi(v) = \bigl( x + \tfrac{1}{2} \sgn y, -\bigl(y - \tfrac{1}{2} \sgn x\bigr) \bigr) .
\]
Thus, we project from $\mathcal{R}$ to $\mathbb{R}^2$, move to the
nearest lattice point in the clockwise direction, and then reflect in
the $x$-axis.
If $v_0$ is any vertex of $H'$, there is a unique oriented path $v_0,
v_1, v_2, \dots$ through $H'$, this time winding clockwise. Take $v_0
= (z_0 - 1/2, 1/2)$.
A little thought shows that the sequence $(X'_i, Y'_i) = \Phi(v_i )$
has the distribution of the leaky simple harmonic urn process. This is
because the choice of edge in $H'$ from $v$ is determined by the choice
of edge in $H$ from the nearest point of $G$ in the clockwise direction.
The map $\Phi$ is not quite a graph homomorphism onto the square
lattice because of its behavior at the axes; for example, it
sends $(3\frac{1}{2},\frac{1}{2})$
and $(3\frac{1}{2}, -\frac{1}{2})$ to $(4,0)$ and $(3, 1)$, respectively.
The decrease of $1$ in the $x$-coordinate corresponds to the leaked
ball in the leaky urn model. If some $v_i$ has projection $(x_i, y_i)$
with $|x_i| + |y_i| = 1$,
then the same is true of all subsequent vertices in the path. This
corresponds to the closed class $\mathcal{C}$.

From results on our urn processes, we will deduce the following quite
subtle properties of the percolation model $H$.
Let $I(v)$ denote the number of vertices in the \textit{in-graph} of the
vertex $v$ in $H$, which is the subgraph of $H$ induced by all
vertices from which it is possible to reach $v$ by following an
oriented path.\vspace*{-3pt}

\begin{thm}
\label{thm:perc} Almost surely, the random oriented graph $H$ is,
ignoring orientations, an infinite tree with a single
semi-infinite end in the out direction. In particular, for any
$v$, $I(v) < \infty$ a.s. and moreover $\E[ I(v)^p ] < \infty$
for any $p<2/3$; however, $\E[ I(v) ] =\infty$.

The dual graph $H'$ is also an infinite tree a.s., with a single
semi-infinite end in the out direction. It has a doubly-infinite
oriented path and the in-graph of any vertex not on this path is finite a.s.
\end{thm}

%s2.5 ###
\subsection{A continuous-time fast embedding of the simple harmonic urn}
\label{S: fast}

$\!\!\!$There is a\vadjust{\goodbreak} natural continuous-time embedding of the simple harmonic urn process.
Let $(A(t),B(t))_{t \in\R_+}$ be a $\mathbb{Z}^2$-valued
continuous-time Markov
chain with $A(0) = a_0$, $B(0) = b_0$, and transition rates
\begin{eqnarray*}
\P\bigl(A(t+ \ud t) = A(t) - \operatorname{sgn}(B(t))\bigr) & =& | B(t) | \,\ud t , \\
\P\bigl(B(t+ \ud t) = B(t) + \operatorname{sgn}(A(t))\bigr) & =& | A(t) | \,\ud t .
\end{eqnarray*}
Given that $(A(t),B(t)) = (a,b)$, the wait until the next jump after
time $t$ is
an exponential random variable with mean $1/(|a|+|b|)$.
The next jump is a~change in the first coordinate with probability
$|b|/(|a|+|b|)$, so the process considered at its sequence
of jump times does indeed follow the law of the simple harmonic urn.
Note that the process does not explode in finite time since the jump rate
at $(a,b)$ is $|a|+|b|$, and $|X_n| + |Y_n| = O(n)$ (as jumps are of
size $1$),
so $\sum_n ( |X_n| + |Y_n| )^{-1} = \infty$ a.s.

The process $(A(t),B(t))$ is an example of a \textit{two-dimensional
linear birth and death process}. The recurrence classification of
such processes defined on $\Z_+^2$ was studied by
Kesten \cite{kesten} and Hutton \cite{hutton}. Our case (which has
$B_{1,1}+B_{2,2}=0$ in their notation) was not covered by the
results in \cite{kesten,hutton}; Hutton
remarks (\cite{hutton}, page 638), that ``we do not yet know whether
this case is recurrent or transient.'' In the $\Z_+^2$ setting
of \cite{kesten,hutton}, the boundaries of the quadrant would
become absorbing in our case. The model on $\Z^2$ considered here
thus seems a natural setting in which to pose the
recurrence/transience question left open by \cite{kesten,hutton}.
Our Theorem \ref{th:main1} implies that $(A(t), B(t))$ is in fact
transient.

We call $(A(t),B(t))$ the \textit{fast} embedding of the urn since
typically many jumps occur in unit time
(the process jumps faster the farther away from the origin it is).
There is another continuous-time
embedding of the urn model that is also very useful in its analysis,
the \textit{slow}
embedding described in Section \ref{secRubin} below.

The mean of the process $(A(t),B(t))$ precisely follows the simple
harmonic oscillation suggested by the name
of the model.
This fact is most neatly expressed in the complex plane $\mathbb{C}$.
Recall that a~complex martingale
is a~complex-valued stochastic process whose real and imaginary parts
are both martingales.

\begin{lemma}
\label{complexmartingale}
The process $(M_t)_{t \in\R_+}$ defined by
\[
M_t := \re^{-it} \bigl( A(t) + i B(t) \bigr)
\]
is a complex martingale. In particular, for $t > t_0$ and $z \in
\mathbb{C}$,
\[
\E[ A(t) + iB(t) \midd A(t_0) + i B(t_0) = z ] = z \re^{i(t-t_0)} .
\]
\end{lemma}

As can be seen directly from the definition,
the continuous-time
Markov chain $(A(t),B(t))$ admits a constant
invariant measure;
this fact is closely related to the ``simple harmonic
flea circus'' that we describe in Section \ref{S: stationary model}.

Returning to the dynamics of the process,
what is the expected time taken to traverse a quadrant in the fast
continuous-time embedding? Define
$\tau_f := \inf\{ t \in\R_+ \dvtx A(t) = 0 \}$.
We use the notation $\P_n ( \cdot)$ for $\P( \cdot\midd A(0) =
n,\allowbreak
B(0) = 0)$,
and similarly for $\E_n$.
Numerical calculations strongly suggest the following:

\begin{conjecture}
\label{conj:fasttime}
Let $n \in\N$. With $\alpha_1 \approx-2.0888$ as in Theorem \ref
{th:zinc} above, as $n \to\infty$,
\[
\E_n [ \tau_f ] = \pi/2 + O\bigl(\re^{\alpha_1 n}/\sqrt{n}\,\bigr) .
\]
\end{conjecture}

We present a possible approach to the
resolution of Conjecture \ref{conj:fasttime} in Section~\ref
{S:polys}; it turns out that $\E_n [\tau_f]$
can be expressed as a rational polynomial of degree $n$ evaluated at
$\re$.
The best result that we have been able to prove along the lines of
Conjecture \ref{conj:fasttime}
is the following,
which shows not only that $\E_n [ \tau_f]$ is close to $\pi/2$ but
also that
$\tau_f$ itself is concentrated about $\pi/2$.

\begin{thm}
\label{th:fasttime}
Let $n \in\N$. For any $\delta>0$, as $n \to\infty$,
%
%e8 ###
%e7 ###
\begin{eqnarray}
\label{fast1} \E_n [ \tau_f ] &=& \pi/2 + O\bigl( n^{\delta-(1/2)} \bigr), \\
\label{fast2}
\E_n [ | \tau_f - (\pi/2) |^2 ] &=& O \bigl(n^{\delta-(1/2)}\bigr) .
\end{eqnarray}
\end{thm}

In the continuous-time fast embedding
the paths of the simple harmonic urn are a discrete stochastic approximation
to continuous circular motion at angular velocity $1$, with the radius
of the motion
growing approximately linearly in line with the transience of the process.
Therefore, a natural quantity to examine is the area enclosed by a path
of the urn across the
first quadrant, together with the two co-ordinate axes. For a typical
path starting at $(n,0)$,
we would expect this to be roughly $\pi n^2 / 4$, this being the area
enclosed by a quarter-circle
of radius $n$ about the origin. We use the percolation model to obtain
an exact relation between the
expected area enclosed and the expected time taken for the urn to
traverse the first quadrant.

\begin{thm}\label{th:timearea} For $n \in\N$, for any $\delta>0$,
\[
\E_n [ \mbox{Area enclosed by a single traversal} ] = \sum_{m=1}^n
m \E_m [ \tau_f ] = \frac{\pi n^2}{4} + O \bigl(n^{(3/2)+\delta}\bigr) .
\]
\end{thm}

In view of the first equality in Theorem \ref{th:timearea} and
Conjecture \ref{conj:fasttime}, we suspect a~sharp version of the asymptotic
expression for the expected area enclosed to be
\[
\E_n [ \mbox{Area enclosed by a single traversal}] = \frac{\pi
n(n+1)}{4} + c + O \bigl(\sqrt{n} \re^{\alpha_1 n} \bigr)
\]
for some constant $c \in\R$.

%s2.6 ###
\subsection{Outline of the paper and related literature}

The outline of the remainder the the paper is as follows. We begin
with a study of the discrete-time embedded process $Z_k$ in the
original urn model. In Section \ref{secRubin}, we use a decoupling
argument to obtain an explicit formula, involving the Eulerian
numbers, for the transition probabilities of $Z_k$. In
Section \ref{S: drift}, we study the drift of the process $Z_k$ and
prove Theorem \ref{th:zinc}. We make use of an attractive coupling
with the renewal process based on the uniform distribution. Then
in Section \ref{S: short proof}, we give a short, stand-alone proof
of our basic result, Theorem~\ref{th:main1}. In
Section \ref{S:zinc}, we study the increments of the process~$Z_k$,
obtaining tail bounds and moment estimates. As a by-product of our
results we obtain (in Lemma \ref{lem: asymptotic} and
Corollary \ref{uniform}) sharp expressions for the first two
moments of the counting function of the uniform renewal process,
improving on existing results in the literature. In Section
\ref{sec:tildez}, we study the asymptotic behavior of the noisy
urn embedded process~$\tZ_k$, building on our results on $Z_k$.
Here, we make use of powerful results of Lamperti and others on
processes with asymptotically zero drift, which we can apply to
the process~$\tZ^{1/2}_k$. Then in Section~\ref{sec:proofs}, we
complete the proofs of
Theorems \ref{th:leaky}--\ref{thm:moments2}, \ref{thm:difflimit},~\ref{thm:area} and \ref{thm:perc}.
In Section~\ref{sec:continuous},
we study the continuous-time fast embedding described in
Section~\ref{S: fast}, and in Sections~\ref{S:fasttime}
and \ref{S:areatime} present proofs of
Theorems~\ref{thm:areamean}, \ref{th:fasttime}
and~\ref{th:timearea}. In Section~\ref{S:polys}, we give some
curious exact formulae for the expected area and time described in
Section~\ref{S: fast}. Finally, in Section~\ref{sec:misc} we
collect some results on several models that are not directly
relevant to our main theorems but that demonstrate further some of
the surprising richness of the phenomena associated with the
simple harmonic urn and its generalizations.

We finish this section with some brief remarks on modeling applications
related to the simple
harmonic urn. The simple harmonic urn model has some similarities to
R. F. Green's urn model for
cannibalism (see, e.g., \cite{pittel}). The cyclic nature of the model
is similar to that of various
stochastic or deterministic models of certain planar systems with
feedback: see for instance~\cite{fgg} and references therein. Finally, one may view the simple
harmonic urn as a gated polling
model with two queues and a single server. The server serves one queue,
while new arrivals are
directed to the other queue. The service rate is proportional to the
ratio of the numbers of
customers in the two queues. Customers arrive at the unserved queue at
times of a Poisson process of
constant rate. Once the served queue becomes empty, the server switches
to the other queue, and a
new secondary queue is started. This model gives a third
continuous-time embedding of the simple
harmonic urn, which we do not study any further in this paper. This
polling model differs from
typical polling models studied in the literature (see, e.g., \cite{mm}) in that the service rate
depends upon the current state of the system. One possible
interpretation of this unusual service
rate could be that the customers in the primary queue are\vadjust{\goodbreak} in fact
served by the waiting customers in
the secondary queue.
The crucial importance of
the behavior at the boundaries for the recurrence classification of certain such processes
was demonstrated already in \cite{mm}.\looseness=-1

%s3 ###
\section{Transition probabilities for $Z_k$}
\label{secRubin}

In this section, we derive an exact formula for the transition probabilities
of the Markov chain $(Z_k)_{k \in\Z_+}$ (see Lemma~\ref{Lemma:
transition probabilities}
below).
We use a coupling (or rather ``decoupling'')
idea that is sometimes attributed to Samuel Karlin
and Herman Rubin.
This construction was used
in \cite{KV}
to study the OK Corral gunfight model, and is
closely related to the embedding of a generic generalized \Polya urn in a
multi-type branching process \cite{AK,AN}. The construction yields
another continuous-time
embedding of the urn process, which, by way of contrast to the
embedding described in Section
\ref{S: fast}, we refer to as the \textit{slow} embedding of the urn.

We couple the segment of the urn process $(X_n,Y_n)$
between
times $\nu_k+1$ and $\nu_{k+1}$ with certain birth and death
processes, as follows. Let $\lambda_k :=1/k$.
Consider two {\it independent} $\Z_+$-valued continuous-time Markov
chains, $U(t)$
and $V(t)$, $t\in\R_+$, where $U(t)$ is a pure death process with
the transition rate
\[
\P\bigl(U(t+ \ud t)=U(t)-1 \midd U(t)=a\bigr)=\lambda_a\,\ud t ,
\]
and $V(t)$ is a pure birth process with
\[
\P\bigl(V(t+ \ud t)=V(t)+1 \midd V(t)=b\bigr)=\lambda_b\, \ud t.
\]
Set $U(0)=z$ and
$V(0)=1$.

From the standard exponential holding-time characterization
for conti\-nuous-time Markov chains and
the properties of independent exponential random\vadjust{\goodbreak}
variables, it follows that the embedded process $(U(t),V(t))$
considered at the times when either of its coordinates changes
has the same distribution as the simple harmonic urn $(X_n,Y_n)$
described above when $(X_n,Y_n)$ is traversing the first quadrant.
More precisely, let $\theta_0:=0$ and define the jump times
of the process $V(t)-U(t)$ for $n \in\N$:
\[
\theta_n :=\inf\{t>\theta_{n-1}\dvtx U(t)< U( \theta_{n-1} )\mbox{ or }V(t)>
V( \theta_{n-1})\}.
\]
Since $\lambda_b \leq1$ for all $b$, the processes $U(t)$, $V(t)$
a.s. do not explode in finite time, so $\theta_n \to\infty$
a.s. as $n \to\infty$. Define $\eta:= \min\{ n \in\N:
U(\theta_n) = 0\}$ and set
\[
T := \theta_{\eta} = \inf\{ t > 0 \dvtx U(t) = 0 \} ,
\]
the extinction time of $U(t)$.
The coupling yields the following result (cf. \cite{AN}, Section V.9.2).

\begin{lemma}
\label{RubinLemma}
Let $k \in\Z_+$ and $z \in\N$.
The sequence $(U(\theta_n),V(\theta_n))$, $n =0,1,\ldots,\eta$,
with $(U(0),V(0))=(z,1)$,
has the same
distribution as each of the following two sequences:
\begin{itemize}[(ii)]
\item[(i)] $(|X_n|,|Y_n|)$, $n = \nu_{k}+1, \ldots, \nu_{k+1}$,
conditioned on $Z_{k} = z$ and $Y_{\nu_k} =0$;
\item[(ii)] $(|Y_n|,|X_n|)$, $n = \nu_{k}+1, \ldots, \nu_{k+1}$,
conditioned on $Z_{k} = z$ and\vadjust{\goodbreak} $X_{\nu_k} =0$.
\end{itemize}
\end{lemma}

Note that we set
$V(0)=1$ since
$(X_{\nu_k+1},Y_{\nu_k+1})$ is always one step in the
``anticlockwise'' lattice direction away from $(X_{\nu_k},Y_{\nu_k})$.
Let
\[
T_w':=\inf\{t>0 \dvtx V(t)=w\}.
\]
We can represent the times $T$ and $T'_w$ as
sums of exponential random variables. Write
%
%e9 ###
\begin{equation}
\label{E: times as sums of exponential rvs}
T_z = \sum_{k=1}^{z} k \xi_k \quad \mbox{and}\quad
T'_w = \sum_{k=1}^{w-1} k \zeta_k,
\end{equation}
where $\xi_1,\zeta_1, \xi_2, \zeta_2, \ldots$ are independent
exponential random variables with mean~$1$. Then setting
$T = T_{U(0)}$, (\ref{E: times as sums of exponential rvs})
gives useful representations of $T$ and~$T'_w$.

As an immediate illustration of the power of this embedding,
observe that $Z_{k+1} \le Z_k$ if and only if
$V$ has not reached $U(0)+1$ by the
time of the extinction of $U$, that is, $T_{U(0)+1}' > T$.
But (\ref{E: times as sums of exponential rvs})
shows that $T_{z+1}'$ and $T_z$ are identically distributed continuous
random variables, entailing the following result.
\begin{lemma}
For $z \in\N$, \mbox{$\P(Z_{k+1} \le Z_k \midd Z_k = z )
= \P(T_{U(0)+1}' > T_{U(0)} ) = \frac{1}{2}$}.
\end{lemma}

We now proceed to derive from the coupling
described in Lemma \ref{RubinLemma}
an exact formula for the transition probabilities
of the Markov chain $(Z_k)_{k \in\Z_+}$.
Define $p(n,m) = \P(Z_{k+1}=m\midd Z_k=n )$. It turns out that
$p(n,m)$ may be expressed in terms of the \textit{Eulerian numbers}
$A(n,k)$, which are the positive integers defined
for $n \in\N$ by
\[
A(n,k) = \sum_{i=0}^k (-1)^i \pmatrix{n+1\cr i}(k-i)^n ,\qquad  k \in\{1,\ldots,n \} .
\]
The Eulerian numbers have several combinatorial interpretations and
have many interesting
properties; see, for example, B\'{o}na \cite{Bona}, Chapter 1.

\begin{lemma}
\label{Lemma: transition probabilities}
For $n,m \in\N$, the transition probability $p(n,m)$ is given by
\begin{eqnarray*}
p(n,m) & = & m \sum_{r=0}^m (-1)^r \frac{(m-r)^{n+m-1}}{r! (n+m-r)!}\\
& = & \frac{m}{(m+n)!} A(n+m-1,n) .
\end{eqnarray*}
\end{lemma}

We give two proofs of Lemma \ref{Lemma: transition probabilities},
both using the coupling of
Lemma \ref{RubinLemma} but in quite different ways. The first uses
moment generating
functions and is similar to calculations in \cite{KV}, while the
second involves a
time-reversal of the death process and makes use of the recurrence
relation satisfied by the Eulerian numbers.
Each proof uses ideas that will be useful later\vadjust{\goodbreak} on.

\begin{pf*}{First proof of Lemma \ref{Lemma: transition probabilities}}
By Lemma \ref{RubinLemma},
the conditional distribution of $Z_{k+1}$ on $Z_k=n$ coincides with the
distribution
conditional of $V(T)$ on $U(0) =n$. So
%
%e10 ###
\begin{equation}
\label{Zs and Ts}
\P(Z_{k+1}> m \midd Z_k=n ) =
\P\bigl( V(T ) > m \midd U(0) = n\bigr) =
\P(T_n > T'_{m+1} ) ,\hspace*{-20pt}
\end{equation}
using the representations in (\ref{E: times as sums of exponential rvs}).
Thus, from (\ref{E: times as sums of exponential rvs}) and (\ref{Zs
and Ts}), writing
\[
R_{n,m} = \sum_{i=1}^{n} i \xi_i -\sum_{j=1}^{m} j
\zeta_j,
\]
we have that $\P(Z_{k+1}> m \midd Z_k=n ) = \P( R_{n,m} > 0)$.
The density of $R_{n,m}$ can be calculated using the moment
generating function and partial fractions; for $t \geq0$,
\begin{eqnarray*}
\E[ \re^{t R_{n,m}} ] &=&\prod_{i=1}^n \frac{1/i}{1/i-t} \times
\prod_{j=1}^m \frac{1/j}{1/j+t}
=\prod_{i=1}^n \frac{1}{1-it} \times\prod_{j=1}^m \frac{1}{1+jt}
\\&=&\sum_{i=1}^n \frac{a_i}{1-it} +\sum_{j=1}^m \frac{b_j}{1+jt}
\end{eqnarray*}
for some coefficients $a_i = a_{i;n,m}$ and $b_j = b_{j;n,m}$.
Multiplying both sides of the last displayed
equality by $\prod_{i=1}^n (1- it) \prod_{j=1}^m (1 + jt)$
and setting $t = 1/i$, we obtain
\begin{eqnarray*}
a_i &=& \prod_{\stackrel{j=1}{j \neq i}}^n \frac{1}{1-(j/i)} \prod
_{k=1}^m \frac{1}{1+(k/i)}\\
&=& (-1)^{n-i} i^{n+m -1} \prod_{j=1}^{i-1} \frac{1}{i-j} \prod
_{j=i+1}^n \frac{1}{j-i} \prod_{k=1}^m \frac{1}{k+i} .
\end{eqnarray*}
Simplifying, and then proceeding similarly but taking $t = -1/j$ to
identify~$b_j$, we obtain
\[
a_i = \frac{(-1)^{n-i} i^{n+m} }{ (n-i)! (m+i)! } \quad \mbox{and}\quad
b_j = \frac{(-1)^{m-j} j^{n+m} }{ (m-j)! (n+j)! }.
\]
Consequently, the density of $R_{n,m}$ is
\[
r(x)=
\cases{
\displaystyle \sum_{i=1}^n a_i i^{-1} \re^{-x/i}, &\quad if
$x\ge0$,\cr
\displaystyle \sum_{j=1}^m b_j j^{-1} \re^{x/j}, &\quad  if $ x< 0$.
}
\]
Thus, we obtain
%
%e11 ###
\begin{eqnarray}
\label{eq:cdfofZ}
\P(Z_{k+1}>m \midd Z_k=n)&=& \P(Z_{k+1}\ge m+1 \midd Z_k=n)\nonumber \\
&=& \P(R_{n,m} \geq0)= \sum_{k=1}^n a_{k;n,m}
\nonumber\\[-8pt]\\[-8pt]
&=&\sum_{k=0}^n \frac{(-1)^{n-k} k^{n+m} }{
(n-k)! (m+k)! }=\sum_{i=0}^n \frac{(-1)^i (n-i)^{n+m}}{i! (m+n-i)!}\nonumber
\\
&=&\frac{1}{(m+n)!}\sum_{i=0}^n (-1)^i \pmatrix{m+n\cr i} (n-i)^{n+m}.\nonumber
\end{eqnarray}
It follows that
\begin{eqnarray*}
p(n,m) &=& \P(Z_{k+1}\ge m \midd Z_k=n)
- \P(Z_{k+1}\ge m+1 \midd Z_k=n)
\\
&=& \sum_{i=0}^n \frac{(-1)^i (n-i)^{n+m-1}}{i! (m-1+n-i)!}
-\sum_{i=0}^n \frac{(-1)^i (n-i)^{n+m}}{i! (m+n-i)!}\\
&=&\sum_{i=0}^n \frac{(-1)^i (n-i)^{n+m-1}} {i! (m+n-i)!}
[(m+n-i)-(n-i) ]
\\
&=&\frac m{(m+n)!}\sum_{i=0}^n (-1)^i \pmatrix{m+n\cr i}
(n-i)^{n+m-1}\\
&=&\frac{m}{(n+m)!}A(m+n-1,n)
\end{eqnarray*}
as required.
\end{pf*}

\begin{pf*}{Second proof of Lemma \ref{Lemma: transition probabilities}}
Consider the birth process $W(t)$ defined by $W(0) = 1$ and
\[
W(t) = \min\{z \in\Z_+ \dvtx T_z > t \} \qquad  (t >0),
\]
where $T_z$ is defined as in (\ref{E: times as sums of exponential rvs}).
The inter-arrival times of $W(t)$ are $(i \xi_i)_{i=1}^z$ and, given
$U(0) = z$,
the death process $U(t)$ has the same inter-arrival
times but taken in the reverse order. The processes
$V(t)$ and $W(t)$ are independent and identically distributed.
Define for $n , m \in\N$,
\[
r(n,m) = \P\bigl( \exists t> 0 \dvtx W(t) = n, V(t) = m \midd V(0) = W(0) = 1 \bigr).
\]

If $Z_k = n$, then $Z_{k+1}$ is the value of $V$ when $W$
first reaches the value \mbox{$n+1$}; $Z_{k+1} = m$ if and only if the
process $(W,V)$ reaches $(n,m)$ \textit{and} then makes the transition to
$(n+1,m)$. Since\vadjust{\goodbreak} $(W,V)$ is Markov, this occurs with
probability $r(n,m) \frac{m}{n+m}$. So for $n,m \in\N$,
%
%e12 ###
\begin{equation}\label{E: p and r}
p(n,m) = \frac{m}{n+m} r(n,m) .
\end{equation}
Conditioning on the site from which $(W, V)$ jumps to $(n,m)$,
we get, for $n,m \in\N$, $n+m \geq3$,
%
%e13 ###
\begin{equation}\label{E: recurrence for r}
r(n,m) = \frac{m}{n+m-1} r(n-1,m) + \frac{n}{n+m-1} r(n,m-1),
\end{equation}
where $r(0,m)= r(n,0) = 0$.
It is easy to check that $r(k,1) = r(1,k) = 1/k!$ for all $k \in\N$.
It will be helpful to define
\[
s(n,m) = (n+m-1)! r(n,m) .
\]
Then we have for $n, m \in\N$, $n+m \geq3$,
\[
s(n,m) = m s(n-1,m) + n s(n,m-1) ,
\]
\[
s(k,1) = s(1,k) = 1 \qquad \mbox{for all }k \in\N.
\]
These constraints completely determine the \textit{positive integers}
$s(n,m)$ for all $m,n \in\N$.
Since the Eulerian numbers $A(n+m-1,m)$ satisfy the same initial
conditions and recurrence relation (\cite{Bona}, Theorem 1.7), we have
$s(n,m) = A(m+n -1, m)$,
which together with (\ref{E: p and r}) gives the desired formula for $p(n,m)$.
\end{pf*}

It is evident from (\ref{E: recurrence for r}) and its initial
conditions that $r(n,m) = r(m,n)$ for all $n,m \in\N$. So
%
%e14 ###
\begin{equation}\label{E: p and Z}
n p(n,m) = m p(m,n) .
\end{equation}
Therefore, the $\sigma$-finite measure $\pi$ on $\mathbb{N}$ defined
by $\pi(n) = n$
satisfies the detailed balance equations and hence is invariant for
$p(\cdot,\cdot)$.
In fact there is a pathwise relation of the same type, which we now describe.
We call a sequence
$\omega=(x_j,y_j)_{j=0}^k$ ($k \geq2$)
of points in $\Z_+^2$ an
\textit{admissible traversal} if
$y_0 = x_k =0$, $x_0 \geq1$, $y_k \geq1$,
each
point $(x_j,y_j)$, $2 \leq j \leq k-1$, is
one of $(x_{j-1}-1,y_{j-1})$, $(x_{j-1},y_{j-1}+1)$, and $(x_1,y_1) =
(x_0,y_0 +1)$,
$(x_k, y_k) = (x_{k-1} - 1,y_{k-1})$. If $\omega$ is an admissible
traversal, then so is the time-reversed and reflected path
$\omega' = (y_{k-j},x_{k-j})_{j=0}^k$. In fact, conditioning on the endpoints,
$\omega$ and $\omega'$ have the same probability of being realized by
the simple harmonic urn.

\begin{lemma}
\label{lem:reverse}
For any admissible traversal
$(x_j,y_j)_{j=0}^k$ with $x_0 = n \in\N$, $y_k =m \in\N$,
\begin{eqnarray*}
&&\P\bigl( ( X_j, Y_j )_{j=0}^{\nu_1} =
(x_j,y_j)_{j=0}^k \midd Z_0 = n, Z_1 = m \bigr) \\
&&\qquad =
\P\bigl( ( X_j ,Y_j )_{j=0}^{\nu_1} =
(y_{k-j},x_{k-j})_{j=0}^k \midd Z_0 = m, Z_1 = n \bigr) .
\end{eqnarray*}
\end{lemma}

\begin{pf}
Let $\omega=(x_j,y_j)_{j=0}^k$ be an admissible traversal,
and define
\[
p = p(\omega) = \P\bigl( ( X_j, Y_j )_{j=0}^{\nu_1} =
(x_j,y_j)_{j=0}^k , Z_1 = m \midd Z_0 = n \bigr)
\]
and
\[
p' = p'(\omega) = \P\bigl( ( X_j, Y_j )_{j=0}^{\nu_1} =
(y_{k-j},x_{k-j})_{j=0}^k , Z_1 = n \midd Z_0 = m \bigr) ,
\]
so that
$p'$ is the probability of the reflected and time-reversed path.
To prove the lemma, it suffices to show that for any $\omega$ with
$(x_0,y_0) = (n,0)$ and $(x_k,y_k)=(0,m)$, $p(\omega) /p(n,m)=
p'(\omega)/p(m,n)$. In light of (\ref{E: p and Z}), it therefore
suffices to show that
$n p = m p'$.
To see this, we use the Markov property along the path $\omega$ to obtain
\[
p = \prod_{j=0}^{k-1} (x_j + y_j )^{-1} \bigl( x_j \1_{\{ x_{j+1} = x_j \}}
+ y_j \1_{\{ y_{j+1} = y_j \}} \bigr) ,
\]
while, using the Markov property along the reflection and reversal of
$\omega$,
\begin{eqnarray*}
p' & = &\prod_{j=0}^{k-1} (x_{k-j}+y_{k-j})^{-1}
\bigl( x_{k-j} \1_{\{ x_{{k-j-1}} = x_{k-j} \}}
+ y_{k-j} \1_{\{ y_{k-j-1} = y_{k-j} \}} \bigr) \\
& = &\prod_{i=0}^{k-1} (x_{i+1}+y_{i+1})^{-1}
\bigl( x_{i} \1_{\{ x_{{i+1}} = x_{i} \}}
+ y_i \1_{\{ y_{i+1} = y_{i} \}} \bigr) ,
\end{eqnarray*}
making the change of variable $i = k-j-1$. Dividing the two products
for $p$ and $p'$ yields, after
cancellation, $p/p' = (x_k+y_k)/ (x_0 +y_0 ) = m/n$,
as required.
\end{pf}

\begin{remarks*}
Of course by summing over paths in the equality $n p(\omega) = m
p'(\omega)$,
we could use the argument in the last proof to prove (\ref{E: p and Z}).
The reversibility and the invariant measure exhibited in Lemma \ref
{lem:reverse}
and (\ref{E: p and Z})
will appear naturally in terms of a stationary model in Section \ref
{S: stationary model}.
\end{remarks*}

%s4 ###
\section{\texorpdfstring{Proof of Theorem \protect\ref{th:zinc} via the uniform renewal process}
{Proof of Theorem 2.2 via the uniform renewal process}}
\label{S: drift}

In this section, we study the
asymptotic behavior of $\E[ Z_{k+1} \midd Z_k = n]$
as $n \to\infty$.
The explicit expression for the distribution of $Z_{k+1}$
given $Z_k =n$ obtained in Lemma
\ref{Lemma: transition probabilities}
turns out not to be very convenient to use directly.
Thus, we proceed somewhat indirectly and exploit
a connection with a renewal process
whose inter-arrival times are uniform on $(0,1)$.
Here and subsequently, we use $U(0,1)$ to denote the uniform
distribution on $(0,1)$.

Let $\chi_1, \chi_2, \chi_3, \ldots$ be an i.i.d. sequence of
$U(0,1)$ random variables. Consider the renewal sequence $S_i$, $i \in
\Z_+$,
defined by $S_0 := 0$ and, for $i \geq1$, $S_i := \sum_{j=1}^i \chi_j$.
For $t \geq0$, define the counting process
%
%e15 ###
\begin{equation}
\label{eq:renewproc}
N(t) := \min\{ i \in\Z_+ \dvtx S_i > t \} = 1+\max\{i \in\Z_+ \dvtx S_i
\le t\} ,
\end{equation}
so a.s., $N(t) \geq t+1$.
In the language of classical renewal theory,
$\E[N(t)]$ is a renewal function (note that we
are counting the renewal at time $0$).
The next result establishes the connection between the uniform renewal
process and the simple harmonic urn.

\begin{lemma}\label{Lemma: same distribution}
For each $n \in\N$, the conditional distribution of $Z_{k+1}$ on $Z_k
= n$ equals the distribution of $N(n)-n$.
In particular, for $n \in\N$, $\E[ Z_{k+1} \midd Z_k = n] = \E[N(n) ]-n$.
\end{lemma}

The proof of Lemma \ref{Lemma: same distribution} amounts to
showing that $\P(N(n) = n+m ) = p (n,m)$ as given by Lemma
\ref{Lemma: transition probabilities}. This equality is Theorem 3
in \cite{sparac}, and it may be verified combinatorially
using the interpretation of $A(n,k)$ as the number of permutations of
$\{1,\dots,n\}$ with exactly $k-1$ falls, together with the
observation that for $n \in\N$,
$N(n)$ is the position of the $n$th fall in
the sequence $\psi_1,\psi_2,\dots,$ where $\psi_k= S_k \bmod1$,
another sequence of i.i.d.\ $U(0,1)$
random variables. Here, we will give a neat proof of Lemma \ref{Lemma:
same distribution} using the coupling
exhibited above in Section \ref{secRubin}.

\begin{pf*}{Proof of Lemma \ref{Lemma: same distribution}}
Consider a doubly-infinite sequence $(\xi_i)_{i \in\mathbb{Z}}$
of independent exponential random variables with mean $1$. Taking
$\zeta_k = \xi_{-k}$, we can write $R_{n,m}$ (as defined
in the first proof of Lemma \ref{Lemma: transition probabilities})
as $\sum_{i= -m}^n i
\xi_i$. Define $S_{n,m} = \sum_{i = -m}^n \xi_i$. For fixed $n
\in\N$, $m \in\Z_+$, we consider normalized partial sums
\[
\chi_j' = \Biggl(\sum_{i = -m}^{j-1-m} \xi_i \Biggr) \Big/ S_{n,m} ,\qquad  j \in\{ 1,
\ldots, n+m \}.
\]
Since $(S_{j-1-m,m})_{j=1}^{n+m}$
are the first $n+m$ points of a unit-rate Poisson process on
$\R_+$, the vector $( \chi_1' , \chi_2', \ldots, \chi_{n+m}') $ is
distributed as the vector of increasing order statistics of the
$n+m$ i.i.d.\ $U(0,1)$ random variables $\chi_1,\ldots,\chi_{n+m}$. In
particular,
\[
\P\bigl(N(n) > n+m \bigr) = \P\Biggl(\sum_{i=1}^{n+m} \chi_i \leq n \Biggr)
= \P\Biggl(\sum_{i=1}^{n+m} \chi_i' \leq n \Biggr) ,
\]
using the fact that, by (\ref{eq:renewproc}),
$\{ N(n) > r\} = \{ S_r \leq n \}$ for $r \in\Z_+$ and $n > 0$.
But
\[
n -\sum_{i=1}^{n+m} \chi_i' = \sum_{i=m+1}^{m+n} (1 - \chi_i') -
\sum_{i=1}^m \chi_i'
= \Biggl(\sum_{i =-m}^n i \xi_i \Biggr)\Big/ S_{n,m}
= R_{n,m} / S_{n,m} .
\]
So, using the equation two lines above (\ref{eq:cdfofZ}),
\[
\P\bigl(N(n) - n > m \bigr) = \P(R_{n,m} \geq0 ) = \P(Z_{k+1} > m \mid Z_k = n ) .
\]
Thus, $N(n)-n$ has the same distribution as $Z_{k+1}$ conditional\vadjust{\goodbreak} on $Z_k=n$.~%
\end{pf*}

In view of Lemma \ref{Lemma: same distribution},
to study $\E[Z_{k+1} \midd Z_k = n]$
we need
to study $\E[N(n)]$.

\begin{lemma}\label{lem:conv} As $n \to\infty$,
\[
\E[N(n)] - \bigl(2n+\tfrac23 \bigr) \to0.
\]
\end{lemma}

\begin{pf}
This is a consequence of the renewal theorem.
For a general nonarithmetic renewal process whose inter-arrival times have
mean $\mu$ and variance $\sigma^2$, let $U(t)$ be the expectation of
the number of arrivals up to
time~$t$, including the initial arrival at time $0$. Then
%
%e16 ###
\begin{equation}
\label{renewalthm}
U(t) - \frac{t}{\mu} \to
\frac{\sigma^2 + \mu^2}{2 \mu^2}  \qquad \mbox{as } t \to\infty.
\end{equation}
We believe this is due to Smith \cite{Smith}. See, for example, Feller
\cite{Feller2}, Section~XI.3, Theorem 1, Cox \cite{Cox}, Section 4, or
Asmussen \cite{Asmussen}, Section V, Proposition~6.1. When the
inter-arrival distribution is $U(0,1)$, we have $U(t) =
\mathbb{E}[N(t)]$ with the notation of (\ref{eq:renewproc}), and in
this case $\mu= 1/2$ and $\sigma^2 = 1/12$.
\end{pf}

Together with Lemma \ref{Lemma: same distribution},
Lemma \ref{lem:conv} gives the following result.

\begin{corollary}\label{cor:Ezk+1}
$\E[Z_{k+1}\midd Z_k=n]-n\to\frac23$ as $n\to\infty$.
\end{corollary}

To obtain the exponential error bound in
(\ref{zinc}) above, we need to know more about the rate of
convergence in Corollary \ref{cor:Ezk+1} and hence in
Lemma~\ref{lem:conv}. The existence of a bound like (\ref{zinc})
for \textit{some} $\alpha_1 < 0$ follows from known results:
Stone \cite{Stone} gave an exponentially small error bound in the
renewal theorem~(\ref{renewalthm}) for inter-arrival distributions
with exponentially decaying tails, and
an exponential bound also
follows from the coupling proof of the renewal theorem (see, e.g.,
Asmussen \cite{Asmussen}, Section VII, Theorem 2.10 and Problem~2.2). However, in this particular case
we can solve the renewal
equation exactly and deduce the asymptotics more precisely, identifying
a (sharp) value for~$\alpha_1$ in (\ref{zinc}). The first step is the
following result.

\begin{lemma}
Let $\chi_1, \chi_2, \ldots$ be an i.i.d. sequence of $U(0,1)$
random variables. For $t \in\R_+$,
\[
\P\Biggl( \sum_{i=1}^k \chi_i \leq t \Biggr) =
\sum_{i=0}^{k \wedge\lfloor t\rfloor} \frac{(t-i)^k (-1)^i}{i!
(k-i)!} ,
\]
and
%
%e17 ###
\begin{eqnarray}\label{eq:robwaters}
\mathbb{E} [N(t) ] = U(t) = \sum_{k=0}^{\infty} \P\Biggl(\sum_{i=1}^k
\chi_i \leq t \Biggr) = \sum_{i=0}^{\lfloor t \rfloor} \frac{(i-t)^i \re
^{t-i}}{i!} .
\end{eqnarray}
\end{lemma}

\begin{pf}
The first formula is classical (see, e.g., \cite{Feller2}, page 27);
according to Feller \cite{Feller1}, page 285, it is due to Lagrange.
The second formula follows from observing (with an empty sum being $0$)
\[
U (t) = \E\sum_{k=0}^\infty\mathbf{1} \Biggl\{ \sum_{i=1}^k \chi_i \leq t
\Biggr\} =
\sum_{k=0}^\infty\P\Biggl( \sum_{i=1}^k \chi_i \leq t \Biggr),
\]
and exchanging the order in the consequent double sum
(which is absolutely convergent).
\end{pf}

We next obtain a more tractable explicit formula for the expression
in~(\ref{eq:robwaters}).
Define for $t \geq0$
\[
f(t) := \sum_{i=0}^{\lfloor t \rfloor} \frac{(i-t)^i \re^{t-i}}{i!} .
\]
It is easy to verify (see also \cite{Asmussen}, page 148)
that $f$ is continuous on $[0,\infty)$ and satisfies
%
%e18 ###
\begin{eqnarray} \label{E: differential-delay}
f(t) &=& \re^t \qquad \hspace*{62.5pt} ( 0 \le t \le1),
\nonumber\\[-8pt]\\[-8pt]
f'(t) &=& f(t) - f(t-1)\qquad   (t \geq1) .\nonumber
\end{eqnarray}

\begin{lemma}\label{lem: asymptotic}
For all $t > 0$,
%
%e19 ###
\begin{equation}\label{E: asymptotic expansion} f(t) = 2t + \frac
{2}{3} + \mathop{\sum_{
\gamma\in\mathbb{C} \dvtx \gamma\neq0,}}_{ \gamma= 1-
\exp(-\gamma)
} \frac{1}{\gamma} \re^{\gamma t} .
\end{equation}
The sum is absolutely convergent, uniformly for $t$ in $(\eps,\infty
)$ for any $\eps> 0$.
\end{lemma}

\begin{pf}
The Laplace transform $\mathcal{L}f(\lambda)$ of $f$ exists
for $\textup{Re}(\lambda) > 0$ since $f(t) = 2t+2/3 + o(1)$ as $t \to
\infty$, by
(\ref{eq:robwaters}) and Lemma
\ref{lem:conv}. Using the differential-delay equation (\ref{E:
differential-delay}), we obtain
\[
\mathcal{L}f(\lambda) = \frac{1}{\lambda- 1 + \re^{-\lambda}} .
\]
The principal part of $\mathcal{L}f$ at $0$ is $\frac{2}{\lambda^2}
+ \frac{2}{3\lambda}$.\vspace*{-1pt}
There are simple poles at the nonzero roots of $\lambda- 1 + \re
^{-\lambda}$,
which occur in complex conjugate pairs $\alpha_n \pm i\beta_n$, where
$\alpha= \alpha_1 > \alpha_2 > \cdots$ and $0 < \beta_1 < \beta_2
< \cdots.$\vspace*{1pt}
In fact, $\alpha_n = -\log(2\pi n) + o(1)$ and $\beta_n = (2n +
\frac{1}{2})\pi+ o(1)$ as $n \to\infty$.
For $\gamma= \alpha_n + i \beta_n$, the absolute value of the term
$\re^{\gamma t}/\gamma$ in the right-hand side of (\ref{E: asymptotic
expansion}) is $1/(|\gamma| |1-\gamma|^t)$, hence the sum converges
absolutely, uniformly on any interval $(\eps, \infty)$, $\eps> 0$.

To establish (\ref{E: asymptotic expansion}), we will
compute the Bromwich integral (inverting the Laplace\vadjust{\goodbreak}
transform), using a carefully chosen sequence of rectangular contours:
\[
f(t) = \lim_{R \to\infty}
\int_{\eps- i R}^{\eps+ i R} \frac{\re^{\lambda t}}{\lambda- 1 +
\exp(-\lambda)} \,d\lambda.
\]
To evaluate this limit for a particular value of $t > 0$, we
will take $\eps= 1/t$ and integrate around a sequence $C_n$
of rectangular contours, with vertices at
$(1/t) \pm(2n - \frac{1}{2})\pi i$ and $-2\log n \pm(2n - \frac
{1}{2})\pi i$.
The integrand along the vertical segment at real part $-2 \log n$\vspace*{1pt}
is bounded by $(1+o(1))/n^2$ and the integrand along the horizontal
segments is bounded by $\re/(2n - \frac{1}{2})\pi$ because the
imaginary parts of $\lambda$ and $\re^{-\lambda}$ have the same
sign there, so $|\lambda- 1 + \re^{-\lambda}| \ge\textup
{Im}(\lambda)$.
It follows that the integrals along these three arcs all tend to
zero as $n \to\infty$. Each pole lies inside all but finitely many of
the contours $C_n$, so the principal value of the Bromwich
integral is the sum of the residues of $\re^{\lambda t}/(\lambda- 1 +
\exp(-\lambda))$.
The residue at $0$ is $2t + 2/3$, and the residue at $\gamma= \alpha
_n + i \beta_n$ is $\re^{\gamma t}/\gamma$.
Thus, we obtain (\ref{E: asymptotic expansion}).
\end{pf}

\begin{pf*}{Proof of Theorem \ref{th:zinc}}
The statement of the theorem follows
from Lemma~\ref{lem: asymptotic}, since by Lemma \ref{Lemma: same
distribution} and (\ref{eq:robwaters})
we have \mbox{$\E[ Z_{k+1} \midd Z_k = n] = f(n) -n$} for $n \in\N$.
\end{pf*}

\begin{remarks*}
According to Feller \cite{Feller2}, Problem 2, page 385,
equation (\ref{eq:robwaters})
``is frequently rediscovered in
queuing theory, but it reveals little about the nature of $U$.'' We
have not
found the formula (\ref{E: asymptotic expansion}) in the literature.
The dominant term in $f(t) - 2t - 2/3$ as $t \to\infty$ is $\re
^{\gamma_1 t}/\gamma_1 + \re^{\overline{\gamma_1}t}/\overline
{\gamma_1}$, that is,
\[
\frac{1}{\alpha_1^2+\beta_1^2} \re^{\alpha_1 t} \bigl(\beta_1\sin
(\beta_1 t) + \alpha_1 \cos(\beta_1 t)\bigr) ,
\]
which changes sign infinitely often. After subtracting this term,
the remainder is $O (\re^{\alpha_2 t} )$. The method that
we have used for analyzing the asymptotic behavior of solutions
to the renewal equation was proposed by A. J. Lotka and was put on
a firm basis by Feller \cite{Feller3}; Laplace
transform inversions of this kind were dealt with by Churchill
\cite{Churchill}.
\end{remarks*}

%s5 ###
\section{\texorpdfstring{Proof of Theorem \protect\ref{th:main1}}{Proof of Theorem 2.1}}\label{S: short proof}

The recurrence relation (\ref{E: recurrence for r}) for $r (n,m)$
permits a direct proof of Theorem
\ref{th:main1} (transience), without appealing to the more general
Theorem \ref{thm:mixed}, via standard
martingale arguments applied\break to~$h (Z_k)$ for a judicious choice of
function $h$. This is
the subject of this section.

Rewriting (\ref{E: recurrence for r}) in terms of $p$
yields the following recurrence relation, which does not seem simple to
prove by
conditioning on a step in the urn model; for $n, m \in\N$, $n + m
\geq3$,
%
%e20 ###
\begin{equation}\label{E: recurrence for p}
\biggl(\frac{n+m}{m} \biggr) p(n,m) = p(n-1,m) +
\biggl(\frac{n}{m-1} \biggr) p(n,m-1) ,
\end{equation}
where if $m=1$ we interpret the right-hand side of (\ref{E: recurrence
for p})
as just $p(n-1,1)$, and where $p(0,m)=p(n,0) = 0$. Note $p(1,1) = 1/2$.
For ease of notation, for any function $F$ we will write $\E_n[F(Z)]$
for $\E[ F(Z_{k+1}) \midd Z_k = n]$,
which, by the Markov property, does not depend on $k$.

\begin{lemma}
\label{L: short}
Let $\alpha_1 \approx- 2.0888$ be as in Theorem \ref{th:zinc}. Then
for $n \geq2$,
\begin{eqnarray*}
\E_n \biggl[\frac{1}{Z} \biggr] & =&
\frac{\E_n[Z] - \E_{n-1}[Z]}{n} = \frac{1}{n} + O (\re^{\alpha_1
n} ) ,\\
\E_n \biggl[\frac{1}{Z^2(Z+1)} \biggr] & =&
\frac{\E_{n-1}[ 1/Z ] - \E_n [ 1/Z] }{n} = \frac{1}{n^2(n-1)} + O
(\re^{\alpha_1 n} ) ,
\end{eqnarray*}
where the asymptotics refer to the limits as $n \to\infty$.
\end{lemma}

\begin{pf}
We use the recurrence relation (\ref{E: recurrence for p}) satisfied
by the transition probabilities of $Z_k$.
First, multiply both sides of (\ref{E: recurrence for p}) by $m$, to
get for $n,m \in\N$, $m +n \geq3$,
\[
(n+m)p(n,m)
%= m p(n-1,m) + \frac{mn}{m-1} p(n,m-1)
= m p(n-1,m) + n p(n,m-1) + \frac{n}{m-1} p(n, m-1) ,
\]
where $p(n,0)=p(0,m)=0$. Summing over $m \in\N$ we obtain for $n
\geq2$,
\[
n + \E_n[Z] = \E_{n-1}[Z] + n + n\E_n [ 1/Z ],
\]
which yields the first equation of the lemma after an application of
(\ref{zinc}).

For the second equation, divide (\ref{E: recurrence for p}) through by
$m$ to get for $n, m \in\N$, $m+ n \geq3$,
\[
\frac{(n+m)}{m^2}p(n,m) = \frac{1}{m} p(n-1,m) + \frac{n}{m(m-1)}
p(n,m-1) .
\]
On summing over $m \in\N$ this gives, for $n \geq2$,
\[
n \E_n [ 1/Z^2 ] + \E_n [ 1/Z ] = \E_{n-1} [1/Z] + n \E_n \biggl[\frac
{1}{(Z+1)Z} \biggr] ,
\]
which gives the second equation when we apply the asymptotic part of
the first equation.
\end{pf}

\begin{pf*}{Proof of Theorem \ref{th:main1}}
Note that $h(x) = \frac{1}{x} - \frac{1}{x^2(x+1)}$ satisfies $h(n) >
0$ for all $n \in\N$
while $h(n) \to0$ as $n \to\infty$. By Lemma \ref{L: short}, we have
\[
\E_n [h(Z)] = \E_n [ 1/Z ]
- \E_n \biggl[ \frac{1}{Z^2 (Z+1) } \biggr] =
\frac{1}{n} - \frac{1}{n^2(n-1)} + O (\re^{\alpha_1 n} ) ,
\]
which is less than $h(n)$ for $n$ sufficiently large. In particular,
$h(Z_k)$ is a~positive supermartingale for $Z_k$ outside a finite set. Hence,
a standard result such as \cite{Asmussen}, Proposition 5.4, page 22,
implies that the Markov chain $(Z_k)$ is transient.
\end{pf*}

%s6 ###
\section{Moment and tail estimates for $Z_{k+1} - Z_k$}
\label{S:zinc}

In order to study the asymptotic behavior
of $(Z_k)_{k \in\Z_+}$, we build on the analysis of Section
\ref{S: drift} to obtain more information
about the increments $Z_{k+1} - Z_k$.
We write $\Delta_k := Z_{k+1} - Z_k$ ($k \in\Z_+)$.
From the relation to the uniform renewal process, by
Lemma \ref{Lemma: same distribution},
we have that
%
%e21 ###
\begin{equation}
\label{eq:coupleDelta}
\P( \Delta_k > x \midd Z_k = n) = \P\bigl( N(n) > 2n + x \bigr)
= \P\Biggl( \sum_{i=1}^{2n+x} \chi_i \leq n \Biggr) ,
\end{equation}
where $\chi_1, \chi_2, \ldots$ are i.i.d. $U(0,1)$ random variables,
using the notation at~(\ref{eq:renewproc}).

Lemma \ref{lem:expotail} below gives a tail bound for $|\Delta_k|$ based
on (\ref{eq:coupleDelta}) and a sharp bound for the moment
generating function of a $U(0,1)$ random variable,
for which we have not been able to find a reference and
which we state first since it may be of interest
in its own right.

\begin{lemma}\label{lem:inequality}
For $\chi$ a $U(0,1)$ variable with moment generating function given
for $\l\in\R$ by
%
%e22 ###
\begin{equation}
\label{phidef}
\phi(\l) = \E[ \re^{\l\chi} ] = \frac{\re^\l- 1}{\l} ,
\end{equation}
we have
\[
\log\phi( -\l) \leq-\frac{\l}{2} + \frac{\l^2}{24}\qquad  ( \l\geq0);\qquad
\log\phi(\l) \leq\frac{\l}{2} + \frac{\l^2}{24}\qquad  ( \l\geq0)
.
\]
\end{lemma}

\begin{pf}
Consider the first of the two stated inequalities. Exponentiating
and multiplying both sides by $\l\re^{\l/2}$, this is equivalent
to
%
%e23 ###
\begin{equation}
\label{shine}
2 \sinh(\l/2) \leq\l\exp( \l^2 /24 )
\end{equation}
for all $\l\geq0$. Inequality (\ref{shine}) is easily verified
since both sides are entire functions with
nonnegative Taylor coefficients and the right-hand
series dominates the left-hand series term by term,
because $6^n n! \leq(2n +1)!$ for all $n \in\N$.
The second stated inequality reduces to (\ref{shine})
also on exponentiating and multiplying through
by $\l\re^{-\l/2}$.
\end{pf}

Now we can state our tail bound for $|\Delta_k|$.
The bound in Lemma \ref{lem:expotail} is a slight improvement
on that provided by Bernstein's inequality in this particular case; the latter
yields a weaker bound with $4x$ instead of $2x$ in the denominator of the
exponential.

\begin{lemma}\label{lem:expotail}
For $n \in\N$ and any integer $x \geq0$, we have
\[
\P( | \Delta_k | >x \midd Z_k=n)\le2 \exp\biggl\{ -\frac
{3x^2}{4n+2x}\biggr \}.
\]
\end{lemma}

\begin{pf}
From (\ref{eq:coupleDelta}) and Markov's
inequality, we obtain for $x \geq0$ and any \mbox{$\l\ge0$},
\begin{eqnarray*}
\P(\Delta_k >x \midd Z_k=n)
&=&\P\Biggl( \exp\Biggl\{ -\l\sum_{i=1}^{2n+x} \chi_i \Biggr\} \geq\re^{-\l n}
\Biggr)\\
&\le&\exp\{ \l n+(2n+x) \log\phi(-\l) \},
\end{eqnarray*}
where $\phi$ is given by (\ref{phidef}).
With $\l=6x/(2n+x)$, the first inequality of Lemma~\ref{lem:inequality} yields
\[
\P(\Delta_k >x \midd Z_k=n) \le\exp\biggl\{ - \frac{x \lambda}{4} \biggr\} =
\exp\biggl\{ -\frac{3x^2}{4n+2x}\biggr\}.
\]
On the other hand,
for $x\in[0,n-1]$, from (\ref{eq:coupleDelta}) and Markov's
inequality once more,
\begin{eqnarray*}
\P(\Delta_k \le-x \midd Z_k=n)
& = & \P\Biggl(\sum_{i=1}^{2n-x} \chi_i > n \Biggr)
= \P\Biggl( \exp\Biggl\{ \l\sum_{i=1}^{2n-x} \chi_i \Biggr\} > \re^{\l n} \Biggr)\\
&\le& \exp\{ -\l n+(2n-x) \log\phi(\l) \}.
\end{eqnarray*}
On setting $\l=6x/(2n-x)$, the second inequality in
Lemma \ref{lem:inequality} yields, for any $x \in[0,n-1]$,
\[
\P(\Delta_k < - x \midd Z_k=n) \le\exp\biggl\{ -\frac{3x^2}{4n-2x} \biggr\} \le
\exp\biggl\{ -\frac{3x^2}{4n+2x} \biggr\},
\]
while $\P( \Delta_k < -n \midd Z_k = n) =0$.
Combining the left and right tail bounds completes the proof.
\end{pf}

Next, from Lemma \ref{lem:expotail},
we obtain the following large deviation and moment bounds
for $\Delta_k$.

\begin{lemma}
\label{cor:Etails}
Suppose that $\eps>0$. Then for some $C < \infty$ and all $n \in\N$,
%
%e24 ###
\begin{equation}
\label{deltatail}
\P\bigl( | \Delta_k | > n^{(1/2)+\eps} \midd Z_k=n \bigr)
\le C \exp\{ -n^{\eps} \}.
\end{equation}
Also for each $r \in\N$, there exists $C(r) < \infty$
such that for any $n \in\N$,
%
%e25 ###
\begin{equation}
\label{deltamoms}
\E[ |\Delta_k|^r \midd Z_k=n ] \le C(r) n^{r/2} .
\end{equation}
\end{lemma}

\begin{pf}
The bound (\ref{deltatail})
is straightforward
from Lemma \ref{lem:expotail}. For $r \in\N$,
\begin{eqnarray}
\label{gammas}
\E[ |\Delta_k|^r \midd Z_k=n ] & \leq&\int_{0}^\infty\P( | \Delta
_k | \geq
\lfloor x^{1/r} \rfloor\midd Z_k = n)\, \ud x\nonumber\hspace*{-30pt}\\[-8pt]\\[-8pt]
& \leq & C \int_{0}^{n^r} \exp\biggl\{ - \frac{x^{2/r}}{2n}\biggr \} \,\ud x
+ C \int_{n^r}^{\infty} \exp\biggl\{ - \frac{x^{1/r}}{2} \biggr\} \,\ud
x\nonumber\hspace*{-30pt}
\end{eqnarray}
for some $C < \infty$, by Lemma \ref{lem:expotail}. With the
substitution $y = x^{1/r}$,
the second integral on the last line of (\ref{gammas})
is seen to\vadjust{\goodbreak} be $O ( n^{r-1} \re^{-n})$ by asymptotics for the
incomplete Gamma function. The first integral
on the last line of (\ref{gammas}), with
the substitution $y = (2n)^{-1} x^{2/r}$, is equal to
\[
\frac{(2n)^{r/2} r}{2} \int_0^{n/2} \re^{-y} y ^{(r/2) -1} \,\ud y
\leq\Gamma(r/2) (2n)^{r/2} r /2.
\]
Combining the last two upper bounds we verify (\ref{deltamoms}).
\end{pf}

The next result gives sharp asymptotics for the first
two moments of
$\Delta_k = Z_{k+1} - Z_k$.

\begin{lemma}
\label{lem:Zincrements}
Let
$\alpha_1 \approx-2.0888$ be as in Theorem \ref{th:zinc}. Then as $n
\to\infty$,
\begin{eqnarray}
\label{zinc1}
\E[ \Delta_k \midd Z_k=n ] & =& \tfrac23+O(\re^{\alpha_1 n}),\\
\label{zinc2}
\E[ \Delta_k ^2 \midd Z_k=n ] & =& \tfrac23 n+ \tfrac23 + O (n \re
^{\alpha_1 n}) .
\end{eqnarray}
\end{lemma}

\begin{pf}
Equation (\ref{zinc1}) is immediate from (\ref{zinc}).
Now we observe that $J_n := X_n^2 + Y_n^2 -n$ is a martingale. Indeed,
for any $(x,y) \in\Z^2$,
\begin{eqnarray*}
&&\E[ J_{n+1} - J_n \midd(X_n, Y_n ) = (x,y) ]\\
&&\qquad = \frac{|x|}{|x|+|y|} \bigl(2 y \sgn(x) + 1\bigr) + \frac{|y|}{|x|+|y|} \bigl(-2 x
\sgn(y) + 1\bigr) -1 = 0 .
\end{eqnarray*}
Between times $\nu_k$ and $\nu_{k+1}$, the urn takes $Z_k+Z_{k+1}$ steps,
so $\nu_{k+1} - \nu_k = Z_k+Z_{k+1}$. Moreover,
$J_{\nu_k} = Z_k^2 - \nu_k$.
Applying the optional stopping theorem at $\nu_k$
and $\nu_{k+1}$, we have that
\begin{eqnarray*}
J_{\nu_k} &=& Z_k^2 - \nu_k =
\E[ J_{\nu_{k+1}} \midd Z_k ] = \E[ Z_{k+1}^2 -\nu_{k+1} \midd Z_k ]\\
&=& \E[ Z_{k+1}^2 - Z_{k+1} \midd Z_k ] - \nu_k - Z_k.
\end{eqnarray*}
The optional stopping theorem is applicable here since a.s. $J_n \leq C
n^2$ for some $C< \infty$
and all $n$, while there is an exponential tail-bound for $\nu_{k+1}
-\nu_k$
(see Lemma~\ref{lem:time} below).
Rearranging the equation in the last display, it follows that for $n
\in\N$,
%
%e26 ###
\begin{equation} \label{eq:sq}
\E[ Z_{k+1}^2 \midd Z_k=n ] = n^2+n+ \E[ Z_{k+1} \midd Z_k=n ].
\end{equation}
Writing $\Delta_k =Z_{k+1}-Z_k$,
we have that
\[
\E[ \Delta_k^2 \midd Z_k = n ] = \E[ Z_{k+1}^2 \midd Z_k=n ] -2n \E[
Z_{k+1} \midd Z_k=n ] + n^2,
\]
which with (\ref{eq:sq}) and (\ref{zinc}) gives (\ref{zinc2}).
\end{pf}

\begin{remark*}In view of Lemma \ref{Lemma: same distribution}, we
could have
used renewal
theory (e.g., \cite{smith2}) to estimate $\E[ \Delta_k^2 \midd Z_k
=n]$. However,
no result we could find in the literature would yield a bound as sharp
as that in (\ref{zinc2}).
\end{remark*}

Lemma \ref{Lemma: same distribution} with (\ref{zinc1}) and (\ref{zinc2})
implies an ancillary result on the $U(0,1)$ renewal\vadjust{\goodbreak} process.
\begin{corollary}
\label{uniform}
Let $N(t)$ be the counting function of
the uniform renewal process, as defined by (\ref{eq:renewproc}). Then with
$\alpha_1 \approx-2.0888$ as in Theorem~\ref{th:zinc}, as $t \to
\infty$,
\[
\E[ N(t)^2 ] = 4t^2 + \tfrac{10}{3} t + \tfrac{2}{3} + O (t \re
^{\alpha_1 t} );\qquad
\Var[ N(t) ] = \tfrac{2}{3} t + \tfrac{2}{9} + O (t \re^{\alpha_1 t}
) .
\]
\end{corollary}

These asymptotic results are both sharper than any we have seen
in the literature; see, for example, \cite{jensen,sparac} in the particular
case of the uniform renewal process or \cite{smith2} for the
general case. We remark that the formula given
in \cite{sparac}, page 231, for $\E[N(t)^2]$ contains an error
(in \cite{sparac} the renewal at $0$ is not counted, so the
notation $m_k(\cdot)$ there is equivalent to our $\E[ (N(\cdot)
-1)^k]$).

%s7 ###
\section{Asymptotic analysis of the noisy urn}
\label{sec:tildez}

%s7.1 ###
\subsection{Connection to Lamperti's problem}
\label{sec:lamperti}

In this section, we study the noisy urn model described in Section \ref
{S:noisy}.
To study the asymptotic behavior of $(\tZ_k)_{k \in\Z_+}$, it
turns out to be more convenient to work with\vspace*{-1.5pt}
the process $(W_k)_{k \in\Z_+}$ defined by $W_k = \tZ_k^{1/2}$,
since the latter
process has asymptotically-zero
drift, in a sense to be made precise shortly,
and such processes have been well-studied
in the literature.

Let $(W_k)_{k \in\Z_+}$ be an irreducible time-homogeneous Markov
chain whose
state-space is an unbounded
countable subset of $\R_+$. Define the increment moment functions
%
%e27 ###
\begin{equation}
\label{mudef}
\mu_r(x) := \E[ (W_{k+1} - W_k)^r \midd W_k = x ];
\end{equation}
by the Markov property, when the corresponding moments exist
the $\mu_r(x)$ are genuine functions of $x$. Given a reasonable choice
of scale for the process~$W_k$, it is common that $\mu_2(x)$ be uniformly
bounded away from $0$ and $\infty$. In this case, under some mild additional
regularity conditions, the regime where $x| \mu_1 (x)| = O(1)$ is critical
from the point of view of the recurrence classification of $W_k$.
For a nearest-neighbor random walk on $\Z_+$ this fact had been
known for a long time (see \cite{harris}), but a study of this and
many other aspects of the problem, in much greater generality
(with absence of the Markovian and countable state-space
assumptions), was carried out by Lamperti \cite{lamp1,lamp2,lamp3}
using martingale techniques. Thus, the analysis of processes with
asymptotically zero drift [i.e., $\mu_1(x) \to0$] is sometimes
known as \textit{Lamperti's problem}.

We will next state some consequences
of Lamperti's results that we will use.
For convenience, we impose conditions that are stronger than
Lamperti's.
We suppose that for each $r \in\N$,
%
%e28 ###
\begin{equation}
\label{mubound}
\sup_{x } | \mu_r(x) | < \infty.
\end{equation}

The recurrence and transience properties of $W_k$
were studied by Lamperti \cite{lamp1,lamp3} and his results were
refined by Menshikov, Asymont and Iasnogorodskii \cite{mai}. Parts (i)
and (ii) of the following\vadjust{\goodbreak}
result are consequences of Theorems 3.1 and 3.2 of \cite{lamp1} with
Theorem 2.1 of \cite{lamp3},
while
part (iii) is a~consequence of
Theorem 3 of \cite{mai} (which is in fact a much sharper result).

\begin{prop}[(\cite{lamp1,lamp3,mai})]
\label{prop:lampclass}
Let $( W_k )$ be an irreducible Markov chain on a countable unbounded
subset of $\R_+$.
Suppose that
(\ref{mubound}) holds, and that there exists $v > 0$ such that
$\mu_2(x) > v$ for all $x$ sufficiently large.
Then the following recurrence criteria are valid:
\begin{itemize}[(iii)]
\item[(i)] $W_k$ is transient if there exist $\delta, x_0 \in
(0,\infty)$ such that
for all $x>x_0$,
\[
2x \mu_1 (x) - \mu_2 (x) > \delta.
\]
\item[(ii)] $W_k$ is positive-recurrent if there exist $\delta, x_0
\in(0,\infty)$ such that
for all $x>x_0$,
\[
2x \mu_1 (x) + \mu_2(x) < - \delta.
\]
\item[(iii)] $W_k$ is null-recurrent if there exists $x_0 \in
(0,\infty)$ such that for all $x>x_0$,
\[
2x | \mu_1 (x) | \leq\biggl( 1 + \frac{1}{\log x} \biggr) \mu_2 (x) .
\]
\end{itemize}
\end{prop}

In \cite{lamp2},
Lamperti proved the existence of weak-sense limiting diffusions for certain
processes satisfying parts (i) or (iii) of Proposition \ref{prop:lampclass}.
To state Lamperti's result,
we need some more notation.
To describe the time-homogeneous
diffusions on $\R_+$
that arise here, it will suffice to describe the
infinitesimal mean $\mu(x)$ and infinitesimal variance
$\sigma^2 (x)$; see, for example, \cite{kt2}, Chapter~15. The transition
functions $p$
of our diffusions will then satisfy the Kolmogorov backward
equation
\[
\frac{ \partial p}{\partial t} = \mu(x) \frac{\partial p}{\partial x}
+\frac{1}{2} \sigma^2 (x) \frac{ \partial^2 p}{\partial x^2} .
\]
Let $(H^{\alpha,\beta}_t)_{t \in[0,1]}$
denote a diffusion process
on $\R_+$ with
infinitesimal mean and variance
%
%e29 ###
\begin{equation}
\label{diffdef}
\mu(x) = \frac{\alpha}{x},\qquad  \sigma^2 (x) = \beta.
\end{equation}

The particular case of a diffusion satisfying (\ref{diffdef}) with
$\beta=1$ and
$\alpha= (\gamma-1)/2$ for some $\gamma\in\R$ is a Bessel process
with parameter $\gamma$; in this case we use the
notation $V^\gamma_t = H^{(\gamma-1)/2,1}_t$.
Recall that for
$\gamma\in\N$, the law of $(V^\gamma_t)_{t \in[0,1]}$
is the same as that of $\| B_t \|_{t \in[0,1]}$ where $(B_t)_{t \in[0,1]}$
is standard $\gamma$-dimensional Brownian motion. In fact, any
$H^{\alpha,\beta}_t$
is related to a Bessel process via simple scaling, as the next result shows.

\begin{lemma}
\label{lem:bessel}
Let $\alpha\in\R$ and $\beta>0$.
The diffusion process $H_t^{\alpha,\beta}$ is a scaled Bessel process:
\[
(H_t^{\alpha,\beta})_{t \in[0,1]}
\mbox{ has the same law as }
( \beta^{1/2} V_t^\gamma)_{t \in[0,1]} \qquad \mbox{with } \gamma= 1 +
\frac{2\alpha}{\beta} .\vadjust{\goodbreak}
\]
\end{lemma}

\begin{pf}
By the It\^o transformation formula
(cf. page 173 of \cite{kt2}), for any $\beta>0$
the process $( \beta^{1/2} V^\gamma_t)_{t \in[0,1]}$
is a diffusion process on $[0,1]$ with infinitesimal mean $\mu(x) =
\beta(\gamma-1)/(2x)$
and infinitesimal variance $\sigma(x) = \beta$,
from which we obtain the result.
\end{pf}

We need the following form
of Lamperti's invariance principle (\cite{lamp2}, Theorems~2.1, 5.1 and A.2).

\begin{prop}[(\cite{lamp2})]
\label{prop:lampinv}
Let $(W_k)$ be an irreducible Markov chain on a~countable unbounded
subset of $\R_+$.
Suppose that
(\ref{mubound}) holds, and that
\[
\lim_{x \to\infty} \mu_2 (x) = \beta>0,\qquad  \lim_{x \to\infty} x
\mu_1 (x) = \alpha> - (\beta/2) .
\]
Let $(H^{\alpha,\beta}_t)_{t\geq0}$ be a diffusion process
as defined at (\ref{diffdef}). Then as $k \to\infty$,
\[
( k^{-1/2} W_{kt} )_{t \in[0,1]} \to( H^{\alpha,\beta}_t )_{t \in
[0,1]}
\]
in the sense of
convergence of finite-dimensional distributions.
Marginally,
\begin{eqnarray*}
\lim_{k \to\infty} \P( k^{-1/2} W_k \leq y ) &=& \frac{2}{(2 \beta
)^{(\alpha/\beta)+(1/2)} \Gamma((\alpha/\beta) + (1/2)) }\\
&&{}\times\int_0^y r^{2\alpha/\beta} \exp\bigl( - r^2/(2 \beta) \bigr) \,\ud r .
\end{eqnarray*}
\end{prop}

%s7.2 ###
\subsection{Increment moment estimates for $W_k$}
\label{sec:increments}

Now consider the\vspace*{-1pt} process $(W_k)_{k \in\Z_+}$ where $W_k = \tZ_k^{1/2}$;
this is a Markov chain with a countable state space (since
$\tZ_k$ is), so fits into the framework described in Section \ref
{sec:lamperti}
above. Lemma \ref{lem:sqrtincrements} below
shows that indeed $W_k$ is an instance of Lamperti's problem in the critical
regime. First we need some simple properties of the random variable
$\kappa$.

\begin{lemma}
\label{kappalem}
If $\kappa$ satisfies (\ref{kappabound}) for $\lambda>0$, then
\begin{equation}
\label{kappatail} \P( | \kappa| \geq x) \leq\exp\{ - \lambda x
\}\qquad (x \geq0)
\end{equation}
and
\begin{equation}
\label{kappamoms}
\E[ | \kappa|^r ] < \infty\qquad  ( r \geq0).
\end{equation}
\end{lemma}

\begin{pf}
(\ref{kappatail}) is immediate from Markov's inequality and (\ref
{kappabound}),
and (\ref{kappamoms}) is also straightforward.
\end{pf}

Now we can start our analysis of the noisy urn and the associated process~$\tilde Z_k$.
Recall that $\tZ_k$ is defined as $\max\{ |\tilde X_{\tilde\nu_k
+1}| , |\tilde Y_{\tilde\nu_k +1 } | \}$.
By definition of the noisy urn process, if we start at unit distance
away from an
axis (in the anticlockwise sense), the path of the noisy urn until it
hits the next axis
has the
same distribution as the corresponding path in the original simple
harmonic urn. Since we refer
to this fact often, we state it as a lemma.

\begin{lemma}
\label{noisylem}
Given $\tilde Z_k = z$, the path $(\tilde X_n, \tilde Y_n)$ for $n =
\tilde\nu_k +1,\ldots, \tilde\nu_{k+1}$
has the same distribution as the
path $(X_n, Y_n)$ for $n = \nu_k +1,\ldots, \nu_{k+1}$ given \mbox{$Z_k = z$}.
In particular,
$\tZ_{k+1}$ conditioned
on $\tZ_k = z$ has the same
distribution as $Z_{k+1} - \min\{ \kappa, Z_{k+1} - 1 \}=Z_{k+1} -
\kappa+ (\kappa+1-Z_{k+1}) \1 \{ \kappa\ge Z_{k+1} \}$ conditioned
on \mbox{$Z_k = z$}.
\end{lemma}

Recall that $\Delta_k = Z_{k+1} - Z_k$, and write $\tD_k = \tZ_{k+1}
- \tZ_k$. The next result
is an analogue of Lemmas \ref{cor:Etails} and \ref{lem:Zincrements}
for $\tD_k$.

\begin{lemma}
\label{lem:tdelta}
Suppose that (\ref{kappabound}) holds.
Let $\eps>0$. Then for some $C < \infty$ and all $n \in\N$,
%
%e30 ###
\begin{equation}
\label{Anbound}
\P\bigl( | \tD_k | > n^{(1/2)+\eps} \midd\tZ_k = n \bigr) \leq C \exp\{ -
n^{\eps/3} \} .
\end{equation}
Also, for any $r \in\N$, there exists $C <\infty$ such that for any
$n \in\N$,
%
%e31 ###
\begin{equation}
\label{tdeltamoms} \E[ | \tD_k |^r \midd\tZ_k = n ] \leq C
n^{r/2} .
\end{equation}
Moreover, there exists $\gamma>0$ for
which, as $n \to\infty$,
\begin{eqnarray}
\label{tzinc1}
\E[ \tD_k \midd\tZ_k = n ]
& =& \tfrac{2}{3} - \E[\kappa] + O (\re^{ - \gamma n}) ,\\
\label{tzinc2}
\E[ \tD_k^2 \midd\tZ_k = n ] & =& \tfrac{2}{3} n + O(1).
\end{eqnarray}
\end{lemma}

\begin{pf}
By the final statement in Lemma \ref{noisylem},
for any $r \geq0$,
\[
\P( | \tilde\Delta_k | > r \midd\tZ_k = n )
\leq
\P( | \Delta_k - \kappa| >r \midd Z_k = n ) .
\]
We have for any $\eps>0$,
\begin{eqnarray*}
&&\P\bigl( | \Delta_k - \kappa| > n^{(1/2)+\eps} \midd Z_k =
n \bigr)\\
&&\qquad \leq\P\bigl( | \Delta_k | > n^{(1+\eps)/2} \midd Z_k = n\bigr)
+ \P\bigl( | \kappa| > n^{(1+\eps)/2} \bigr)
\end{eqnarray*}
for all $n$ large enough. Using the bounds in (\ref{deltatail}) and
(\ref{kappatail}), we obtain (\ref{Anbound}).
For $r \in\N$,
\[
\E[ | \tD_k |^r \midd\tZ_k = n ] \leq\E[ ( | \Delta_k | + |\kappa
| )^r \midd Z_k = n ] .
\]
Then with Minkowski's inequality, (\ref{deltamoms})
and (\ref{kappamoms}) we obtain (\ref{tdeltamoms}).

Next, we have from Lemma \ref{noisylem} and (\ref{zinc1}) that
\begin{eqnarray*}
&&\E[ \tD_k \midd\tZ_k = n ]\\
 &&\qquad  =  \E[ \Delta_k -\kappa
+ (\kappa+1-Z_{k+1}) \1 \{ \kappa\ge Z_{k+1}\} \midd Z_k = n ] \\
&&\qquad  = \tfrac{2}{3} + O (\re^{\alpha_1 n}) -\E[\kappa] + \E[
(\kappa+1-Z_{k+1}) \1 \{ \kappa\ge Z_{k+1}\} \midd Z_k = n ].
\end{eqnarray*}
By the Cauchy--Schwarz inequality, (\ref{kappamoms}) and the
bound $0\le\kappa+1-Z_{k+1}\le\kappa$, the last term here is
bounded by a constant times the square-root of
\begin{eqnarray*}
\P( \kappa\geq Z_{k+1} \midd Z_k =n ) &\leq&\P( | \Delta_k | \geq n/2
\midd Z_k = n ) + \P( | \kappa| > n/2 )\\
&=& O ( \exp\{ - \lambda n /2 \} ),
\end{eqnarray*}
using the bounds (\ref{deltatail}) and (\ref{kappatail}). Hence, we
obtain (\ref{tzinc1}).
Similarly, from~(\ref{zinc2}), we obtain (\ref{tzinc2}).
\end{pf}

Now, we can give the main result of this section on the increments
of the process $(W_k)_{k \in\Z_+}$.

\begin{lemma}
\label{lem:sqrtincrements}
Suppose that $\kappa$ satisfies (\ref{kappabound}).
With $\mu_r(x)$ as defined by (\ref{mudef}), we have
%
%e32 ###
\begin{equation}
\label{wlem1}
\sup_x | \mu_r (x)| < \infty
\end{equation}
for each $r \in\N$. Moreover as $x \to\infty$,
%
%e33 ###
\begin{equation}
\label{wlem2}
\mu_1 (x) = \frac{1-2\E[\kappa] }{4x} + O(x^{-2} );\qquad  \mu_2 (x) =
\frac{1}{6} + O(x^{-1}) .
\end{equation}
\end{lemma}

\begin{pf}
For the duration of this proof, we write $\E_{x^2} [ \cdot]$
for $\E[ \cdot\midd\tZ_k = x^2 ] = \E[ \cdot\midd W_k = x ]$.
For $r \in\N$ and $x \geq0$, from (\ref{mudef}),
\begin{equation}
\label{murbound}
| \mu_r(x) | \leq\E_{x^2} [ | \tZ_{k+1}^{1/2} - \tZ_k^{1/2} |^r ]
= x^r \E_{x^2} [ | ( 1+ x^{-2} \tD_k )^{1/2} - 1 |^r ].
\end{equation}
Fix $\eps>0$ and write
$A(n) := \{ | \tD_k | > n^{(1/2)+\eps}\}$ and $A^c(n)$ for the
complementary event.
Now for some $C < \infty$ and all $x \geq1$, by Taylor's theorem,
\[
| ( 1+ x^{-2} \tD_k )^{1/2} - 1 |^r \1_{A^c(x^2)}
\leq C x^{-2r} | \tD_k |^r \1_{A^c(x^2)} .
\]
Hence,
\begin{equation}
\label{murbound1}
\quad\E_{x^2} \bigl[ | ( 1+ x^{-2} \tD_k )^{1/2} - 1 |^r \1_{A^c(x^2)} \bigr]
\leq C x^{-2r} \E_{x^2} [ | \tD_k |^r ] = O( x^{-r} )
\end{equation}
by (\ref{tdeltamoms}).
On the other hand, using the fact that for $y \geq-1$,
$0 \leq(1+y)^{1/2} \leq1 + (y/2)$, we have
\begin{eqnarray*}
&&\E_{x^2} \bigl[ | ( 1+ x^{-2} \tD_k )^{1/2} - 1 |^r \1_{A(x^2)} \bigr]\\
&&\qquad \leq\E_{x^2} \bigl[ \bigl(1 + (1/2) x^{-2} | \tD_k | \bigr)^r \1_{A(x^2)} \bigr]
\\
&&\qquad \leq\bigl( \E_{x^2} [ (1 + | \tD_k |)^{2r} ] \bigr)^{1/2}
\bigl( \P\bigl( A(x^2) \midd\tZ_k = x^2 \bigr) \bigr)^{1/2}
\end{eqnarray*}
for $x \geq1$,
by Cauchy--Schwarz. Using (\ref{tdeltamoms})
to bound the expectation
here and (\ref{Anbound})
to bound the probability, we obtain, for any $r \in\N$,
%
%e34 ###
\begin{equation}
\label{murbound2}
\E_{x^2} \bigl[ | ( 1+ x^{-2} \tD_k )^{1/2} - 1 |^r \1_{A(x^2)} \bigr]
= O ( \exp\{ - x^{\eps/2} \} ).
\end{equation}
Combining (\ref{murbound1}) and (\ref{murbound2}) with (\ref{murbound}),
we obtain (\ref{wlem1}).

Now, we prove (\ref{wlem2}).
We have that for $x \geq0$,
\begin{eqnarray}
\label{zmu0}
\mu_1 (x) & =& \E_{x^2} [ W_{k+1} - W_k ]
= x \E_{x^2} [ (1 + x^{-2} \tD_k )^{1/2} -1 ]\nonumber\\[-8pt]\\[-8pt]
& =& x \E_{x^2} \bigl[ \bigl( (1 + x^{-2} \tD_k )^{1/2} -1 \bigr) \1_{A^c (x^2)} \bigr] +
O ( \exp\{ - x^{\eps/3} \} ) ,\nonumber
\end{eqnarray}
using (\ref{murbound2}).
By Taylor's theorem with Lagrange form for the remainder, we have
\begin{eqnarray}
\label{zmu1}
&&\quad x \E_{x^2} \bigl[ \bigl( (1 + x^{-2} \tD_k )^{1/2} -1 \bigr) \1_{A^c(x^2)} \bigr]
\nonumber\\[-8pt]\\[-8pt]
&&\quad\qquad  = \frac{1}{2x} \E_{x^2} \bigl[ \tD_k \1_{A^c(x^2)} \bigr] -
\frac{1}{8x^3} \E_{x^2} \bigl[ \tD_k^2 \1_{A^c(x^2)} \bigr]
+ O ( x^{-5} \E_{x^2} [ | \tD_k|^3 ] ) .\nonumber
\end{eqnarray}
Here we have that
$x^{-5} \E_{x^2} [ | \tD_k|^3 ] = O( x^{-2} )$,
by (\ref{tdeltamoms}),
while for $r \in\N$, we obtain
\[
\E_{x^2} \bigl[ \tD_k^r \1_{A^c(x^2)} \bigr] = \E_{x^2} [ \tD_k^r ]
+ O \bigl( \bigl(\E_{x^2} [ | \tD_k |^{2r} ] \P\bigl( A (x^2) \midd\tZ_k = x^2 \bigr)
\bigr)^{1/2} \bigr)
\]
by Cauchy--Schwarz. Using (\ref{Anbound}) again and combining
(\ref{zmu0}) with (\ref{zmu1}), we obtain
\[
\mu_1 (x) = \frac{1}{2x} \E_{x^2} [ \tD_k ] -
\frac{1}{8x^3} \E_{x^2} [ \tD_k^2 ]
+ O ( x^{-2} ) .
\]
Thus, from (\ref{tzinc1}) and (\ref{tzinc2}),
we obtain
the expression for $\mu_1$ in (\ref{wlem2}).
Now, we use the fact that
\begin{eqnarray*}
(W_{k+1} - W_k)^2 &=& W_{k+1}^2 - W_k^2 - 2 W_k (W_{k+1} - W_k) \\
&=& \tZ_{k+1} - \tZ_k - 2W_k (W_{k+1} - W_k )
\end{eqnarray*}
to obtain
$\mu_2 (x) = \E_{x^2} [ \tD_k ] - 2 x \mu_1 (x)$,
which with (\ref{tzinc1})
yields the expression for $\mu_2$ in (\ref{wlem2}).
\end{pf}

%s8 ###
\section{Proofs of theorems}
\label{sec:proofs}

%s8.1 ###
\subsection{\texorpdfstring{Proofs of Theorems \protect\ref{th:leaky}, \protect\ref{thm:mixed}, \protect\ref
{thm:moments} and \protect\ref{thm:difflimit}}
{Proofs of Theorems 2.3, 2.4, 2.5, and 2.8}}

First, we work with the noisy urn model of Section \ref{S:noisy}.
Given the moment estimates of Lemma \ref{lem:sqrtincrements}, we
can now apply the general results described in Section
\ref{sec:lamperti} and \cite{aim}.

\begin{pf*}{Proof of Theorem \ref{thm:mixed}}
First, observe that $(\tZ_k)_{k \in\Z_+}$ is transient, null-, or
positive-recurrent
exactly when $(W_k)_{k \in\Z_+}$ is.
From Lemma \ref{lem:sqrtincrements}, we have that
\begin{eqnarray*}
2x \mu_1(x) - \mu_2 (x) &=& \tfrac{1}{3} - \E[\kappa]+ O(x^{-1});\\
2 x \mu_1 (x) + \mu_2 (x) &=& \tfrac{2}{3} - \E[\kappa] + O(x^{-1}).
\end{eqnarray*}
Now apply Proposition \ref{prop:lampclass}.
\end{pf*}

\begin{pf*}{Proof of Theorem \ref{thm:moments}}
By the definition of
$\tau_q$ at (\ref{returntime}),
$\tau_q$ is also the first hitting time of $1$
by $(W_k)_{k \in\N}$.
Then with Lemma \ref{lem:sqrtincrements} we can
apply results of Aspandiiarov, Iasnogorodski and Menshikov \cite{aim}, Propositions 1 and 2, which generalize
those of Lamperti \cite{lamp3} and give
conditions on $\mu_1$ and $\mu_2$ for existence and nonexistence of
passage-time moments,
to obtain the stated result.
\end{pf*}

\begin{pf*} {Proof of Theorem \ref{thm:difflimit}}
First, Proposition \ref{prop:lampinv} and Lemma \ref{lem:sqrtincrements}
imply that, as $n \to\infty$,
\[
( n^{-1/2} W_{nt} )_{t \in[0,1]} \to( H^{\alpha,\beta}_t )_{t \in[0,1]}
\]
in the sense of finite-dimensional distributions,
where $\alpha= (1-2\E[\kappa])/4$ and $\beta=1/6$, provided $\E
[\kappa] < 2/3$.
By the It\^o transformation formula
(cf. page~173 of \cite{kt2}), with $H^{\alpha,\beta}_t$ as defined at
(\ref{diffdef}),
$(H^{\alpha,\beta}_t)^2$ is a diffusion process
with infinitesimal mean $\mu(x) = \beta+2\alpha$ and infinitesimal
variance $\sigma^2 (x) = 4 \beta x$. In particular,
$(H^{\alpha,\beta}_t)^2$ has the same law as the process
denoted~$D_t$ in the statement of Theorem \ref{thm:difflimit}.
Convergence of finite-dimensional distributions for
$( n^{-1} W_{nt}^2 )_{t \in[0,1]}
= ( n^{-1} \tZ_{nt} )_{t \in[0,1]}$ follows. The final statement in
the theorem
follows from Lemma~\ref{lem:bessel}.
\end{pf*}

Next, consider the leaky urn model of Section \ref{S:leaky}.

\begin{pf*} {Proof of Theorem \ref{th:leaky}}
This is an immediate
consequence of the $\P(\kappa=1)=1$ cases of Theorems \ref
{thm:mixed} and \ref{thm:moments2}.
\end{pf*}

\begin{remark*}
There is a short proof of the first part of Theorem
\ref{th:leaky} due to the existence of a particular martingale.
Consider the process
$Q'_n$ defined by
$Q'_n = Q(X'_n, Y'_n)$, where
\[
Q(x,y)  := \bigl(x + \tfrac{1}{2} \sgn(y ) - \tfrac{1}{2} \mathbf{1}_{\{y =
0\}} \sgn(x ) \bigr)^2
+ \bigl(y - \tfrac{1}{2} \sgn(x ) - \tfrac{1}{2} \1_{ \{x = 0 \}} \sgn(y )
\bigr)^2 .
\]
It turns out that $Q'_n$ is a (nonnegative) martingale.
Thus, it converges a.s. as $n \to\infty$. But since
$Q(x,y) \to\infty$
as $\| (x,y)\| \to\infty$, we must have
that eventually $(X'_n,Y'_n)$ gets trapped
in the closed class $\mathcal{C}$. So $\P(\tau< \infty)=1$.
\end{remark*}

%s8.2 ###
\subsection{\texorpdfstring{Proofs of Theorems \protect\ref{thm:moments2} and \protect\ref{thm:area}}
{Proofs of Theorems 2.6 and 2.9}}

The proofs of Theorems~\ref{thm:moments2} and~\ref{thm:area}
that we give in this section both rely on the
good estimates we have for the embedded process $\tilde Z_k$ to
analyze the noisy urn $(\tilde X_n, \tilde Y_n)$.
The main additional ingredient is to
relate the two different time-scales. The first result
concerns the time to traverse a quadrant.

\begin{lemma}
\label{lem:time}
Let $k \in\Z_+$.
The distribution of $\tilde\nu_{k+1} - \tilde\nu_k$ given $\tZ_k =n$
coincides with that of $Z_{k+1} + Z_k$ given $Z_k =n$.
In addition,
%
%e35 ###
\begin{equation}
\label{noisytime}
\tilde\nu_{k+1} - \tilde\nu_k = | \tilde X_{\tilde\nu_{k+1}} |
+ | \tilde Y_{\tilde\nu_{k+1}} | + \tilde Z_k .
\end{equation}
Moreover,
%
%e36 ###
\begin{equation}
\label{timetail}
\P( \tilde\nu_{k+1} - \tilde\nu_k > 3n \midd\tilde Z_k = n ) = O (
\exp\{ - n^{1/2} \} ) .
\end{equation}
\end{lemma}

\begin{pf}
Without loss of generality, suppose that we are traversing the first
quadrant. Starting at time $\tilde\nu_k + 1$, Lemma \ref{noisylem}
implies that the time until hitting the next axis,
$\tilde\nu_{k+1} - \tilde\nu_k -1$, has the same distribution as
the time
taken for the original simple harmonic urn to hit the next axis, starting
from $(Z_k, 1)$.
In this time, the simple harmonic urn must
make $Z_k$ horizontal jumps and $Z_{k+1} -1$ vertical jumps. Thus,
$\tilde\nu_{k+1} - \tilde\nu_k -1$ has the same
distribution as $Z_{k+1} + Z_k -1$, conditional
on $Z_k = \tilde Z_k$. Thus, we obtain the first statement in the lemma.
For equation (\ref{noisytime}), note that between times $\tilde\nu_k
+1$ and
$\tilde\nu_{k+1}$ the noisy urn must make $\tZ_k$ horizontal steps and
(in this case) $|\tilde Y_{\tilde\nu_{k+1}} |-1$ vertical steps.
Finally, we have from the first statement of the lemma that
\[
\P( \tilde\nu_{k+1} - \tilde\nu_k > 3n \midd\tilde Z_k = n )
= \P( Z_{k+1} > 2n \midd Z_k = n ) ,
\]
and then (\ref{timetail}) follows from (\ref{deltatail}).
\end{pf}

Very roughly speaking, the key to our Theorems \ref{thm:moments2}
and \ref{thm:area} is the fact that
$\tau\approx\sum_{k=0}^{\tau_q} W^2_k$ and $A \approx\sum
_{k=0}^{\tau_q} W^4_k$.
Thus, to study $\tau$ and $A$ we need to look at sums of powers of
$W_k$ over a \textit{single excursion}. First, we will give results
for $S_\alpha:= \sum_{k=0}^{\tau_q} W^\alpha_k$, $\alpha\geq0$.
Then we quantify the
approximations
``$\approx$'' for $\tau$ and $A$ by a series of bounds.

Let $M := \max_{0 \leq k \leq\tau_q} W_k$ denote the maximum of
the first excursion of~$W_k$. For ease of notation, for the
rest of this section we set $r := 6 \E[ \kappa] - 3$.

\begin{lemma}
\label{max} Suppose that $r > -1$. Then for any $\eps>0$, for all
$x$ sufficiently large
\[
x^{-1-r-\eps} \leq\P( M \geq x ) \leq x^{-1-r+\eps} .
\]
In particular, for any $\eps>0$, $\E[ M^{1+r+\eps} ]= \infty$ but
$\E[ M^{1+r-\eps} ]< \infty$.
\end{lemma}

\begin{pf}
It follows from Lemma \ref{lem:sqrtincrements} and
some routine Taylor's theorem computations
that for any $\eps>0$ there exists $w_0 \in[1, \infty)$ such that for
any $x \geq w_0$,
\begin{eqnarray*}
\E[ W_{k+1}^{1+r+\eps} - W_k ^{1+r+\eps} \midd W_k = x] & \geq&0, \\
\E[ W_{k+1}^{1+r-\eps} - W_k ^{1+r-\eps} \midd W_k = x] & \leq&0.
\end{eqnarray*}
Let $\eta:= \min\{k \in\Z_+ \dvtx W_k \leq w_0 \}$ and $\sigma_x :=
\min\{ k \in\Z_+ \dvtx W_k \geq x \}$. Recall that $(W_k)_{k \in
\Z_+}$ is an irreducible time-homogeneous Markov chain on a
countable subset of $[1,\infty)$. It follows that to prove the
lemma it suffices to show that, for some $w \geq2w_0$, for any
$\eps>0$,
%
%e37 ###
\begin{equation}
\label{eq0}
x^{-1-r-\eps} \leq\P( \sigma_x < \eta\midd W_0 = w )
\leq x^{-1-r+\eps}
\end{equation}
for all $x$ large enough.

We first prove the lower bound in (\ref{eq0}). Fix $x > w$. We
have that $W_{k \wedge\eta\wedge\sigma_x}^{1+r+\eps}$\vspace*{2pt} is a
submartingale, and, since $W_k$ is an irreducible Markov chain,
$\eta< \infty$ and $\sigma_x < \infty$ a.s. Hence
\[
\P( \sigma_x < \eta) \E[ W_{\sigma_x}^{1+r+\eps} ]
+ \bigl(1 - \P( \sigma_x < \eta) \bigr) \E[ W_{\eta}^{1+r+\eps} ] \geq
w^{1+r+\eps} .
\]
Here $W_\eta\leq w_0$ a.s., and for some $C \in(0,\infty)$ and all
$x > w$,
\[
\E[ W_{\sigma_x}^{1+r+\eps} ] \leq
\E\bigl[ \bigl( x + ( W_{\sigma_x} - W_{\sigma_x -1} )
\bigr)^{1+r+\eps} \bigr] \leq C x^{1+r+\eps} ,
\]
since $\E[ (
W_{\sigma_x} - W_{\sigma_x -1} )^{1+r+\eps} ]$ is uniformly
bounded in $x$, by equation (\ref{wlem1}). It follows that
\[
\P( \sigma_x < \eta) [ C x^{1+r+\eps} - w_0^{1+r+\eps} ] \geq
w^{1+r+\eps} - w_0^{1+r+\eps} > 0 ,
\]
which yields the lower bound
in (\ref{eq0}). The upper bound follows by a similar argument
based on the supermartingale property of $W_{k \wedge\eta\wedge
\sigma_x}^{1+r-\eps}$.
\end{pf}

The next result gives the desired moment bounds for $S_\alpha$.

\begin{lemma}
\label{slem} Let $\alpha\geq0$ and $r > -1$. Then $\E[ S_\alpha^p] <
\infty$ if $p< \frac{1+r}{\alpha+2}$ and $\E[ S_\alpha^p] =
\infty$ if $p> \frac{1+r}{\alpha+2}$.
\end{lemma}

\begin{pf}
First we prove the upper bound. Clearly, $S_\alpha\leq
(1+\tau_q) M^\alpha$. Then, by H\"older's inequality,
\[
\E[ S_\alpha^p ] \leq\bigl( \E\bigl[ (1+\tau_q)^{(2+\alpha)p/2} \bigr] \bigr)^{\bfrac
{2}{2+\alpha}} \bigl( \E\bigl[ M^{(2+\alpha)p} \bigr] \bigr)^{\bfrac{\alpha}{2+\alpha
}} .
\]
For $p < \frac{1+r}{\alpha+2}$ we have $(2+\alpha)p/2 < (1+r)/2 =
3\E[\kappa] -1$
and $(2+\alpha)p < 1+r$ so that Lemma \ref{max} and Theorem
\ref{thm:moments} give the upper bound.

For the lower bound, we claim that there exists $C \in(0,\infty)$
such that
%
%e38 ###
\begin{equation}
\label{claim3}
\P( S_\alpha\geq x )
\geq\tfrac{1}{2} \P\bigl( M \geq C x^{\bfrac{1}{\alpha+2}} \bigr)
\end{equation}
for all $x$ large enough.
Given the claim (\ref{claim3}), we have, for any $\eps>0$,
\[
\E[ S_\alpha^p ] \geq\frac{p}{2} \int_1^\infty x^{p-1}\P\bigl( M \geq
C x^{\bfrac{1}{\alpha+2}}
\bigr)\, \ud x
\geq\frac{p}{2} \int_1^\infty x^{p-1} x^{-\afrac{1+r}{\alpha+2}
-\eps} \,\ud x
\]
by Lemma \ref{max}. Thus, $\E[ S_\alpha^p]=\infty$ for $p > \frac
{1+r}{\alpha+2}$.
It remains to verify (\ref{claim3}).
Fix $y > 2$. Let $\F_k=\sigma(W_1,\dots,W_k)$, and define stopping times
\[
\sigma_1 = \min\{ k \in\N\dvtx W_k \geq y \} ;\qquad
\sigma_2 = \min\{ k \geq\sigma_1 \dvtx W_k \leq y/2 \} .
\]
Then $\{ \sigma_1 < \tau_q \}$, that is, the event that
$W_k$ reaches $y$ before $1$, is $\F_{\sigma_1}$-measurable.
Now
\begin{equation}
\label{eq5}
\qquad \P( \{ \sigma_1 < \tau_q \} \cap\{ \sigma_2 \geq
\sigma_1 + \delta y^2 \} )
= \E[ \1 \{ \sigma_1 < \tau_q \} \P( \sigma_2 \geq\sigma_1 +
\delta y^2 \midd\F_{\sigma_1} ) ] .
\end{equation}
We claim that there exists $\delta>0$ so that
%
%e39 ###
\begin{equation}
\label{claim1a} \P( \sigma_2 \geq\sigma_1 + \delta y^2 \midd\F
_{\sigma_1} )
\geq\tfrac{1}{2} \qquad  \mbox{a.s.}
\end{equation}
Let $D_k = (y - W_k )^2 \1 \{ W_k < y\}$.
Then, with $\Delta_k = W_{k+1} - W_k$,
\[
\E[ D_{k+1} - D_k \midd\F_k ] \leq2(W_k-y) \E[ \Delta_k \midd\F_k
] + \E[ \Delta_k^2 \midd\F_t ] .
\]
Lemma \ref{lem:sqrtincrements} implies that on $\{W_k > y/2\}$ this
last display is bounded above by some
$C <\infty$ not depending on $y$.
Hence, an appropriate maximal inequality (\cite{mvw}, Lemma 3.1), implies
(since $D_{\sigma_1} =0$)
that
$\P( \max_{0 \leq s \leq k} D_{(\sigma_1+s) \wedge\sigma_2} \geq w
) \leq C k/w$.
Then, since $D_{\sigma_2} \geq y^2 /4$, we have
\begin{eqnarray*}
\P( \sigma_2 \leq\sigma_1 + \delta y^2 \midd\F_{\sigma_1} )
&\leq&\P\Bigl( \max_{1 \leq s \leq\delta y^2} D_{(\sigma_1 + s) \wedge
\sigma_2} \geq(y^2/4) \bigm|\F_{\sigma_1} \Bigr)\\
&\leq&\frac{C \delta y^2}{(y^2/4) } \leq\frac{1}{2}\qquad  \mbox{a.s.}
\end{eqnarray*}
for
$\delta>0$ small enough. Hence, (\ref{claim1a}) follows. Combining
(\ref{eq5}) and (\ref{claim1a}), we get
\[
\P( \{ \sigma_1 < \tau_q \} \cap\{ \sigma_2 \geq\sigma_1 +
\delta y^2 \} )
\geq\tfrac{1}{2} \P( \sigma_1 < \tau_q ) = \tfrac{1}{2} \P( M
\geq y ) .
\]
Moreover, on $\{ \sigma_1 < \tau_q \} \cap\{
\sigma_2 \geq\sigma_1 + \delta y^2 \}$ we have that $W_s \geq
y/2$ for all $\sigma_1 \leq s < \sigma_2$, of which there are at least
$\delta y^2$ values; hence $S_\alpha\geq\delta y^2 \times
(y/2)^\alpha$. Now taking
$x = 2^{-\alpha} \delta y^{2+\alpha}$, we obtain (\ref{claim3}), and
so complete the proof.
\end{pf}

Next, we need a technical lemma.

\begin{lemma}
\label{kappasums}
Let $p \geq0$. Then for any $\eps>0$ there exists $C<\infty$ such that
%
%e40 ###
\begin{equation} \label{claim2}
\E\Biggl[ \Biggl( \sum_{k=1}^{\tau_q} |\kappa_{\tilde\nu_{k}}| \Biggr)^p \Biggr]
\leq C \E[ \tau_q^{p+\eps} ] .
\end{equation}
\end{lemma}

\begin{pf}
For any $s \in(0,1)$,
\[
\P\Biggl( \sum_{k=1}^{\tau_q} |\kappa_{\tilde\nu_{k}}| > x \Biggr)
\leq\P( \tau_q > x^s ) + \P\Biggl( \sum_{k=1}^{x^s} |\kappa_{\tilde
\nu_{k}}| > x \Biggr) .
\]
For any random variable $X$,
\[
\E[ X^p] = p \int_0^\infty x^{p-1} \P(
X > x)\,\ud x
\leq1 + p \int_1^\infty x^{p-1}  \P( X >x ) \,\ud x;
\]
 so
%
%e41 ###
\begin{eqnarray}
\label{eq3}
\E\Biggl[ \Biggl( \sum_{k=1}^{\tau_q} |\kappa_{\tilde\nu_{k}}| \Biggr)^p \Biggr]
&\leq&1 + p \int_1^\infty x^{p-1} \P( \tau_q > x^s ) \,\ud x \nonumber\\[-8pt]\\[-8pt]
&&{}+ p
\int_1^\infty x^{p-1} \P\Biggl( \sum_{k=1}^{x^s} |\kappa_{\tilde\nu
_{k}}| > x \Biggr)\, \ud x .\nonumber
\end{eqnarray}
Here, we have that
\[
\P\Biggl( \sum_{k=1}^{x^s} |\kappa_{\tilde\nu_{k}}| > x \Biggr) \leq
\P\Biggl( \bigcup_{k=1}^{x^s} \{ |\kappa_{\tilde\nu_{k}}| > x^{1-s} \} \Biggr)
\leq\sum_{k=1}^{x^s} \P( |\kappa| > x^{1-s} )
\]
by Boole's inequality. Then Markov's inequality and the
moment bound (\ref{kappabound}) yield
%
%e42 ###
\begin{equation}
\label{eq4} \P\Biggl( \sum_{k=1}^{x^s} |\kappa_{\tilde\nu_{k}}| > x \Biggr)
\leq x^s \P\bigl( \re^{\lambda| \kappa|} > \re^{x^{1-s}} \bigr)
\leq x^s \E\bigl[ \re^{\lambda| \kappa| } \bigr] \re^{-x^{1-s}} .
\end{equation}
It follows that, since $s<1$, the final integral in (\ref{eq3})
is finite for any $p$. Also, from Markov's inequality,
for any $\eps>0$,
\[
\int_1^\infty x^{p-1} \P( \tau_q > x^s ) \,\ud x \leq\E[ \tau
_q^{p+\eps} ] \int_1^\infty
x^{p-1-s(p+\eps)} \,\ud x;
\]
taking $s$ close to $1$ this last integral is finite, and (\ref
{claim2}) follows
(noting $\tau_q \geq1$ by definition).
\end{pf}

\begin{pf*}{Proof of Theorem \ref{thm:moments2}}
By the definitions of $\tau$ and $\tau_q$, we have that
$\tau= \tilde\nu_{\tau_q} = -1+ \sum_{k=1}^{\tau_q} ( \tilde\nu
_k - \tilde\nu_{k-1} )$, recalling $\tilde\nu_0 = -1$.
Hence, by Lemma \ref{lem:time},
\[
\tau= -1 + \sum_{k=1}^{\tau_q} ( W_{k-1}^2 + W_{k}^2 ) + R
\]
for $R$ a random variable such that $| R | \leq\sum_{k=1}^{\tau_q}
|\kappa_{\tilde\nu_{k}}|$.
It follows that
%
%e43 ###
\begin{equation}
\label{bounds1} -1 + \sum_{k=0}^{\tau_q} W_k^2 - | R | \leq\tau
\leq2 \sum_{k=0}^{\tau_q} W_k^2 + |R| .
\end{equation}
Lemma \ref{kappasums} implies that for any $\eps>0$ there exists $C<
\infty$ such that
$\E[ |R|^p ] \leq C \E[ \tau_q^{p+\eps} ]$.
The $\E[ \kappa] > 1/3$ case of the theorem now follows from (\ref
{bounds1}) with
Theorem \ref{thm:moments}, Lemma \ref{slem} and Minkowski's
inequality. In the $\E[\kappa] =1/3$ case,
it is required to prove that $\E[\tau^p] = \infty$ for any $p>0$;
this follows from the
$\E[\kappa] = 1/3$ case of Theorem \ref{thm:moments} and the fact
that $\tau\geq\tau_q$ a.s.
\end{pf*}

\begin{pf*}{Proof of Theorem \ref{thm:area}}
First, note that we can write
\[
A = \sum_{n=1}^\tau T_n = \sum_{n=1}^{\tilde\nu_{\tau_q}} T_n =
\sum_{k=1}^{\tau_q} A_k ,
\]
where $A_1 = \sum_{n=1}^{\tilde\nu_1} T_n$ and
$A_k = \sum_{n= \tilde\nu_{k-1} +1}^{\tilde\nu_k} T_n$ ($k \geq2$)
is the area swept
out in traversing a quadrant for the $k$th time. Since $A_k \geq1/2$,
part (i) of the
theorem is immediate from part (i) of Theorem \ref{thm:mixed}. For
part (ii), we have
that
\begin{eqnarray*}
A_k & \leq&( \tilde Z_k + | \kappa_{\tilde\nu_k} | )
( \tilde Z_{k-1} + | \kappa_{\tilde\nu_{k-1}} | ) \\
& \leq& W_k^4 + W_{k-1}^4 + W_{k-1}^2 | \kappa_{\tilde\nu_{k}} | +
W_{k}^2 | \kappa_{\tilde\nu_{k-1}} |
+ | \kappa_{\tilde\nu_{k-1}} || \kappa_{\tilde\nu_{k}} | .
\end{eqnarray*}
Thus,
\[
A \leq2 \sum_{k=0}^{\tau_q} W_k^4 + R_1 + R_2 + R_3 ,
\]
where $R_1\!=\!\sum_{k=1}^{\tau_q}
W_{k-1}^2 | \kappa_{\tilde\nu_{k}} |$,
$R_2\!=\!\sum_{k=1}^{\tau_q}
W_{k}^2 | \kappa_{\tilde\nu_{k-1}} |$ and
$R_3\!=\!\sum_{k=1}^{\tau_q}
| \kappa_{\tilde\nu_{k-1}} || \kappa_{\tilde\nu_{k}} |$.
Here $\sum_{k=0}^{\tau_q} W_k^4$ has finite $p$th moment
for $p< \frac{3 \E[ \kappa] -1}{3}$, by Lemma \ref{slem}.
Next we deal with the terms $R_1, R_2$ and $R_3$.
Consider $R_1$. We have that, by H\"older's inequality, $ \E[ |R_1 |^p
]$ is at most
\begin{eqnarray*}
&&\E\Biggl[ \Biggl(\sum_{k=1}^{\tau_q}
W_{k-1}^2 \Biggr)^{3p/2} \Biggr]^{2/3} \E\Biggl[ \Biggl(
\sum_{k=0}^{\tau_q} | \kappa_{\tilde\nu_{k}} | \Biggr)^{3p}
\Biggr]^{1/3}\\
&&\qquad \leq C' \E\Biggl[ \Biggl(\sum_{k=0}^{\tau_q}
W_{k}^2 \Biggr)^{3p/2} \Biggr]^{2/3} \E[ \tau_q^{3p+\eps} ]^{1/3}
\end{eqnarray*}
for any $\eps>0$, by (\ref{claim2}). Lemma \ref{slem} and Theorem
\ref{thm:moments}
show that this is finite\vspace*{1pt} provided $p < \frac{3 \E[ \kappa] -1}{3}$
(taking $\eps$ small enough).
A similar argument holds for~$R_2$. Finally,
\[
\E[ |R_3|^p ] \leq\E\Biggl[ \Biggl( \sum_{k=0}^{\tau_q}
| \kappa_{\tilde\nu_{k}} | \Biggr)^{2p} \Biggr]
\leq C'' \E[ \tau_q^{2p+\eps} ]
\]
for any $\eps>0$, by (\ref{claim2}). For $\eps$ small enough, this
is also finite when $p < \frac{3 \E[ \kappa] -1}{3}$
by Theorem \ref{thm:moments}. These estimates and Minkowski's
inequality then complete
the proof.
\end{pf*}

%s8.3 ###
\subsection{\texorpdfstring{Proof of Theorem \protect\ref{thm:perc}}{Proofs of Theorem 2.11}}

We now turn to the percolation model
described in Section \ref{S: percolation}.

\begin{lemma}
\label{lem:perc1}
Let $v$ and $v'$ be any two vertices of $G$. Then with probability~1
there exists a vertex $w \in G$ such that the unique semi-infinite
oriented paths in $H$ from $v$ and $v'$ both pass through $w$.
\end{lemma}

\begin{pf}
Without loss of generality,
suppose $v, v'$ are distinct vertices in~$G$
on the positive $x$-axis on the same sheet of $\mathcal{R}$. Let $Z_0 =
|v| < Z'_0 = |v'|$.
The two paths in $H$ started at $v$ and $v'$, call them
$P$ and $P'$, respectively, lead to instances
of processes $Z_k$ and $Z'_k$, each a copy
of the simple harmonic urn embedded process $Z_k$.
Until $P$ and $P'$ meet,
the urn processes they instantiate are independent. Thus, it suffices
to take $Z_k$, $Z'_k$ to be independent and show that they
eventually cross with probability 1, so that the
underlying paths must meet. To do this, we consider the
process $(H_k)_{k \in\Z_+}$ defined by $H_k := \sqrt{Z'_k} - \sqrt
{Z_k}$ and show that it is eventually less than
or equal to $0$.

For convenience, we use the notation $W_k = (Z_k)^{1/2}$ and $W'_k =
(Z'_k)^{1/2}$. Since $H_{k+1} - H_k = (W'_{k+1}-W'_k) - (W_{k+1} -
W_k)$, we have that for $x<y$,
\[
\E[ H_{k+1} - H_k \midd W_k = x, W'_k = y] = \frac{1}{4y} - \frac
{1}{4x} + O( x^{-2} )
= -\frac{(y-x)}{4xy} + O(x^{-2} )
\]
by the $\E[\kappa] =0$ case of
(\ref{wlem2}). Similarly,
\[
\E[ (H_{k+1} - H_k)^2 \midd W_k = x, W'_k = y] = \tfrac{1}{3} + O(
x^{-1} ) ,
\]
from (\ref{wlem2}) again. Combining these, we see that
\begin{eqnarray*}
&&2 (y-x) \E[ H_{k+1} - H_k \midd W_k = x, W'_k = y] - \E[ (H_{k+1} -
H_k)^2 \midd W_k = x, W'_k = y] \\
%= \frac{2}{3} - \frac{x^2 +y^2}{2xy} + O(x^{-1} )
&&\qquad \leq- \tfrac{1}{3} + O (x^{-1} ) < 0
\end{eqnarray*}
for $x > C$, say.
However, we know from Theorem \ref{th:main1}
that $W_k$ is transient, so in particular
$W_k > C$ for all $k > T$ for some finite $T$.
Let
$\tau= \min\{ k \in\Z_+ : H_k \leq0 \}$. Then we have that
$H_{k} \1\{ k < \tau\}$, $k >T$, is a process on $\R_+$
satisfying Lamperti's recurrence criterion (cf. Proposition \ref
{prop:lampclass}). Here
$H_{k \wedge\tau}$ is not a Markov process but
the general form of Proposition \ref{prop:lampclass} applies (see~\cite{lamp1}, Theorem 3.2)
so we can conclude that
$\P( \tau< \infty) = 1$.
\end{pf}

\begin{lemma}
\label{lem:perc2}
The in-graph of any individual vertex in $H$ is almost surely finite.
\end{lemma}

\begin{pf}
We work in the dual percolation model $H'$.
As we have seen, the oriented paths through $H'$ simulate the leaky
simple harmonic urn via the mapping~$\Phi$. The path in $H'$ that
starts from a vertex over $(n+1/2,1/2)$
explores the outer boundary of the in-graph in $H$ of a lift of the set
$\{(i,0)\dvtx 1 \le i \le n\}$.
The leaky urn a.s. reaches the
steady state with one ball, so every oriented path in $H'$ a.s.
eventually joins the infinite path cycling immediately around the origin.
It follows that the in-graph of any vertex over a~co-ordinate axis is
a.s. finite. For any vertex $v$ of $H$, the oriented path
from $v$ a.s. contains a vertex $w$ over an axis, and the in-graph
of~$v$ is contained in the in-graph of $w$, so it too is a.s. finite.
\end{pf}

All that remains to complete the proof of Theorem \ref{thm:perc} is to
establish the two statements about the moments of $I(v)$. For $p < 2/3$,
$\E[I(v)^p]$ is bounded above by $\E[A^p]$, where $A$ is the area
swept out by a path of the leaky simple harmonic urn, or equivalently
by a path of the noisy simple harmonic urn with $\P(\kappa= 1) = 1$
up to the hitting time $\tau$.
$\E[A^p]$ is finite, by Theorem~\ref{thm:area}(ii).
The final claim $\E[I(v)] = \infty$ will be proven in the next
section as equation (\ref{in-graph}), using a connection with expected
exit times from quadrants.

%s9 ###
\section{Continuous-time models}
\label{sec:continuous}

%s9.1 ###
\subsection{\texorpdfstring{Expected traversal time: Proof of Theorem \protect\ref{th:fasttime}}
{Expected traversal time: Proof of Theorem 2.14}}
\label{S:fasttime}\mbox{}

\begin{pf*}{Proof of Lemma \ref{complexmartingale}}
A consequence of Dynkin's formula for a con\-tin\-u\-ous-time
Markov chain $X(t)$ on a countable state-space $S$ with infinitesimal
(generator) matrix $Q = (q_{ij})$
is that for a function $g\dvtx \R_+ \times S \to\R$ with continuous
time-derivative to be such that $g ( t, X(t) )$ is a local martingale,
it suffices that
%
%e44 ###
\begin{equation}
\label{dynkin}
\frac{ \partial g (t , x) }{\partial t} + Q ( g(t, \cdot) ) (t,x) =
0
\end{equation}
for all $x \in S$ and $t \in\R_+$: see, for example, \cite{robert}, page 364. In our case, $S = \mathbb{C} \setminus\{ 0 \}$,
$X(t) = A(t) + i B(t)$,
and for $z = x+ i y \in\mathbb{C}$,
\begin{eqnarray*}
Q ( f) (z ) &=& \sum_{w \in S, w \neq z} q_{z w}[ f(w) - f(z) ]\\
&=& | x| \bigl[ f \bigl( z + \sgn(x) i \bigr) - f(z) \bigr] + |y| \bigl[ f\bigl ( z - \sgn(y) \bigr) -
f(z) \bigr] .
\end{eqnarray*}
Taking $f(x+iy) = g(t, x+iy)$ to be first $x \cos t + y \sin t$
and second $y \cos t -x \sin t$, we verify the identity (\ref{dynkin})
in each case.
Thus, the real and imaginary parts of $M_t$ are local martingales, and hence
martingales since it is not hard to see that $\E| A(t) + B(t) | <
\infty$.
\end{pf*}

To prove Theorem \ref{th:fasttime}, we need the following bound on the
deviations of $\tau_f$ from $\pi/2$.

\begin{lemma}
\label{lem:fastdev}
Suppose $\eps_n>0$ and $\eps_n \to0$ as $n \to\infty$. Let $\phi
_n \in[0,\pi/2]$.
Then as $n \to\infty$,
\[
\P\biggl( \biggl| \tau_f - \frac{\pi}{2} + \phi_n \biggr|
\geq\eps_n \bigm| A(0) =n \cos\phi_n, B(0) = n \sin\phi_n \biggr) = O(
n^{-1} \eps_n^{-2})
\]
uniformly in $(\phi_n)$.
\end{lemma}

\begin{pf}
First note that $M_0 = n \re^{i \phi_n}$ and, by the martingale property,
\[
\E[ |M_t - M_0 |^2 ] = \E[ |M_t |^2 ] - |M_0|^2 =
\E[ A(t)^2 + B(t)^2 ] - n^2.
\]
We claim that for all $n \in\N$ and $t \in\R_+$,
%
%e45 ###
\begin{equation}
\label{claim1}
\E[ A(t)^2 + B(t)^2 ] - n^2 \leq\frac{t^2}{2} + 2^{1/2} n t .
\end{equation}
Since $M_t-M_0$ is a (complex) martingale,
$|M_t-M_0|^2$ is a submartingale. Doob's maximal
inequality therefore implies that, for any $r>0$,
\[
\P\Bigl( \sup_{0 \leq s \leq t} | M_s - M_0 | \geq r \Bigr) \leq r^{-2}
\E[ |M_t - M_0|^2 ] \leq2 t (t+n) r^{-2}
\]
by (\ref{claim1}).
Set $t_0 = (\pi/2) -\phi_n + \theta$ for $\theta\in(0,\pi/2)$.
Then on $\{ t_0 < \tau_f \}$, $A(t_0) + i B(t_0)$
has argument in $[ \phi_n , \pi/2]$, so that
$M_{t_0}$ has argument in $[ 2 \phi_n\,{-}\,(\pi/2)\,{-}\,\theta,\allowbreak \phi_n -
\theta]$.
All points with argument in the latter interval
are at distance at least $n \sin\theta$ from $M_0$.
Hence, on $\{ t_0 < \tau_f \}$,
\[
\sup_{0 \leq s \leq t_0} | M_s - M_0 |
\geq| M_{t_0} - M_0 | \geq n \sin\theta.
\]
It follows that for $\eps_n >0$ with $\eps_n \to0$,
\begin{eqnarray*}
\P\bigl( \tau_f > (\pi/2) -\phi_n +\eps_n \bigr) & \leq&\P\Bigl( \sup_{0 \leq
s \leq(\pi/2) -\phi_n +\eps_n } | M_s - M_0 | \geq n
\sin\eps_n \Bigr) \\
& =& O( n^{-1} (\sin\eps_n)^{-2} ) = O( n^{-1} \eps_n^{-2} ).
\end{eqnarray*}
A similar argument yields
the same bound for $\P_n ( \tau_f < (\pi/2) - \phi_n - \eps_n )$.
It remains to prove the claim (\ref{claim1}). First, note that
\begin{eqnarray*}
&&\E\bigl[ A(t+\Delta t)^2 + B(t +\Delta t)^2 - \bigl( A (t)^2 + B(t)^2 \bigr) \midd
A(t)=x, B(t)=y \bigr]\\
&&\qquad = (|x| + |y|) \Delta t + O ((\Delta t)^2 ),
\end{eqnarray*}
and $(|x| + |y|)^2 \leq2 (x^2 + y^2)$.
Writing $g(t) = \E[ A (t)^2 + B(t)^2 ]$, it follows that
\[
\frac{\ud}{\ud t} g(t) \leq\sqrt{2} g(t) ^{1/2}
\]
with $g(0) = n^2$. Hence, $g(t)^{1/2} \leq n + 2^{-1/2} t$. Squaring
both sides
yields (\ref{claim1}).
\end{pf}

A consequence of Lemma \ref{lem:fastdev} is that $\tau_f$
has finite moments of all orders, uniformly in the initial point.

\begin{lemma}
\label{taufmoms}
$\!\!\!\!\!$For any $r\!>\!0$, there exists $C\!<\!\infty$ such that \mbox{$\max_{n \in\N}
\E_n [ \tau_f^r ]\!\leq\!C$}.
\end{lemma}

\begin{pf}
By Lemma \ref{lem:fastdev}, we have that there exists $n_0 < \infty$
for which
%
%e46 ###
\begin{equation}
\label{taufeq}
\sup_{x >0, y>0 \dvtx | x+ i y | \geq n_0 }
\P\bigl( \tau_f - t > 2n_0 \midd A(t) + i B(t) = x + iy \bigr) \leq1/2
.
\end{equation}
On the other hand, if $| A(t) + iB(t) | < n_0$, we
have that $\tau_f - t$ is stochastically dominated by a sum of
$n_0$ exponential random variables with mean~$1$. Thus, by Markov's
inequality, the bound (\ref{taufeq}) holds for \textit{all} \mbox{$x>0,
y>0$}. Then, for $t > 1$, by conditioning on the path of the
process at times $2n_0, 4n_0, \ldots, 2n_0 (t-1)$ and using the
strong Markov property we have
\begin{eqnarray*}
&&\P_n ( \tau_f > 2n_0 t )\\
&&\qquad \leq\prod_{j=1}^{t-1} \sup_{x_j>0, y_j>0 } \P\bigl( \tau_f - 2n_0 j >
2n_0 \midd A(2n_0 j) + i B(2n_0 j) = x_j + i y_j \bigr)\\
&&\qquad \leq2^{1-t}
\end{eqnarray*}
by (\ref{taufeq}). Hence, $\P_n ( \tau_f > t )$ decays faster than
any power of
$t$, uniformly in~$n$.
\end{pf}

\begin{pf*}{Proof of Theorem \ref{th:fasttime}}
For now fix $n \in\N$.
Suppose $A(0) =Z_0 = n$, $B(0) = 0$. Note that $A(\tau_f) = 0$,
$B(\tau_f) = Z_1$.
The stopping time $\tau_f$ has all moments, by Lemma \ref{taufmoms},
while $\E_n [ | M_t |^2 ] = O(t^2)$ by (\ref{claim1}),
and $\E_n [ | M_{\tau_f} |^2 ] = \E_n [ Z_1^2 ] < \infty$.
It follows that the real and imaginary parts of the
martingale~$M_{t \wedge\tau_f}$ are
uniformly integrable. Hence,
we can apply the optional stopping theorem to
any linear combination of the real and imaginary parts of~$M_{t \wedge
\tau_f}$ to obtain
\[
\E_n [ Z_1 ( \alpha\sin\tau_f + \beta\cos\tau_f ) ] = \alpha n
\]
for any $\alpha, \beta\in\R$. Taking $\alpha= \cos\theta$,
$\beta= \sin\theta$ this says
\begin{eqnarray*}
n \cos\theta&=& \E_n [ Z_1 \sin(\theta+\tau_f) ] \\
&=& \E_n [ (Z_1 -
\E_n Z_1 ) \sin(\theta+\tau_f) ]
+ \E_n [ Z_1 ] \E_n [ \sin(\theta+ \tau_f) ]
\end{eqnarray*}
for any $\theta$.
By Cauchy--Schwarz,
the first term on the right-hand side here is bounded in absolute value
by $\sqrt{\Var_n(Z_1)}$, so on rearranging we have
\[
\biggl| \E_n [ \sin(\theta+\tau_f) ] - \frac{n\cos\theta}{\E_n [Z_1]}
\biggr| \le\frac{( \Var_n(Z_1) )^{1/2} }{\E_n [Z_1 ]}
\le\frac{( \E_n [ \Delta_1^2] )^{1/2} }{\E_n [Z_1 ]} ,
\]
and then using (\ref{zinc1}) and (\ref{zinc2}) we obtain, as $n \to
\infty$,
%
%e47 ###
\begin{equation}
\label{eqtheta}
| \E_n [ \sin(\theta+\tau_f) ] - \cos\theta| = O (n^{-1/2})
\end{equation}
uniformly in $\theta$. This strongly suggests that $\tau_f$ is concentrated
around $\pi/2,\allowbreak 5\pi/2, \ldots.$ To rule out the larger values, we need
to use Lemma \ref{lem:fastdev}. We proceed as follows.

Define the event $E_n := \{ | \tau_f - (\pi/2) | < \eps_n \}$ where
$\eps_n \to0$.
From the $\theta=-\pi/2$ case of (\ref{eqtheta}) we have that
$\E_n [\sin( \tau_f - (\pi/2)) ] = O(n^{-1/2})$.\vadjust{\goodbreak} Since $\sin x = x
+ O(x^3 )$
as $x \to0$ we have
\begin{eqnarray*}
&&\E_n \bigl[ \1_{E_n} \sin\bigl( \tau_f - (\pi/2)\bigr) \bigr]\\
&&\qquad = \E_n [ \tau_f \1_{E_n} ] - \frac{\pi}{2} + O(\eps_n^3) + O ( \P
_n (E_n^c) ) \\
&&\qquad = \E_n [ \tau_f ] - \frac{\pi}{2} + O(\eps_n^3)+ O \bigl( (\E_n [
\tau_f^r] )^{1/r} (\P_n (E_n^c))^{1-(1/r)} \bigr)
\end{eqnarray*}
for any $r>1$, by H\"older's inequality. Here $\E_n [ \tau_f^r ] =
O(1)$, by Lemma \ref{taufmoms},
so that for any $\delta>0$, choosing $r$ large enough we see that the
final term in the last display is $O(n^{\delta-1} \eps_n^{-2} )$ by Lemma
\ref{lem:fastdev}. Hence, for any $\delta>0$,
\[
O(n^{-1/2}) = \E_n [ \tau_f ] - \frac{\pi}{2} + O (n^{\delta-1}
\eps_n^{-2} ) + O(\eps_n^3)
+ \E_n \bigl[ \1_{E^c_n} \sin\bigl( \tau_f - (\pi/2)\bigr) \bigr] ,
\]
and this last expectation is $O( n^{-1} \eps_n^{-2} )$ by Lemma \ref
{lem:fastdev}
once more. Taking $\eps_n = n^{-1/4}$ yields (\ref{fast1}).
Next, from the $\theta=0$ case of (\ref{eqtheta}) we have that
$\E_n | 1- \cos( \tau_f - (\pi/2)) | = O(n^{-1/2})$.
This time
\[
\E_n \bigl[
\bigl| 1- \cos\bigl( \tau_f - (\pi/2)\bigr) \bigr| \1_{E_n} \bigr]
= \E_n [ | \tau_f - (\pi/2) |^2 \1_{E_n} ] + O( \eps_n^4) .
\]
Following a similar argument to that for (\ref{fast1}),
we obtain (\ref{fast2}).
\end{pf*}

%s9.2 ###
\subsection{\texorpdfstring{Traversal time and area enclosed: Proofs of Theorems
\protect\ref{thm:areamean} and \protect\ref{th:timearea}}
{Traversal time and area enclosed: Proofs of Theorems 2.10 and 2.15}}
\label{S:areatime}

Our proofs of Theorems \ref{thm:areamean} and \ref{th:timearea}
both use the percolation model of Section \ref{S: percolation}.

\begin{pf*}{Proof of Theorem \ref{th:timearea}}
The asymptotic statement in the theorem is a consequence of Theorem
\ref{th:fasttime}. Thus it remains to prove
the exact formula.
For $x > 0$, and $y \ge0$, let $T(x,y)$ denote $\E[ \tau_f \midd A(0)
= x, B(0) = y ]$.
Also, set $T(0,y) = 0$ for $y > 0$. Note that $T(n,0) = \E_n [ \tau_f
] $. Conditioning on the first step shows that for $x > 0$ and $y \ge0$,
\[
T(x,y) = \frac{1}{x+y} + \frac{x}{x+y} T(x,y+1) + \frac{y}{x+y}
T(x-1, y) .
\]
For fixed $x$, $T(x,y) \to0$ as $y \to\infty$. Indeed, for $y \ge1$
the time to make $x$ horizontal jumps is
stochastically dominated by the sum of $x$ exponential random variables
with mean $1/y$.

We now consider the percolation model restricted to the first quadrant.
More precisely, we consider the induced graph on the set of sites
$(x,y)$ with $x \ge0$ and $y > 0$, on a single sheet of $\mathcal{R}$.
Let $I(x,y)$ denote the expected number of sites in the in-graph of
$(x,y)$ in this restricted model. This count includes the site $(x,y)$ itself.
For $x > 0$, we also set $I(x,0) = 0$.
Considering the two possible directed edges into the site $(x,y)$, we obtain
\[
I(x,y) = 1 + \frac{y}{x+y+1} I(x+1, y) + \frac{x}{x+ y-1} I(x, y-1) .
\]
Dividing through by $(x+y)$, we have
\[
\frac{I(x,y)}{x+y} = \frac{1}{x+y} + \frac{y}{x+y} \biggl(\frac
{I(x+1,y)}{(x+1) + y} \biggr) + \frac{x}{x+y}\biggl (\frac{I(x,y-1)}{x+(y-1)} \biggr)
.\vadjust{\goodbreak}
\]
We now claim that for each fixed $y$, $I(x,y)$ is bounded as $x \to
\infty$. Indeed, the number of sites in the in-graph of $(x,y)$ is at
most $y$ plus $y$ times the number of horizontal edges in this
in-graph. The number of horizontal edges may be stochastically bounded
above by the sum of $y$ geometric random variables with mean $1/x$, so
its mean tends to $0$ as $x \to\infty$.

We see that $I(y,x)/(x+y)$ and $T(x,y)$ satisfy the same recurrence
relation with the same boundary conditions; their difference satisfies
a homogeneous recurrence relation with boundary condition $0$ at $x=0$
and limit $0$ as $y \to\infty$ for each fixed $x$. An induction with
respect to $x$ shows that the difference is identically zero.
In particular, taking $x = m$ and $y=0$, for any $m \ge1$, we find
\[
\ I(0,m) = m T(m,0) .
\]
The union of the in-graphs of the sites $(0,m)$, for $1 \le m \le n$,
is the set of all sites $(x,y)$ with $x \ge0$ and $y > 0$ that lie
under the oriented path of the dual percolation graph $H'$ that starts
at $(-1/2, n+1/2)$. Each of these sites lies at the center of a unit
square with vertices $(x \pm1/2, y \pm1/2)$, and the union of these
squares is the region bounded by the dual percolation path and the
lines $x = -1/2$ and $y = 1/2$. Reflecting this region in the line $y =
x+ 1/2$, we obtain a sample of the region bounded by a simple harmonic
urn path and the coordinate axes.
The expected number of unit squares in this region is therefore $\sum
_{m=1}^n I(0,m)$, so we are done.
\end{pf*}

\begin{pf*}{Proof of Theorem \ref{thm:areamean}}
The argument uses a similar idea to the proof of Theorem \ref
{th:timearea}, this time
for the percolation model on the whole of~$\mathcal{R}$. Choose a
continuous branch of the argument function on $\mathcal{R}$.
Let $I_+(v)$ denote the expected number of points $w$ with $\arg(w) >
0$ in the in-graph of~$v$ in $H$, including~$v$ itself if $\arg(v) >
0$. Arguing as before, if the projection of~$v$ to $\mathbb{Z}^2$ is
$(x,y)$, then $I_+(v)$ satisfies the boundary condition $I_+(v) = 0$
for $\arg(v) \le0$, and the recurrence relation
\begin{eqnarray*}
I_+(v) &=& 1 + \frac{|y|}{|x +\sgn(y) | + |y| }
I_+\bigl( v + (\sgn(y),0) \bigr)\\
&&{} + \frac{|x|}{|x| + |y -\sgn(x)| }
I_+\bigl( v + ( 0 , -\sgn(x))\bigr) ,
\end{eqnarray*}
where on the right-hand side $I_+$ is evaluated at two of the
neighbors of $v$ in the graph $G$.
Setting $J_+(v) := I_+(v) / ( |x(v)| + |y(v)| )$, we have a recurrence
relation for $J_+$:
\begin{eqnarray*}
J_+(v) &=& \frac{1}{|x| + |y|}
+ \frac{|y|}{|x|+ |y| }
J_+\bigl(v + (\sgn(y) , 0)\bigr) \\
&&{}+ \frac{|x|}{|x| + |y| }
J_+\bigl( v + ( 0 , -\sgn(x))\bigr) .
\end{eqnarray*}

The same recurrence relation and boundary conditions hold for
$T_+(\overline{v})$, where $T_+(w)$ is the expected time to hit the
set $\arg z \ge0$ in $\mathcal{R}$ in the fast embedding, starting
from a vertex $w$. Here, $\overline{v}$ is the vertex of $G$ at the
same distance from the origin as $v$, satisfying $\arg(\overline{v})
= -\arg(v)$. The reasoning of the previous proof shows that
$T_+(\overline{v}) = J_+(v)$ for all vertices $v$ with $\arg v \le\pi
/2$, and the argument may be repeated on the subsequent quadrants to
show by induction that $T_+(\overline{v}) = J_+(v)$ for all vertices $v$.
We therefore have the lower bound
\[
J_+(v) = T_+(\overline{v}) \geq\biggl\lfloor\frac{\arg(v)}{\pi/2}
\biggr\rfloor\inf_n \E_n [ \tau_f ] .
\]
The asymptotic expression (\ref{fast1}),
together with trivial lower bounds for small~$n$,
implies that $\inf_n \E_n [ \tau_f] > 0$. Therefore, as $v$ varies
over the set of vertices of~$G$ with a given projection $(x,y)$, both
$J_+(v)$ and $I_+(v)$ tend to infinity with $\arg(v)$.
Note that $I_+(v)$ is a lower bound for $I(v)$, and $I(v)$ depends only
on the projection $(x,y)$. It follows that
%
%e48 ###
\begin{equation}
\label{in-graph}
\E[I(v) ] = \infty.
\end{equation}

Recall that the oriented path of $H'$ starting at $(m+\frac{1}{2},
-\frac{1}{2})$ explores the outer boundary of the in-graph of the set
$S$ of vertices with $\arg(v) = 0$ and $ x \le m$, and that it can be
mapped via $\Phi$ onto a path of the leaky simple harmonic urn. Let $A$
denote the area swept out by this path up until time $\tau$ (the
hitting time of $\{(x,y) : |x| + |y| = 1\}$). The mapping $\Phi$ from
the vertices of $G'$ to $\mathbb{Z}^2$ can be extended by affine
interpolation to a locally area-preserving map from $\mathcal{R}$ to
$\mathbb{R}^2 \setminus(0,0)$. So $A$ is equal to the area swept out
by the dual percolation path until its projection hits the set $\{(\pm
\frac{1}{2}, \pm\frac{1}{2})\}$. Since the expected number of points
in the in-graph that it surrounds is infinite, we have $\E[A]
=\infty$.
\end{pf*}

%s9.3 ###
\subsection{Exact formulae for expected traversal time and enclosed area}
\label{S:polys}

In this section, we present some explicit, if mysterious, formulae for
the expected area enclosed
by a quadrant-traversal of the urn process and the expected
quadrant-traversal time in the fast embedding.
We obtain these formulae in a similar way to our first proof of Lemma
\ref{Lemma: transition probabilities},
and they are reminiscent, but more involved than, the formulae for the
Eulerian numbers. There is thus
some hope that the asymptotics of these formulae can be handled as in
the proof
of Lemma \ref{lem: asymptotic}, which gives a possible approach to the
resolution
of Conjecture \ref{conj:fasttime}.

\begin{lemma}
$\E_n[\mbox{Area enclosed} ]$ and $\E_n [ \tau_f]$ are
rational polynomials of degree $n$ evaluated at $\re$:
\begin{eqnarray}
\label{poly:area}
\E_n [\mbox{Area enclosed} ] & =& \sum_{i=1}^n \sum_{x = 1}^i
\frac{i^{n-x-i} i! (-1)^{n-i}}{(n-i)! (i-x)!} \Biggl(\re^i - \sum
_{k=0}^{i-1} \frac{i^k}{k!} \Biggr), \\
\label{poly:time}
\E_n [ \tau_f] & =& \sum_{i=1}^n
\sum_{x=1}^i \frac{ i^{n-x-i-1} i! (-1)^{n-i}}{(n-i)!(i-x)!} \Biggl(\re^i
- \sum_{k=0}^{i-1} \frac{i^k}{k!} \Biggr) .
\end{eqnarray}
\end{lemma}

\begin{pf}
The expected area enclosed can be obtained by summing the probabilities
that each unit square of the first quadrant is enclosed; that is,
\[
\E_n [ \mbox{Area enclosed} ] = \sum_{x = 1}^{n} \sum_{y =
1}^\infty\P_n\bigl((x,y) \mbox{ lies on or below the urn path}\bigr) .
\]
In terms of the slow continuous-time embedding of Section \ref
{secRubin}, $(x,y)$ lies on or below the urn path if and only if
$\sum_{j=1}^{y-1} j \zeta_j < \sum_{i=x}^n i \xi_i$.
Let
\[
R_{n,x,y} = \sum_{i=x}^n i \xi_i - \sum_{j=1}^{y-1} j \zeta_j ,
\]
so that
\[
\E_n [\mbox{Area enclosed} ] = \sum_{x=1}^n \sum_{y = 1}^\infty
\P(R_{n,x,y} > 0 ) .
\]
The moment generating function of $R_{n,x,y}$ is
\[
\E[\exp(\theta R_{n,x,y} ) ] = \prod_{i=x}^n \frac{1}{1 - i \theta
} \prod_{j=1}^{y-1} \frac{1}{1 + j \theta} = \sum_{i=x}^n \frac
{\alpha_i}{1 - i \theta} + \sum_{j=1}^{y-1} \frac{\beta_j}{1 + j
\theta} ,
\]
where
\[
\alpha_i = \frac{i^{n-x+y-1} i! (-1)^{(n-i)}}{(i+y-1)! (i-x)! (n-i)!} .
\]
Now the density of $R_{n,x,y}$ at $w > 0$ is
\[
\sum_{i=x}^n \alpha_i \frac{\exp(w/i)}{i} ,
\]
so that $\P(R_{n,x,y} > 0 ) = \sum_{i = x}^n \alpha_i$.
Therefore
\[
\E_n [\mbox{Area enclosed} ] = \sum_{x=1}^n \sum_{y = 1}^\infty
\sum_{i = x}^n \frac{i^{n-x+y-1} i! (-1)^{(n-i)}}{(i+y-1)! (i-x)!
(n-i)!} .
\]
The series converges absolutely so we can rearrange to obtain (\ref
{poly:area}).
By the first equality in Theorem \ref{th:timearea}, we find that
$\E_n[\tau_f]$ is also a rational polynomial of degree $n$ evaluated at
$\re$. After some simplification, we obtain~(\ref{poly:time}).
\end{pf}

A remarkable simplification occurs in the derivation of
(\ref{poly:time}) from (\ref{poly:area}), so it is natural to try
the same step again, obtaining
\[
\frac{1}{n}(\E_n [ \tau_f] - \E_{n-1} [ \tau_f ]) = \sum_{i=1}^n
\sum_{x=1}^i \frac{ i^{n-x-i-2} i! (-1)^{n-i}}{(n-i)!(i-x)!} \Biggl(\re^i
- \sum_{k=0}^{i-1} \frac{i^k}{k!} \Biggr) .
\]
In light of Theorem \ref{th:fasttime} and Conjecture \ref
{conj:fasttime}, we would like to prove that this expression decays
exponentially as $n \to\infty$.
Let us make one more observation that might be relevant
to Conjecture \ref{conj:fasttime}. Define $F(i) = \sum_{x=1}^i \frac
{i!}{(i-x)! i^x}$, which\vspace*{1pt}
can be interpreted as the expected number of distinct balls drawn
if we draw from an urn containing $i$ distinguishable balls, with
replacement, stopping when we first draw some ball for the second time.
We have already seen, in equation (\ref{poly:time}), that
\[
\E_n [\tau_f] = \sum_{i=1}^n F(i) \frac{(-1)^{n-i}
i^{n-i-1}}{(n-i)!} \Biggl(\re^i - \sum_{k=0}^{i-1}\frac{i^k}{k!} \Biggr) ;
\]
perhaps one could exploit the resemblance to the formula
\[
\E_n [ 1 / Z_{k+1} ] = \sum_{i=1}^n \frac{(-1)^{n-i}
i^{n-i-1}}{(n-i)!} \Biggl(\re^i - \sum_{k=0}^{i}\frac{i^k}{k!} \Biggr),
\]
but we were unable to do so.

%s10 ###
\section{Other stochastic models related to the simple harmonic urn}
\label{sec:misc}

%s10.1 ###
\subsection{A stationary model: The simple harmonic flea circus}
\label{S: stationary model}

In Section \ref{secRubin}, we saw that
the Markov chain $Z_k$ has an infinite invariant measure $\pi(n) = n$.
We can understand this in the probabilistic setting by considering the
formal sum of infinitely many independent copies of the
fast embedding.
Here is an informal description of the model.
At time $0$, populate each vertex of $\mathbb{Z}^2$ with
an independent Poisson-distributed number of fleas with mean $1$.
Each flea performs a copy of the
process $(A(t),B(t))$, independently of all the other fleas. Let
$N_t(m,n)$ denote the number of fleas at location $(m,n)$ at time $t$.

As we make no further use of this process in this paper, we do not
define it more formally. Instead we
just state the following result and sketch the proof: compare the lemma
in \cite{hoffrose}, Section 2, which the authors
attribute to Doob.

\begin{lemma}
The process $\{N_t(m,n)\dvtx m,n \in\mathbb{Z}\}$ is stationary. That is,
for each fixed
time $t > 0$, the array $N_t(m,n)$, $m,n \in\mathbb{Z}$ consists of
independent Poisson(1) random variables.
The process is reversible in the sense that the ensemble of random
variables $N_t(m,n), 0 \le t \le c$,
has the same law as the ensemble $N_{c-t}(m, -n)$, $0 \le t \le c$, for
any $c > 0$.
\end{lemma}

This skew-reversibility allows us to extend the stationary process to
all times $t \in\R$.

To see that the process has stationary means, note that the
expectations $\E[N_t(m,n)]$ satisfy a system of coupled differential equations:
\begin{eqnarray*}
\frac{\ud}{\ud t} \E[N_t(m,n)] &=& -(|m| + |n|)\E[N_t(m,n)] + |m|\E
\bigl[N_t\bigl(m,n-\sgn(m)\bigr)\bigr] \\
&&{}+ |n|\E\bigl[N_t\bigl(m+\sgn(n),n\bigr)\bigr],
\end{eqnarray*}
the solution to which is simply $\E[N_t(m,n)] = 1$ for all $t$, $m$
and $n$.

To establish the independence of the variables $N_t(a,b)$ when $t > 0$
is fixed, we use a Poisson thinning argument. That is, we construct
each variable
$N_0(m,n)$ as an infinite sum of independent Poisson random variables
$N(m,n,a,b)$ with means
\[
\E[N(m,n,a,b)] = \P\bigl( (A(t),B(t)) = (a,b) \midd(A(0),B(0)) = (m,n) \bigr) .
\]
The variable $N(m,n,a,b)$ gives the number of fleas that start at
$(m,n)$ at time $0$
and are at $(a,b)$ at time $t$. Then $N_t(a,b)$ is also a sum of
infinitely many independent
Poisson random variables, whose means sum to $1$, so it is a Poisson
random variable with mean $1$.
Moreover, for $(a,b) \neq(a',b')$, the corresponding sets of summands
are disjoint, so $N_t(a,b)$ and $N_t(a',b')$ are independent.

%s10.2 ###
\subsection{The Poisson earthquakes model}

We saw how the percolation model of Section \ref{S: percolation} gives
a static grand coupling of
many instances of (paths of) the simple harmonic urn. In this section,
we describe a model, based on ``earthquakes,''
that gives a dynamic grand coupling of many instances of simple
harmonic urn processes with particularly
interesting geometrical properties.

The earthquakes model is defined as a continuous-time Markov chain
taking values in the group of area-preserving homeomorphisms of the
plane, which we will write as
\[
\mathfrak{S}_t \dvtx \mathbb{R}^2 \to\mathbb{R}^2 ,\qquad  t \in\mathbb{R}.
\]
It will have the properties:
\begin{itemize}
\item$\mathfrak{S}_0$ is the identity,
\item$\mathfrak{S}_t(0,0) = (0,0)$,
\item$\mathfrak{S}_t$ acts on $\mathbb{Z}^2$ as a permutation,
\item$\mathfrak{S}_s \circ\mathfrak{S}_t^{-1}$ has the same
distribution as $\mathfrak{S}_{s-t}$, and
\item for each pair $(x_0, y_0) \neq(x_1, y_1) \in\mathbb{Z}^2$,
the displacement vector
\[
\mathfrak{S}_t (x_1,y_1) - \mathfrak{S}_t (x_0, y_0)
\]
has the distribution of the continuous-time
fast embedding of the simple harmonic urn, starting at $(x_1-x_0, y_1 - y_0)$.
\end{itemize}

In order to construct $\mathfrak{S}_t$, we associate a unit-rate
Poisson process to each horizontal
strip $ H_n: = \{(x,y) \in\mathbb{R}^2 \dvtx n < y < n + 1\}$, and to
each vertical strip $V_n: = \{(x,y) \in\mathbb{R}^2 \dvtx n < x < n + 1\}$
(where $n$ ranges over $\mathbb{Z}$). All these Poisson processes
should be independent.
Each Poisson process determines the sequence of times at which an \emph
{earthquake} occurs along the corresponding strip. An earthquake is a
homeomorphism of the plane that translates one of the complementary
half-planes of the given strip through a unit distance parallel to the
strip, fixes the other complementary half-plane, and shears the strip
in between them. The fixed half-plane is always the one containing the
origin, and the other half-plane always moves in the anticlockwise
direction relative to the origin.

Consider a point $(x_0,y_0) \in\mathbb{R}^2$. We wish to define
$\mathfrak{S}_t(x_0,y_0)$
for all \mbox{$t \ge0$}. We will define inductively a sequence of stopping
times $\eps_i$, and points $(x_i,\allowbreak y_i) \in\mathbb{R}^2$,
for $i \in\Z_+$. First, set $\eps_0 = 0$. For $i \in\N$, suppose
we have defined $(x_{i-1},y_{i-1})$
and $\eps_{i-1}$. Let $\eps_i$ be the least point greater than $\eps
_{i-1}$ in
the union of the Poisson processes associated to those strips for which
$(x_{i-1},y_{i-1})$ and $(0,0)$ do not both lie in one or other
complementary half-plane. This is a.s. well defined since there are
only finitely many such strips, and a.s. there
is only one strip for which an earthquake occurs at time $\eps_i$.
That earthquake moves $(x_{i-1},y_{i-1})$ to
$(x_i,y_i)$. Note that $\eps_i - \eps_{i-1}$ is an exponential random
variable with
mean $1/(\lceil|x_{i-1}|\rceil+ \lceil|y_{i-1}|\rceil)$,
conditionally independent of all previous jumps, given this mean. Since
each earthquake increases the distance between any two points by at
most $1$, it follows that a.s. the process does not explode in finite
time. That is, $\eps_i \to\infty$ as $i \to\infty$. Define
$\mathfrak{S}_t (x_0,y_0)$ to be $(x_i,y_i)$, where $\eps_i \le t <
\eps_{i+1}$. The construction of $\mathfrak{S}_t$ for $t < 0$ is
similar, using the inverses of the earthquakes.

Note that we cannot simply define $\mathfrak{S}_t$ for $t > 0$ to be
the composition of all the earthquakes that occur between times $0$ and
$t$, because almost surely infinitely many earthquakes occur during
this time; however, any bounded subset of the plane will only be
affected by finitely many of these, so the composition makes sense locally.

The properties listed above follow directly from the construction. For
$(x_0,y_0)$, $(x_1, y_1) \in\mathbb{Z}^2$, the displacement vector $
(\Delta x_t, \Delta y_t) = \mathfrak{S}_t (x_1,y_1) - \mathfrak{S}_t
(x_0, y_0)$ only changes when an earthquake occurs along a strip that
separates the two endpoints; the waiting time after $t$ for this to
occur is exponentially distributed with mean $1/(|\Delta x_t| + |\Delta
y_t| )$, and conditionally independent of $\mathfrak{S}_t$ given
$(\Delta x_t, \Delta y_t)$.

%f3 ###
\begin{figure}

\includegraphics{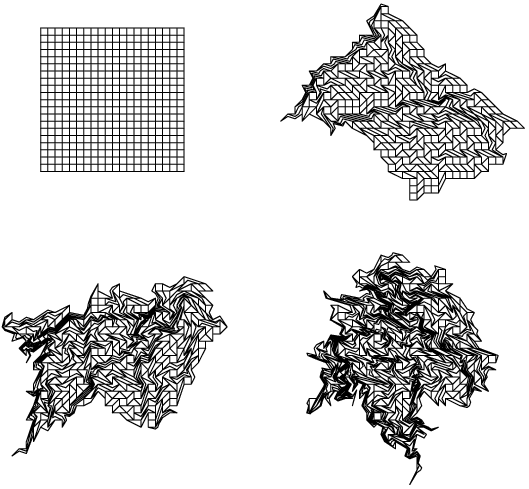}

\caption{A simulation of $\mathfrak{S}_t$, shown at times $t = 0$, $\pi
/4$, $\pi/2$, $3\pi/4$ acting on a $20 \times20$ box.}
\label{quakes}
\end{figure}

The model is spatially homogeneous in the following sense. Fix some
$(a,b) \in\mathbb{Z}^2$ and define
\[
\tilde{\mathfrak{S}}_t(x,y) = \mathfrak{S}_t (x+a,y+b) - \mathfrak
{S}_t(a,b) .
\]
Then $\tilde{\mathfrak{S}}_t$ has the same distribution as $\mathfrak{S}_t$.

\begin{lemma}
Define an oriented polygon $\Gamma$ by the cyclic sequence of vertices
\[
((x_1, y_1), \dots, (x_{n-1}, y_{n-1}), (x_n, y_n), (x_1,y_1)) ,\qquad  (x_i,
y_i) \in\mathbb{Z}^2 .
\]
The signed area enclosed by the polygon $\Gamma_t$, given by
\[
(\mathfrak{S}_t(x_1,y_1), \mathfrak{S}_t(x_2, y_2), \dots, \mathfrak
{S}_t(x_n, y_n), \mathfrak{S}_t(x_1, y_1) ) ,
\]
is a martingale.
\end{lemma}

\begin{pf}
For convenience, we write $(x_i(t), y_i(t)) = \mathfrak
{S}_t((x_i,y_i))$. The area enclosed by the oriented polygon $\Gamma
_t$ is
given by the integral $\frac{1}{2}\int_{\Gamma_t} x \,{d}y - y
\,{d}x$, which we can write as
\[
\frac{1}{2}\sum_{i=1}^n \bigl(x_i(t) y_{i+1}(t) - x_{i+1}(t) y_i(t) \bigr) ,
\]
where $(x_{n+1},y_{n+1})$ is taken to mean $(x_1, y_1)$. So it suffices
to show that each term in this sum is itself a martingale; let us
concentrate on the term $x_1(t) y_2(t) - x_2(t) y_1(t)$, considering
the first positive time at which either of $(x_1(t),y_1(t))$ or
$(x_2(t),y_2(t))$ jumps. There appear to be at least $36$ cases to
consider, depending on the ordering of $\{0, x_1, x_2\}$ and $\{0, y_1,
y_2\}$, but we can reduce this to four by taking advantage of the
spatial homogeneity of the earthquakes model, described above. By
choosing $(a,b)$ suitably, and replacing $\mathfrak{S}$ by
$\mathfrak{S}'$, we can assume that $x_i, y_i > 0$, for $i = 1, \dots,
n$. Furthermore, swapping the indices $1$ and $2$ only changes the sign
of $x_1(t) y_2(t) - x_2(t) y_1(t)$, so we may also assume that $x_1 \le
x_2$. Suppose that the first earthquake of interest is along a vertical
line. Then with probability $x_1/x_2$ it increments both $y_1$ and
$y_2$ and otherwise it increments only $y_2$. The expected jump in
$x_1(t) y_2(t) - x_2(t) y_1(t)$ conditional on the first relevant
earthquake being parallel to the $y$-axis is therefore
\begin{eqnarray*}
&&\frac{x_1}{x_2} \bigl(\bigl(x_1 (y_2+1) - x_2(y_1 +1) \bigr) - (x_1 y_2 - x_2 y_1)
\bigr)\\
&&\qquad {}+ \frac{x_2 - x_1}{x_2} \bigl( \bigl(x_1 (y_2 + 1) - x_2 y_1\bigr) - (x_1 y_2 - x_2
y_1) \bigr) = 0 .
\end{eqnarray*}
A similar argument shows that the expected jump in $x_1(t) y_2(t) -
x_2(t) y_1(t)$ conditional on the first relevant earthquake being
parallel to the $x$-axis is also zero.
\end{pf}

%s10.3 ###
\subsection{Random walks across the positive quadrant}

In this section, we describe another possible generalization of the simple
harmonic urn that has some independent interest.
We define a discrete-time process $(A_n,B_n)_{n \in\Z_+}$ on $\R^2$
based on the distribution of an underlying nonnegative,
nonarithmetic
random variable $X$
with $\E[ X ] =\mu\in(0,\infty)$ and $\Var[ X ] = \sigma^2 \in
(0,\infty)$.
Let $X_1, X_2, \ldots$ and $X'_1, X'_2, \ldots$ be independent copies
of $X$.
Roughly speaking, the walk starts on the horizontal axis and takes
jumps $(-X_i',X_i)$
until its first component is negative. At this point, suppose the walk
is at $(-r,s)$. Then the
walk starts again at $(s,0)$ and the process repeats.
We will see (Lemma \ref{quadshu}) that in the case when $X \sim U(0,1)$,
this process is closely related to the simple harmonic urn and is
consequently transient.
It is natural to study the same question for general distributions $X$.
It turns out that the
recurrence classification depends only on $\mu$ and $\sigma^2$. Our
proof uses
renewal theory.

We now formally define the model. With $X, X_n, X'_n$ as above,
we suppose that $\E[X^4] < \infty$.
Let
$(A_0,B_0)=(a,0)$, for $a>0$. Define the random process for $n \in\Z
_+$ by
\[
(A_{n+1},B_{n+1})=
\cases{
( A_n-X'_n, B_n+X_n ), &\quad if $A_n \geq0$,\cr
(B_n , 0), &\quad if $A_n < 0 $.
}
\]

\begin{thm}
\label{th:quadwalks}
Suppose $\E[X^4]<\infty$.
The walk $(A_n,B_n)$ is transient if and only if $\mu^2 > \sigma^2$.
\end{thm}

Set
$\tau_0 :=-1$ and for $k \in\N$,
\[
\tau_k :=\min\{ n > \tau_{k-1} \dvtx A_n<0\}.
\]
Define $T_k:= \tau_k - (\tau_{k-1} + 1)$. That is,
$T_k$ is the number of steps that the random walk takes to cross the
positive quadrant for the $k$th time.

\begin{lemma}
\label{quadshu}
If $X \sim U(0,1)$
and the initial value $a$ is distributed as
the sum of $n$ independent $U(0,1)$ random variables,
independent of the $X_i$ and~$X'_i$, then the distribution of the process
$(T_k)_{k \in\N}$ coincides with that of the embedded simple harmonic urn
process $(Z_k)_{k \in\N}$ conditional on $Z_0 = n$.
\end{lemma}

\begin{pf}
It suffices to show that $T_1 = \tau_1$ has the distribution of $Z_1$
conditional on $Z_0 = n$ and that conditional on $\tau_1$ the new
starting point $A_{1+\tau_1}$, which is $B_{\tau_1}$, has the
distribution of the sum of
$\tau_1$ independent $U(0,1)$ random variables.
Then the lemma will follow since the two processes $(\tau_k, B_{\tau
_k})$ and $(Z_k)$ are both Markov.
To achieve this, we couple the process $(A_n, B_n)$ up to time $\tau
_1$ with the renewal process described
in Section \ref{S: drift}. To begin,
identify $a$ with the sum $(1-\chi_1) + \cdots+ (1-\chi_n)$.
Then for $k \in\{ 1, \ldots, N(n)-n\}$, where $N(n) > n$
is as defined at (\ref{eq:renewproc}),
we identify $X'_{k}$
with $\chi_{n+k}$. For $m \le\tau_1$ we have
\[
A_m = a - \sum_{i=1}^m X'_i = n - \sum_{i=1}^{n+m} \chi_i ,
\]
so in particular we have $A_{N(n) -n-1} \geq0$
and $A_{N(n)-n} < 0$ by definition of~$N(n)$. Hence, $\tau_1 = N(n)-n$
has the distribution
of $Z_1$ by Lemma \ref{Lemma: same distribution}.
Moreover, $A_{1+\tau_1} = B_{\tau_1}$ is the sum of
the independent $U(0,1)$ random variables~$X_i$, $i = 1, \ldots, \tau_1$.
\end{pf}

Thus, by Theorem \ref{th:main1}, in the case where $X$ is
$U (0,1)$, the process
$(A_n,B_n)$ is transient, which is consistent with Theorem \ref{th:quadwalks}
since in the uniform case $\mu= 1/2$ and $\sigma^2 = 1/12$.
To study the general case, it is helpful to
rewrite the definition of $(A_n,B_n)$ explicitly
in the language of renewal theory.
Let $S_0 = S'_0 = 0$ and for $n \in\N$ set $S_n = \sum_{i=1}^n X_i$, $S'_n
= \sum_{i=1}^n X'_i$. Define the renewal counting function for $S'_n$
for $a>0$ as
\[
N(a) := \min\{ n \in\Z_+ \dvtx S'_n > a \} = 1 + \max\{ n \in\Z_+\dvtx
S'_n \leq a \}.
\]
Then starting at $(A_0,B_0) = (a,0)$, $a>0$, we see $\tau_1 =
N(a)$ so that $B_{\tau_1} = S_{N(a)}$. To study the recurrence and
transience of $(A_n,B_n)$, it thus suffices to study the process
$(R_n)_{n \in\Z_+}$ with $R_0 := a$ and $R_n$ having the
distribution of $S_{N(x)}$ given $R_{n-1} = x$. The increment of
the process $R_n$ starting from~$x$ thus is distributed as
$\Delta(x) := S_{N(x)} -x$. It is this random quantity that we
need to analyze.

\begin{lemma}
\label{quadjumps}
Suppose that $\E[ X^4] < \infty$. Then as $x \to\infty$, $\E[
|\Delta(x)|^4 ] = O(x^{2})$ and
\begin{eqnarray*}
\E[ \Delta(x) ] & = &\frac{\sigma^2 +\mu^2}{2 \mu}
+ O(x^{-1}), \\
\E[ \Delta(x)^2 ] & =& \frac{2 x \sigma^2}{\mu} + O(1) .
\end{eqnarray*}
\end{lemma}

\begin{pf}
We make use of results on higher-order renewal theory
expansions due to Smith \cite{smith2} (note that in \cite{smith2}
the renewal at $0$ is not counted). Conditioning on $N(x)$ and
using the independence of the $X_i$, $X'_i$, we obtain the Wald
equations:
\[
\E\bigl[ S_{N(x)} \bigr] = \mu\E[ N(x) ] ;\qquad  \Var\bigl[ S_{N(x)} \bigr] = \sigma^2 \E[
N(x) ] + \mu^2 \Var[ N(x) ] .
\]
Assuming $\E[X^3]<\infty$, \cite{smith2}, Theorem 1, shows that
\begin{eqnarray*}
\E[ N(x) ] & = &\frac{x}{\mu} + \frac{\sigma^2+\mu
^2}{2 \mu^2} + O(x^{-1}) ,\\
\Var[ N(x) ] & =& \frac{x \sigma^2}{\mu^3} + O(1).
\end{eqnarray*}
The expressions in the lemma for $\E[ \Delta(x)]$ and $\E[ \Delta
(x)^2]$ follow.

It remains to prove the bound for $\E[ |\Delta(x)|^4 ]$. Write
$\Delta(x)$ as
%
%e49 ###
\begin{equation}
\label{divideandconkout}
\qquad S_{N(x)} - x = \bigl( S_{N(x)} - \mu N(x) \bigr) + \bigl( \mu N(x) - \mu\E[ N(x) ] \bigr)
+ \bigl(\mu\E[ N(x) ] - x\bigr) .
\end{equation}
Assuming $\E[X^2]<\infty$, a result of Smith \cite{smith2}, Theorem 4,
implies that the
final bracket on the right-hand side of (\ref{divideandconkout})
is $O(1)$. For the first bracket
on the right-hand side of (\ref{divideandconkout}), it follows
from the Marcinkiewicz--Zygmund inequalities
(\cite{gut}, Corollary 8.2, page 151), that
\[
\E\bigl[ \bigl( S_{N(x)} - \mu N(x) \bigr)^4\bigr] \leq C \E[ N(x)^2 ] ,
\]
provided $\E[X^4] < \infty$. This last upper bound is $O(x^2)$ by the
computations in the first part of this proof. It remains to deal
with the second bracket on the right-hand side of (\ref{divideandconkout}).
By the algebra relating central moments to cumulants, we have
\[
\E\bigl[ \bigl( \mu N(x) - \mu\E[ N(x) ] \bigr)^4 \bigr]
= \mu^4 \bigl( k_4(x) + 3 k_2(x)^2 \bigr) ,
\]
where $k_r(x)$ denotes the
$r$th cumulant of $N(x)$. Again appealing to a result of Smith
(\cite{smith2}, Corollary 2, page 19), we have that $k_2(x)$ and $k_4(x)$
are both~$O(x)$ assuming $\E[X^4]<\infty$. (The fact
that \cite{smith2} does not count the renewal at $0$ is
unimportant here, since the $r$th cumulant of $N(x) \pm1$ differs
from $k_r(x)$ by a constant depending only on $r$.) Putting these
bounds together, we obtain from (\ref{divideandconkout}) and
Minkowski's inequality that $\E[ (S_{N(x)} - x)^4 ] = O(x^2)$.
\end{pf}

To prove Theorem \ref{th:quadwalks}, we basically need to compare $\E
[ \Delta(x)]$ to
$\E[ \Delta(x)^2]$.
As in our analysis of $\tZ_k$, it is most convenient to work on the square-root
scale. Set $V_n := R_n^{1/2}$.

\begin{lemma}
\label{vmoms}
Suppose that $\E[ X^4] < \infty$. Then there exists $\delta>0$ such
that as $y \to\infty$,
\begin{eqnarray*}
\E[ V_{n+1} - V_n \midd V_n = y ] & = &\frac{\E[ \Delta
(y^2) ]}{2y} - \frac{\E[\Delta(y^2)^2 ]}{8y^3}
+ O(y^{-1 -\delta}) ,\\
\E[ ( V_{n+1} - V_n)^2 \midd V_n = y ] & =& \frac{\E[\Delta(y^2)^2
]}{4y^2} + O (y^{-\delta} ),
\\
\E[ | V_{n+1} - V_n |^{3} \midd V_n = y ] & =& O ( 1 ) .
\end{eqnarray*}
\end{lemma}

\begin{pf}
The proof is similar to the proof of Lemma \ref{lem:sqrtincrements},
except that here we must work a little harder as we have weaker tail bounds
on $\Delta(x)$. Even so, the calculations will be familiar, so we
do not give all the details. Write $\E_x [ \cdot]$
for $\E[ \cdot\midd R_n = x ]$ and similarly for $\P_x$. From
Markov's inequality
and the fourth moment bound in Lemma \ref{quadjumps}, we have for $\eps
\in(0,1)$ that
%
%e50 ###
\begin{equation}
\label{bigjumptail}
\P_x \bigl( | \Delta(x) | > x^{1-\eps} \bigr) = O ( x^{4 \eps- 2} ) .
\end{equation}
We have that for $x \geq0$,
\[
\E[ V_{n+1} - V_n \midd V_n = x^{1/2} ] = \E_x [ R_{n+1}^{1/2} -
R_n^{1/2} ] = \E_x \bigl[ \bigl(x + \Delta(x)\bigr)^{1/2} - x^{1/2} \bigr] .
\]
Here we can write
%
%e51 ###
\begin{eqnarray}
\label{bigjump}
\bigl(x + \Delta(x)\bigr)^{1/2} - x^{1/2} &=& \bigl[ \bigl(x + \Delta(x)\bigr)^{1/2} - x^{1/2} \bigr]
\1 \{ | \Delta(x) | \leq x^{1-\eps} \}\nonumber\\[-8pt]\\[-8pt]
&&{}+ R_1 + R_2 \nonumber
\end{eqnarray}
for remainder terms $R_1$, $R_2$ that we define shortly.
The main term on the right-hand side admits a Taylor expansion
and analysis
(whose details we omit) in a similar
manner to the proof of Lemma \ref{lem:sqrtincrements}, and
contributes to the main terms in the statement of the present lemma.
The remainder terms in~(\ref{bigjump}) are
\begin{eqnarray*}
R_1 &=& \bigl[ \bigl(x + \Delta(x)\bigr)^{1/2} - x^{1/2} \bigr] \1 \{ \Delta(x) > x^{1-\eps
} \} ,\\
R_2 &=& \bigl[ \bigl(x + \Delta(x)\bigr)^{1/2} - x^{1/2} \bigr] \1 \{ \Delta(x) < -
x^{1-\eps} \} .
\end{eqnarray*}
For the second of these, we have $| R_2 | \leq x^{1/2} \1 \{ \Delta
(x) < - x^{1-\eps} \}$,
from which we obtain, for $r <4$, $\E_x [ |R_2 |^r ] = O( x^{4\eps+(r-4)/2})$,
by (\ref{bigjumptail}). Taking $\eps$ small enough,
this term contributes only to the negligible
terms in our final expressions. For $R_1$, we have the bound
\[
| R_1 | \leq C \bigl( 1 + | \Delta(x) | \bigr)^{(1/2) + \eps} \1 \{ \Delta(x)
> x^{1-\eps} \}
\]
for some $C \in(0,\infty)$ not depending on $x$, again for $\eps$
small enough. An application of H\"older's inequality and the
bound (\ref{bigjumptail}) implies that, for $r<4$, for any
$\eps>0$,
\begin{eqnarray*}
\E_x [ | R_1 |^r] &\leq& C \bigl( \E_x \bigl[ \bigl( 1+ | \Delta(x) |\bigr)^{4} \bigr] \bigr)^{\gfrac
{r(1+ 2\eps)}{8}}
\bigl( \P_x \bigl( \Delta(x) > x^{1-\eps} \bigr) \bigr)^{1-\gfrac{r(1+2\eps)}{8}}\\
& = & O\bigl(x^{6\eps+(r-4)/2 } \bigr) .
\end{eqnarray*}
It is now routine to complete the proof.
\end{pf}

\begin{pf*}{Proof of Theorem \ref{th:quadwalks}}
For the recurrence classification,
the crucial quantity is
\begin{eqnarray*}
&&2 y \E[ V_{n+1} - V_n \midd V_n = y ] - \E[ ( V_{n+1} - V_n)^2 \midd
V_n = y ] \\
&&\qquad = \E[ \Delta(y^2) ] - \frac{\E
[ \Delta(y^2)^2]}{2 y^2} + O(y^{-\delta} )
\end{eqnarray*}
by Lemma \ref{vmoms}. Now by Lemma \ref{quadjumps}, this last
expression is seen to be equal to
\[
\frac{\mu^2 - \sigma^2}{2 \mu} + O (y^{-\delta} ).
\]
Now \cite{lamp1}, Theorem 3.2, completes the proof.
\end{pf*}

\begin{remarks*}(i) To have some examples, note that if $X$ is
exponential, the
process is recurrent, while
if $X$ is the sum of two independent exponentials, it is transient. We
saw that if $X$
is $U(0,1)$ the process is transient; if $X$ is the square-root of a
$U(0,1)$ random variable,
it is recurrent.

(ii)
Another special case of the model that has some interesting features is
the case where $X$
is exponential with mean $1$.
In this particular case, a calculation shows that the distribution of
$T_{k+1}$ given $T_k=m$ is negative binomial
$(m+1,1/2)$, that is,
\[
\P(T_{k+1}=j \midd T_k=m) = \pmatrix{j+m\cr m}2^{-m-j-1}\qquad  ( j \in\Z_+).
\]
Since $\mu^2 =\sigma^2$, this case is in some sense
critical, a fact supported by the following branching process interpretation.
\end{remarks*}

Consider a version of the gambler's ruin problem. The gambler begins
with an initial stake, a pile of $m_0$ chips. A sequence of
independent
tosses of a fair coin is made; when the coin comes\vadjust{\goodbreak} up heads, a chip is
removed from the gambler's pile, but when it comes up tails, a chip is
added to a second pile by the casino.
The game ends when the gambler's original pile of chips is exhausted;
at this point the gambler receives the second pile of chips as his prize.
The total number of chips in play is a martingale;
by the optional stopping theorem, the expectation of the prize equals
the initial stake.
As a loss leader, the casino announces that it will add one extra chip
to each gambler's initial stake, so that the game is now in favour of
the gambler.
Suppose a gambler decides to play this game repeatedly, each time
investing his prize as the initial stake of the next game.
If the casino were to allow a zero stake (which of course it does not),
then the sequence of augmented stakes would form an irreducible Markov
chain $S_k$ on $\N$.
Conditional on $S_k = m$, the distribution of $S_{k+1}$ is negative
binomial $(m+1/2, 1/2)$.
So by the above results, this chain is recurrent. It follows that with
probability one the gambler will eventually lose everything.

We can interpret the sequence of prizes as a Galton--Watson process in
which each generation corresponds to one game, and individuals in the
population correspond to chips
in the gambler's pile at the start of the game. Each individual has a
Geo($1/2$) number of offspring
(i.e., the distribution that puts mass $2^{-1-k}$ on each $k \in\Z
_+$), being the chips that are added to the prize pile while that individual
is on top of the gambler's pile, and at each generation there is
additionally a Geo($1/2$) immigration, corresponding to the chips added
to the prize pile while the casino's bonus chip is
on top of the gambler's pile. This is a critical case of the
Galton--Watson process with immigration. By a result of Zubkov \cite{Zub}, if we start at time $0$ with population $0$, the time $\tau$ of
the next visit to $0$ has pgf
\[
\E[ s^\tau] = \frac{1}{s} + \frac{1}{\log(1-s)} .
\]
Since this tends to $1$ as $ s \nearrow1$, we have $\mathbb{P}(\tau<
\infty) = 1$. In fact, this can be deduced in an elementary way as follows.
The pgf of the Geo($1/2$) distribution is $f(s) = 1/(2-s)$, and its
$n$th iterate, the pgf of the $n$th
generation starting from one individual, is $f(s) =
(n-(n-1)s)/((n+1)-ns)$. In particular, the probability that an
individual has no descendants at the $n$th generation is
$n/(n+1)$. If $S_0 = 1$, then $S_k = 1$ if and only if for each $j=0,
\dots,k-1$ the bonus chip from game $j$ has no descendants at the
$(k-j)$th generation. These events are independent, so
\[
\P(S_k = 1 \midd S_0 = 1) = \prod_{j=0}^{k-1} \frac{k-j}{k-j-1} =
\frac{1}{k+1} ,
\]
which sums to $\infty$ over $k \in\N$ so that
the Markov chain is recurrent (see, e.g.,~\cite{Asmussen}, Proposition 1.2, Section I).

The results of Pakes \cite{Pakes} on the critical Galton--Watson
process with immigration show that the casino should certainly not
add two bonus chips to each stake, for then the process becomes
transient, and gambler's ruin will no longer apply.

\section*{Acknowledgments}

The authors are grateful to John Harris, Sir John Kingman, Iain MacPhee and
James Norris for helpful discussions, and to an anonymous referee for
constructive comments.
Most of this work was
carried out while Andrew Wade was at the Heilbronn Institute for
Mathematical Research, Department of Mathematics, University of
Bristol.

% imsref loaded by dianan, 2010-12-01 16:50:32
%

\printaddresses


\begin{thebibliography}{44}

%b1 ###
\bibitem{Asmussen}
%
\begin{bbook}[vtex]
\bauthor{\bsnm{Asmussen}, \bfnm{S{\o}ren}\binits{S.}}
(\byear{2003}).
\btitle{Applied Probability and Queues},
\bedition{2nd} ed.
\bseries{Applications of Mathematics (Stochastic Modelling and Applied Probability)}
\bvolume{51}.
\bpublisher{Springer}, \baddress{New York}.
\bid{mr={1978607}}
\end{bbook}
%
\endbibitem

%b2 ###
\bibitem{ai}
%
\begin{barticle}[mr]
\bauthor{\bsnm{Aspandiiarov}, \bfnm{S.}\binits{S.}} \AND
\bauthor{\bsnm{Iasnogorodski}, \bfnm{R.}\binits{R.}}
(\byear{1997}).
\btitle{Tails of passage-times and an application to stochastic
processes with
boundary reflection in wedges}.
\bjournal{Stochastic Process. Appl.}
\bvolume{66}
\bpages{115--145}.
\bid{doi={10.1016/S0304-4149(96)00118-4}, mr={1431874}}
\end{barticle}
%
\endbibitem

%b3 ###
\bibitem{aim}
%
\begin{barticle}[mr]
\bauthor{\bsnm{Aspandiiarov}, \bfnm{S.}\binits{S.}},
\bauthor{\bsnm{Iasnogorodski}, \bfnm{R.}\binits{R.}} \AND
\bauthor{\bsnm{Menshikov}, \bfnm{M.}\binits{M.}}
(\byear{1996}).
\btitle{Passage-time moments for nonnegative stochastic processes and an
application to reflected random walks in a quadrant}.
\bjournal{Ann. Probab.}
\bvolume{24}
\bpages{932--960}.
\bid{doi={10.1214/aop/1039639371}, mr={1404537}}
\end{barticle}
%
\endbibitem

%b4 ###
\bibitem{AK}
%
\begin{barticle}[mr]
\bauthor{\bsnm{Athreya}, \bfnm{Krishna B.}\binits{K. B.}} \AND
\bauthor{\bsnm{Karlin}, \bfnm{Samuel}\binits{S.}}
(\byear{1968}).
\btitle{Embedding of urn schemes into continuous time {M}arkov branching
processes and related limit theorems}.
\bjournal{Ann. Math. Statist.}
\bvolume{39}
\bpages{1801--1817}.
\bid{mr={0232455}}
\end{barticle}
%
\endbibitem

%b5 ###
\bibitem{AN}
%
\begin{bbook}[vtex]
\bauthor{\bsnm{Athreya}, \bfnm{Krishna B.}\binits{K. B.}} \AND
\bauthor{\bsnm{Ney}, \bfnm{Peter E.}\binits{P. E.}}
(\byear{1972}).
\btitle{Branching Processes}.
\bseries{Die Grundlehren der mathematischen Wissenschaften}
\bvolume{196}.
\bpublisher{Springer}, \baddress{New York}.
\bid{mr={0373040}}
\end{bbook}
%
\endbibitem

%b6 ###
\bibitem{Bona}
%
\begin{bbook}[vtex]
\bauthor{\bsnm{B{\'o}na}, \bfnm{Mikl{\'o}s}\binits{M.}}
(\byear{2004}).
\btitle{Combinatorics of Permutations}.
\bpublisher{Chapman \& Hall/CRC}, \baddress{Boca Raton, FL}.
\bid{doi={10.1201/9780203494370}, mr={2078910}}
\end{bbook}
%
\endbibitem

%b7 ###
\bibitem{Churchill}
%
\begin{barticle}[mr]
\bauthor{\bsnm{Churchill}, \bfnm{R. V.}\binits{R. V.}}
(\byear{1937}).
\btitle{The inversion of the {L}aplace transformation by a direct
expansion in
series and its application to boundary-value problems}.
\bjournal{Math. Z.}
\bvolume{42}
\bpages{567--579}.
\bid{doi={10.1007/BF01160095}, mr={1545692}}
\end{barticle}
%
\endbibitem

%b8 ###
\bibitem{Cox}
%
\begin{bbook}[vtex]
\bauthor{\bsnm{Cox}, \bfnm{D. R.}\binits{D. R.}}
(\byear{1962}).
\btitle{Renewal Theory}.
\bpublisher{Methuen}, \baddress{London}.
\bid{mr={0153061}}
\end{bbook}
%
\endbibitem

%b9 ###
\bibitem{Feller3}
%
\begin{barticle}[vtex]
\bauthor{\bsnm{Feller}, \bfnm{Willy}\binits{W.}}
(\byear{1941}).
\btitle{On the integral equation of renewal theory}.
\bjournal{Ann. Math. Statist.}
\bvolume{12}
\bpages{243--267}.
\bid{mr={0005419}}
\end{barticle}
%
\endbibitem

%b10 ###
\bibitem{Feller1}
%
\begin{bbook}[vtex]
\bauthor{\bsnm{Feller}, \bfnm{William}\binits{W.}}
(\byear{1968}).
\btitle{An Introduction to Probability Theory and Its Applications. {V}ol.~{I}},
\bedition{3rd} ed.
\bpublisher{Wiley}, \baddress{New York}.
\bid{mr={0228020}}
\end{bbook}
%
\endbibitem

%b11 ###
\bibitem{Feller2}
%
\begin{bbook}[vtex]
\bauthor{\bsnm{Feller}, \bfnm{William}\binits{W.}}
(\byear{1971}).
\btitle{An Introduction to Probability Theory and Its Applications. {V}ol.~{II}},
\bedition{2nd}~ed.
\bpublisher{Wiley}, \baddress{New York}.
\bid{mr={0270403}}
\end{bbook}
%
\endbibitem

%b12 ###
\bibitem{fgg}
%
\begin{bincollection}[mr]
\bauthor{\bsnm{Fisch}, \bfnm{Robert}\binits{R.}},
\bauthor{\bsnm{Gravner}, \bfnm{Janko}\binits{J.}} \AND
\bauthor{\bsnm{Griffeath}, \bfnm{David}\binits{D.}}
(\byear{1991}).
\btitle{Cyclic cellular automata in two dimensions}.
In \bbooktitle{Spatial Stochastic Processes}.
\bseries{Progress in Probability}
\bvolume{19}
\bpages{171--185}.
\bpublisher{Birkh\"auser}, \baddress{Boston, MA}.
\bid{mr={1144096}}
\end{bincollection}
%
\endbibitem

%b13 ###
\bibitem{Fl}
%
\begin{barticle}[mr]
\bauthor{\bsnm{Flajolet}, \bfnm{Philippe}\binits{P.}},
\bauthor{\bsnm{Gabarr{\'o}}, \bfnm{Joaquim}\binits{J.}} \AND
\bauthor{\bsnm{Pekari}, \bfnm{Helmut}\binits{H.}}
(\byear{2005}).
\btitle{Analytic urns}.
\bjournal{Ann. Probab.}
\bvolume{33}
\bpages{1200--1233}.
\bid{doi={10.1214/009117905000000026}, mr={2135318}}
\end{barticle}
%
\endbibitem

%b14 ###
\bibitem{gut}
%
\begin{bbook}[mr]
\bauthor{\bsnm{Gut}, \bfnm{Allan}\binits{A.}}
(\byear{2005}).
\btitle{Probability: A Graduate Course}.
\bpublisher{Springer}, \baddress{New York}.
\bid{mr={2125120}}
\end{bbook}
%
\endbibitem

%b15 ###
\bibitem{harris}
%
\begin{barticle}[mr]
\bauthor{\bsnm{Harris}, \bfnm{T. E.}\binits{T. E.}}
(\byear{1952}).
\btitle{First passage and recurrence distributions}.
\bjournal{Trans. Amer. Math. Soc.}
\bvolume{73}
\bpages{471--486}.
\bid{mr={0052057}}
\end{barticle}
%
\endbibitem

%b16 ###
\bibitem{hoffrose}
%
\begin{barticle}[mr]
\bauthor{\bsnm{Hoffman}, \bfnm{John R.}\binits{J. R.}} \AND
\bauthor{\bsnm{Rosenthal}, \bfnm{Jeffrey S.}\binits{J. S.}}
(\byear{1995}).
\btitle{Convergence of independent particle systems}.
\bjournal{Stochastic Process. Appl.}
\bvolume{56}
\bpages{295--305}.
\bid{doi={10.1016/0304-4149(94)00075-5}, mr={1325224}}
\end{barticle}
%
\endbibitem

%b17 ###
\bibitem{hutton}
%
\begin{barticle}[mr]
\bauthor{\bsnm{Hutton}, \bfnm{J.}\binits{J.}}
(\byear{1980}).
\btitle{The recurrence and transience of two-dimensional linear birth
and death
processes}.
\bjournal{Adv. in Appl. Probab.}
\bvolume{12}
\bpages{615--639}.
\bid{doi={10.2307/1426423}, mr={0578840}}
\end{barticle}
%
\endbibitem

%b18 ###
\bibitem{Ja}
%
\begin{barticle}[mr]
\bauthor{\bsnm{Janson}, \bfnm{Svante}\binits{S.}}
(\byear{2004}).
\btitle{Functional limit theorems for multitype branching processes and
generalized {P}\'olya urns}.
\bjournal{Stochastic Process. Appl.}
\bvolume{110}
\bpages{177--245}.
\bid{doi={10.1016/j.spa.2003.12.002}, mr={2040966}}
\end{barticle}
%
\endbibitem

%b19 ###
\bibitem{jensen}
%
\begin{barticle}[mr]
\bauthor{\bsnm{Jensen}, \bfnm{U.}\binits{U.}}
(\byear{1984}).
\btitle{Some remarks on the renewal function of the uniform distribution}.
\bjournal{Adv. in Appl. Probab.}
\bvolume{16}
\bpages{214--215}.
\bid{doi={10.2307/1427232}, mr={0732138}}
\end{barticle}
%
\endbibitem

%b20 ###
\bibitem{JK}
%
\begin{bbook}[vtex]
\bauthor{\bsnm{Johnson}, \bfnm{Norman L.}\binits{N. L.}} \AND
\bauthor{\bsnm{Kotz}, \bfnm{Samuel}\binits{S.}}
(\byear{1977}).
\btitle{Urn Models and Their Application: An Approach to Modern Discrete Probability Theory}.
\bpublisher{Wiley}, \baddress{New York}.
\bid{mr={0488211}}
\end{bbook}
%
\endbibitem

%b21 ###
\bibitem{kt2}
%
\begin{bbook}[mr]
\bauthor{\bsnm{Karlin}, \bfnm{Samuel}\binits{S.}} \AND
\bauthor{\bsnm{Taylor}, \bfnm{Howard M.}\binits{H. M.}}
(\byear{1981}).
\btitle{A Second Course in Stochastic Processes}.
\bpublisher{Academic Press},
\baddress{New York}.
\bid{mr={0611513}}
\end{bbook}
%
\endbibitem

%b22 ###
\bibitem{kesten}
%
\begin{barticle}[vtex]
\bauthor{\bsnm{Kesten}, \bfnm{Harry}\binits{H.}}
(\byear{1976}).
\btitle{Recurrence criteria for multi-dimensional {M}arkov chains and
multi-dimensional linear birth and death processes}.
\bjournal{Adv. in Appl. Probab.}
\bvolume{8}
\bpages{58--87}.
\bid{mr={0426187}}%
\end{barticle}
%
\endbibitem%

%b23 ###
\bibitem{K99}
%
\begin{barticle}[mr]
\bauthor{\bsnm{Kingman}, \bfnm{J. F. C.}\binits{J. F. C.}}
(\byear{1999}).
\btitle{Martingales in the {OK} {C}orral}.
\bjournal{Bull. London Math. Soc.}
\bvolume{31}
\bpages{601--606}.
\bid{doi={10.1112/S0024609399006098}, mr={1703841}}
\end{barticle}
%
\endbibitem

%b24 ###
\bibitem{KV}
%
\begin{barticle}[mr]
\bauthor{\bsnm{Kingman}, \bfnm{J. F. C.}\binits{J. F. C.}} \AND
\bauthor{\bsnm{Volkov}, \bfnm{S. E.}\binits{S. E.}}
(\byear{2003}).
\btitle{Solution to the {OK} {C}orral model via decoupling of {F}riedman's
urn}.
\bjournal{J. Theoret. Probab.}
\bvolume{16}
\bpages{267--276}.
\bid{doi={10.1023/A:1022294908268}, mr={1956831}}
\end{barticle}
%
\endbibitem

%b25 ###
\bibitem{BK}
%
\begin{bincollection}[vtex]
\bauthor{\bsnm{Kotz}, \bfnm{Samuel}\binits{S.}} \AND
\bauthor{\bsnm{Balakrishnan}, \bfnm{N.}\binits{N.}}
(\byear{1997}).
\btitle{Advances in urn models during the past two decades}.
In \bbooktitle{Advances in Combinatorial Methods and Applications to
Probability and Statistics}
\bpages{203--257}.
\bpublisher{Birkh\"auser}, \baddress{Boston, MA}.
\bid{mr={1456736}}
\end{bincollection}
%
\endbibitem

%b26 ###
\bibitem{lamp1}
%
\begin{barticle}[mr]
\bauthor{\bsnm{Lamperti}, \bfnm{John}\binits{J.}}
(\byear{1960}).
\btitle{Criteria for the recurrence or transience of stochastic
process. {I}.}
\bjournal{J. Math. Anal. Appl.}
\bvolume{1}
\bpages{314--330}.
\bid{mr={0126872}}
\end{barticle}
%
\endbibitem

%b27 ###
\bibitem{lamp2}
%
\begin{barticle}[mr]
\bauthor{\bsnm{Lamperti}, \bfnm{John}\binits{J.}}
(\byear{1962}).
\btitle{A new class of probability limit theorems}.
\bjournal{J. Math. Mech.}
\bvolume{11}
\bpages{749--772}.
\bid{mr={0148120}}
\end{barticle}
%
\endbibitem

%b28 ###
\bibitem{lamp3}
%
\begin{barticle}[mr]
\bauthor{\bsnm{Lamperti}, \bfnm{John}\binits{J.}}
(\byear{1963}).
\btitle{Criteria for stochastic processes. {II}. {P}assage-time moments}.
\bjournal{J.~Math. Anal. Appl.}
\bvolume{7}
\bpages{127--145}.
\bid{mr={0159361}}
\end{barticle}
%
\endbibitem

%b29 ###
\bibitem{mm}
%
\begin{barticle}[mr]
\bauthor{\bsnm{MacPhee}, \bfnm{I. M.}\binits{I. M.}} \AND
\bauthor{\bsnm{Menshikov}, \bfnm{M. V.}\binits{M. V.}}
(\byear{2003}).
\btitle{Critical random walks on two-dimensional complexes with
applications to
polling systems}.
\bjournal{Ann. Appl. Probab.}
\bvolume{13}
\bpages{1399--1422}.
\bid{doi={10.1214/aoap/1069786503}, mr={2023881}}
\end{barticle}
%
\endbibitem

%b30 ###
\bibitem{Mahmoud}
%
\begin{bbook}[mr]
\bauthor{\bsnm{Mahmoud}, \bfnm{Hosam M.}\binits{H. M.}}
(\byear{2009}).
\btitle{P\'olya Urn Models}.
\bpublisher{CRC Press}, \baddress{Boca Raton, FL}.
\bid{mr={2435823}}
\end{bbook}
%
\endbibitem

%b31 ###
\bibitem{mai}
%
\begin{barticle}[mr]
\bauthor{\bsnm{Menshikov}, \bfnm{M. V.}\binits{M. V.}},
\bauthor{\bsnm{{\`E}{\u\i}symont}, \bfnm{I. M.}\binits{I. M.}}
\AND
\bauthor{\bsnm{Yasnogorodski{\u\i}}, \bfnm{R.}\binits{R.}}
(\byear{1995}).
\btitle{Markov processes with asymptotically zero drift}.
\bjournal{Problemy Peredachi Informatsii}
\bvolume{31}
\bpages{60--75};
\bnote{translated in \textit{Probl. Inf. Transm.} \textbf{31} 248--261.}
\bid{mr={1367920}}
\end{barticle}
%
\endbibitem

%b32 ###
\bibitem{mvw}
%
\begin{barticle}[mr]
\bauthor{\bsnm{Menshikov}, \bfnm{M. V.}\binits{M. V.}},
\bauthor{\bsnm{Vachkovskaia}, \bfnm{M.}\binits{M.}} \AND
\bauthor{\bsnm{Wade}, \bfnm{A. R.}\binits{A. R.}}
(\byear{2008}).
\btitle{Asymptotic behaviour of randomly reflecting billiards in unbounded
tubular domains}.
\bjournal{J. Stat. Phys.}
\bvolume{132}
\bpages{1097--1133}.
\bid{doi={10.1007/s10955-008-9578-z}, mr={2430776}}
\end{barticle}
%
\endbibitem

%b33 ###
\bibitem{Pakes}
%
\begin{barticle}[mr]
\bauthor{\bsnm{Pakes}, \bfnm{A. G.}\binits{A. G.}}
(\byear{1971}).
\btitle{On the critical {G}alton--{W}atson process with immigration}.
\bjournal{J.~Austral. Math. Soc.}
\bvolume{12}
\bpages{476--482}.
\bid{mr={0307370}}
\end{barticle}
%
\endbibitem

%b34 ###
\bibitem{Review}
%
\begin{barticle}[mr]
\bauthor{\bsnm{Pemantle}, \bfnm{Robin}\binits{R.}}
(\byear{2007}).
\btitle{A survey of random processes with reinforcement}.
\bjournal{Probab. Surv.}
\bvolume{4}
\bpages{1--79 (electronic)}.
\bid{doi={10.1214/07-PS094}, mr={2282181}}
\end{barticle}
%
\endbibitem

%b35 ###
\bibitem{PV1999}
%
\begin{barticle}[mr]
\bauthor{\bsnm{Pemantle}, \bfnm{Robin}\binits{R.}} \AND
\bauthor{\bsnm{Volkov}, \bfnm{Stanislav}\binits{S.}}
(\byear{1999}).
\btitle{Vertex-reinforced random walk on {${\bf Z}$} has finite range}.
\bjournal{Ann. Probab.}
\bvolume{27}
\bpages{1368--1388}.
\bid{doi={10.1214/aop/1022677452}, mr={1733153}}
\end{barticle}
%
\endbibitem

%b36 ###
\bibitem{pittel}
%
\begin{barticle}[mr]
\bauthor{\bsnm{Pittel}, \bfnm{B.}\binits{B.}}
(\byear{1987}).
\btitle{An urn model for cannibal behavior}.
\bjournal{J. Appl. Probab.}
\bvolume{24}
\bpages{522--526}.
\bid{mr={0889816}}
\end{barticle}
%
\endbibitem

%b37 ###
\bibitem{robert}
%
\begin{bbook}[vtex]
\bauthor{\bsnm{Robert}, \bfnm{Philippe}\binits{P.}}
(\byear{2003}).
\btitle{Stochastic Networks and Queues},
\bedition{French} ed.
\bseries{Applications of Mathematics (Stochastic Modelling and Applied Probability)}
\bvolume{52}.
\bpublisher{Springer}, \baddress{Berlin}.
\bid{mr={1996883}}
\end{bbook}
%
\endbibitem

%b38 ###
\bibitem{Smith}
%
\begin{barticle}[vtex]
\bauthor{\bsnm{Smith}, \bfnm{Walter L.}\binits{W. L.}}
(\byear{1954}).
\btitle{Asymptotic renewal theorems}.
\bjournal{Proc. Roy. Soc. Edinburgh Sect. A}
\bvolume{64}
\bpages{9--48}.
\bid{mr={0060755}}
\end{barticle}
%
\endbibitem

%b39 ###
\bibitem{smith2}
%
\begin{barticle}[mr]
\bauthor{\bsnm{Smith}, \bfnm{Walter L.}\binits{W. L.}}
(\byear{1959}).
\btitle{On the cumulants of renewal processes}.
\bjournal{Biometrika}
\bvolume{46}
\bpages{1--29}.
\bid{mr={0104300}}
\end{barticle}
%
\endbibitem

%b40 ###
\bibitem{Stone}
%
\begin{barticle}[mr]
\bauthor{\bsnm{Stone}, \bfnm{Charles}\binits{C.}}
(\byear{1965}).
\btitle{On moment generating functions and renewal theory}.
\bjournal{Ann. Math. Statist.}
\bvolume{36}
\bpages{1298--1301}.
\bid{mr={0179857}}
\end{barticle}
%
\endbibitem

%b41 ###
\bibitem{sparac}
%
\begin{barticle}[vtex]
\bauthor{\bsnm{Ugrin-{\v{S}}parac}, \bfnm{G.}\binits{G.}}
(\byear{1990}).
\btitle{On a distribution encountered in the renewal process based on uniform
distribution}.
\bjournal{Glas. Mat. Ser. III}
\bvolume{25}
\bpages{221--233}.
\bid{mr={1108966}}
\end{barticle}
%
\endbibitem

%b42 ###
\bibitem{watson}
%
\begin{barticle}[mr]
\bauthor{\bsnm{Watson}, \bfnm{R. K.}\binits{R. K.}}
(\byear{1976}).
\btitle{An application of martingale methods to conflict models}.
\bjournal{Operations Res.}
\bvolume{24}
\bpages{380--382}.
\bid{mr={0418198}}
\end{barticle}
%
\endbibitem

%b43 ###
\bibitem{WM}
%
\begin{barticle}[mr]
\bauthor{\bsnm{Williams}, \bfnm{David}\binits{D.}} \AND
\bauthor{\bsnm{McIlroy}, \bfnm{Paul}\binits{P.}}
(\byear{1998}).
\btitle{The {OK} {C}orral and the power of the law (a curious {P}oisson-kernel
formula for a parabolic equation)}.
\bjournal{Bull. London Math. Soc.}
\bvolume{30}
\bpages{166--170}.
\bid{doi={10.1112/S0024609397004062}, mr={1489328}}
\end{barticle}
%
\endbibitem

%b44 ###
\bibitem{Zub}
%
\begin{barticle}[mr]
\bauthor{\bsnm{Zubkov}, \bfnm{A. M.}\binits{A. M.}}
(\byear{1972}).
\btitle{Life-periods of a branching process with immigration}.
\bjournal{Teor. Verojatnost. i Primenen.}
\bvolume{17}
\bpages{179--188};
\bnote{translated in \textit{Theor. Probab. Appl.}
\textbf{17}
174--183}.
\bid{mr={0300351}}
\end{barticle}
%
\endbibitem

\end{thebibliography}
\end{document}